\documentclass[10pt,twoside]{article}

\usepackage{amsmath}
\usepackage{amsfonts}
\usepackage{amssymb}
\usepackage{amscd}
\usepackage{amsthm}
\usepackage{amsbsy}
\usepackage{graphicx}
\usepackage{bm}

\def\@oddhead{\hfill \shorttitle \hfill \thepage}
\def\@evenhead{\thepage \hfill \shortauthor \hfill}
\def\@oddfoot{}
\def\@evenfoot{}

\date{}
\title{\ \\[0.4cm] \ \\ \bf  Spectral rigidity of group actions on  homogeneous spaces }
\author{Bachir Bekka\footnote{ IRMAR, UMR-CNRS 6625,  Universit\'e de  Rennes 1, 
Campus Beaulieu, F-35042  Rennes Cedex, 
 France. E-mail: bachir.bekka@univ-rennes1.fr. Support in part 
  by the  ANR (French Agence Nationale de la Recherche)
through the projects Labex Lebesgue (ANR-11-LABX-0020-01) and   GAMME  (ANR-14-CE25-0004).}}
\begin{document}
%-------------------

\maketitle

%------------------------

\thispagestyle{empty}

%--------------------------------------
\begin{abstract}
\vskip 3mm\footnotesize{Actions of a locally compact group $G$ on  a measure space $X$ 
give rise to unitary representations of $G$ on Hilbert spaces.
We review results on  the rigidity of these actions from the spectral point of view,  that is, 
results about  the existence of a spectral gap for associated averaging operators 
and their consequences.
 We will deal  both with  spaces $X$ with an infinite measure 
 as well as with spaces with an invariant probability measure.
 The  spectral gap property has several striking applications
 to group theory, geometry, ergodic theory, operator algebras,
 graph theory, theoretical computer science, etc.

\vskip 4.5mm
\noindent

\vspace*{2mm} \noindent{\bf 2000 Mathematics Subject Classification:
} 22D40,  37A30, 28D05, 43A07 

\vspace*{2mm} \noindent{\bf Keywords and Phrases: } Ergodic group actions, spectral gap property,  Kazhdan property, co-amenable actions} 
\end{abstract}
% ----------------------------------------------------------------

\newtheorem{theorem}{Theorem}[section]
\newtheorem{lemma}[theorem]{Lemma}
\newtheorem{corollary}[theorem]{Corollary}
\newtheorem{proposition}[theorem]{Proposition}
\newtheorem{obs}[theorem]{Observation}
\newenvironment{observation}{\begin{obs}\rm}{\end{obs}}
 \newtheorem{defi}[theorem]{Definition}
\newenvironment{definition}{\begin{defi}\rm}{\end{defi}}
\newtheorem{exa}[theorem]{Example}
\newenvironment{example}{\begin{exa}\rm}{\end{exa}}
\newtheorem{rem}[theorem]{Remark}
\newenvironment{remark}{\begin{rem}\rm}{\end{rem}}
\newtheorem{rems}[theorem]{Remarks}
\newenvironment{remarks}{\begin{rems}\rm}{\end{rems}}
\newtheorem{exer}[theorem]{Exercise}
\newenvironment{exercise}{\begin{exer}\rm}{\end{exer}}
\newtheorem{ack}[theorem]{Acknowlegment}
\newenvironment{acknowlegment}{\begin{ack}\rm}{\end{ack}}

\def\proof{\noindent\textbf{Proof}\quad}
\def\bsq{\blacksquare\medskip}
\def\n{\noindent}

\def\Z{\mathbf{Z}}
\def\N{\mathbf{N}}
\def\P{\mathbf{P}}
\def\Q{\mathbf{Q}}
\def\ZZ{\mathbf{Z}}
\def\NN{\mathbf{N}}
\def\PP{\mathbf{P}}
\def\QQ{\mathbf{Q}}
\def\CC{\mathbf{C}}
\def\PP{\mathbf{P}}
\def\Un{{\mathbf{1}}}
\def\proof{\noindent\textbf{Proof}\quad}
\def\bsq{\blacksquare\medskip}
\def\n{\noindent}
%%%%%%%%%%%%%%%%DEFINITIONS SYMBOLES MATHEMATIQUES%%%%%%%%%%%%%%%%%%%%%%%%%%%%%%%%%%%%%%%%%%%%%
\def\H{\mathcal H}
\def\L{\mathcal L}
\def\K{\mathcal K}
\def\P{\mathcal P}
\def\M{\mathcal M}
\def\F{\mathcal F}
\def\A{\mathcal A}
\def\B{\mathcal B}
\def\C{\mathcal C}
\def\V{\mathcal V}
\def\X{\mathcal X}
\def\G{\mathcal G}
\def\U{\mathcal U}
\def\R{\mathcal R}
\def\CalS{\mathcal S}
\def\Hpi{(\pi, {\mathcal H})}
\def\unirep{unitary representation}
\def\supp{{\rm supp}}
\def\NN{{\mathbf N}}
\def\ZZ{{\mathbf Z}}
\def\CCC{{\mathbf C}}
\def\RRR{{\mathbf R}}
\def\QQ{\mathbf Q}
\def\FF{\mathbf F}
\def\RR+{{\mathbf R}^*}
\def\TT{\mathbf T}
\def\KK{\mathbf K}
\def\GG{\mathbb G}
\def\EE{\mathbf E}
\def\HS{\rm HS}
\def\Q_p{{\mathbf Q}_p}
\def\S1{{\mathbf S}^1}

\def\ext{\rm ext}
\def\hG{\widehat G}
\def\Proj{\rm Proj}
\def\Ksigma{(\sigma, {\mathcal K})}
\def\ind{{\rm Ind}}
\def\sgn{\rm sgn}
\def\Tr{{\rm Trace}}
\def\PT{Property (T)}
\def\PTB{Property (${\rm T}_B$)}
\def\eps{\varepsilon}
\def\Ga{\Gamma}
\def\la{\lambda}
\def\ga{\gamma}
\def\La{\Lambda}
\def\vfi{\varphi}
\def\Sp{Sp(2n,\KK)}
\def\SL{SL(2,\KK)}
\def\chk{{\rm char} (\KK)}
\def\tous{\qquad\text{for all}\quad}
\def\bs{\backslash}
\def\nil{\La\backslash H}
\def\tor{\La Z\backslash H}
\def\ind{{\mathrm Ind}}
\def\Im{{\mathrm Im}}
\def\Ker{{\mathrm Ker}}
\def\Ad{{\mathrm Ad}}
\def\hH{\widehat H}
\def\rhonil{\rho_{\nil}}
\def\norm#1{\left\Vert #1 \right\Vert}
\def\abs#1{\left\vert #1 \right\vert}
\def\tout{\qquad\text{for all}\quad}
\def\Heis{H_{2n+1}(\RRR)}
\def\Proj{\rm Proj}
\def\Ksigma{(\sigma, {\mathcal K})}
\def\Ind{{\mathrm Ind}}
\def\sgn{\mathrm sgn}
\def\vol{\mathrm vol}
\def\Tr{{\mathrm Trace}}
\def\PT{Property (T)}
\def\PTB{Property (${\mathrm T}_B$)}
\def\eps{\varepsilon}
\def\Ga{\Gamma}
\def\ga{\gamma}
\def\la{\lambda}
\def\vfi{\varphi}
\def\Sp{Sp(2n,\KK)}
\def\SL{SL(2,\KK)}
\def\chk{{\rm char} (\KK)}
\def\tous{\qquad\text{for all}\quad}
\def\ExAut{{\mathrm Aut}_{\rm e}}
\def\Aut{{\mathrm Aut}}
\def\Inn{{\mathrm Inn}}
\def\Out{{\mathrm Out}}
\def\Ad{{\mathrm Ad}}
\def\VN{{\mathrm VN}}
\def\Hom{{\mathrm Hom}}
\def\spec{{\mathrm spec}}
\def\End{{\mathrm End}}
\def\H{\mathcal H}
\def\O{\mathcal O}
\def\L{\mathcal L}
\def\K{\mathcal K}
\def\P{\mathcal P}
\def\N{\mathcal N}
\def\F{\mathcal F}
\def\A{\mathcal A}
\def\B{\mathcal B}
\def\C{\mathcal C}
\def\V{\mathcal V}
\def\X{\mathcal X}
\def\G{\mathcal G}
\def\U{\mathcal U}
\def\R{\mathcal R}
\def\Z{\mathcal Z}
\def\CalS{\mathcal S}
\def\Hpi{(\pi, {\mathcal H})}
\def\unirep{unitary representation}
\def\supp{{\rm supp}}
\def\NN{{\mathbf N}}
\def\ZZ{{\mathbf  Z}}
\def\CCC{{\mathbf  C}}
\def\RRR{{\mathbf R}}
\def\QQ{\mathbf Q}
\def\RR+{{\mathbf R}^*}
\def\TT{\mathbf T}
\def\HH{\mathbf H}
\def\KK{\mathbf K}
\def\kk{\mathbf k}
\def\PP{\mathbf P}
\def\GG{\mathbf G}
\def\HS{\rm HS}
\def\Q_p{{\mathbf Q}_p}
\def\OO{\mathbf O}
\def\UU{\mathbf U}

\def\vfi{{\varphi}}
\def\tout{\quad \quad \text{for all}\quad }

\def\norm#1{\left\Vert #1 \right\Vert}
\def\abs#1{\left\vert #1 \right\vert}

\def\CF{{\cal F}}

\def\nil{\La\backslash N}
\def\tor{\La [N,N]\backslash N}
\def\Aut{{\mathrm Aut}}
\def\Aff{{\mathrm Aff}}
\def\Affnil{{\mathrm Aff}(\nil)}
\def\Autnil{\Aut (\nil)}
\def\AffT{{\mathrm Aff}(T)}
\def\AutT{\Aut (T)}
\def\paut{p_{\rm a}}
\def\ind{{\mathrm Ind}}
\def\Ker{{\mathrm Ker}}
\def\End{{\mathrm End}}
\def\Tr{{\mathrm Tr}}
\def\hN{\widehat N}
\def\rhonil{\rho_{\nil}}
\def\Utor{U_{\rm tor}}
\def\ZC{{\rm Zc}}
\def\deg{{\mathrm deg}}

\section{Introduction}
Let $G$ be a separable  locally compact group.  The study of actions of $G$ on 
various spaces is of course of fundamental importance, both for
the understanding of properties of the groups and the spaces under consideration.
When $G$ acts on such a space $X$, which might be a 
manifold or a graph, there is often a   positive measure $m$ on 
the measurable subsets of $X$
which is quasi-invariant under $G$ and sometimes even invariant.
We will consider two kinds of actions, which require different approaches: 
\begin{itemize}
\item[(A)] the case where $m(X)=\infty,$ that is, $m$ is an infinite measure;
\item[(B)] the case where $m$ is a finite measure which is $G$-invariant.
\end{itemize} 
In case (B), we   may of course  assume that $m$  is  a probability measure.
Attached to these data, there is a natural unitary representation $\pi_X$ of $G$ 
on the Hilbert space $L^2(X,m)$; see Section~\ref{S:UnirepAction}.
In case $m$  is a probability measure and is $G$-invariant, the space $\CCC {\mathbf 1}_{X}$ of the constant functions on $X$
is contained in $L^2(X,m)$ and is  $G$-invariant as well as its orthogonal
complement 
$$
L^2_0(X,m)=\left\{f\in L^2(X,m)\, :\, \int_{X} f(x) dm(x)=0\right\}.
$$
%In case $\mu$ is infinite, we set $L^2_0(X,m)= L^2(X,\mu).$

%\begin{definition}
In this survey, we  will be concerned with  the study of a spectral rigidity property 
for actions $G\curvearrowright X$,   in each of the  two situations (A) and (B) above.
 We say that  the action of $G$  on
 $(X,m)$ has the Spectral Gap Property 
 if the representation $\pi_X$ on $L^2_0(X,m)$ 
 does not  have almost invariant vectors, 
that is, if there is no sequence of unit vectors
 $\xi_n$  in  $ L^2_0(X,m)$  such
 that $\lim_n\Vert \pi_X(g)\xi_n-\xi_n\Vert=0$ 
 for all  $g$ in $G.$
 Here, we have denoted the space $L^2(X,m)$ by $L^2_0(X,m)$
 in the case where  $m$ is an infinite measure.
 
 Two classes of groups are distinguished with respect to this Spectral Gap Property: 
 amenable groups and Kazhdan groups. Indeed,
 an amenable group $G$
 never has the Spectral Gap Property (in case (B), we have
 to assume that  $G$ is countable and $X$ non-atomic) and Kazhdan groups
 always have the Spectral Gap Property for ergodic actions 
 (see Corollary~\ref{Pro-AmenableCoam} and Theorem~\ref{Theo-Rosenblatt}
 as well as Corollary~\ref{Pro-KazhdanCoam}).
 We will first recall some basic facts on these classes of groups.
 Then, we turn to  the negation of the Spectral Gap Property 
 in case (A) above.  Actions $G\curvearrowright X$ 
 without the Spectral Gap Property are called \emph{co-amenable,}
 as they are characterized by the existence
 of an invariant mean on $L^\infty(X,m).$ Such actions have been studied,  
 with various degrees of generality, by several authors
 (see \cite{Eymard}, \cite{Greenleaf},  \cite{GuivAsym}, among others).
 We will give a comprehensive account about the main characterizations of such actions
 in Section~\ref{S:Actions-Coam}.
 %%Ajout
 %%
  In the case where $G$ is a non compact simple 
 Lie group, the co-amenable proper subgroups $H$ (that is, such that  $G\curvearrowright G/H$
 is co-amenable) are Zariski dense in $G$ and are characterized  
 as the discrete subgroups with the maximal critical exponent (Theorems~
 \ref{Theo-CoAmZariski} and \ref{Theo-CoAm-CritExp}).
 %%%
 
 Next, we will deal with the Spectral Gap Property
 in case (B) above, that is, for actions with an invariant probability measure.
 We will review some  specific examples 
 for  actions with the Spectral Gap Property  on a homogeneous $X=H/\La$, 
 where $\La$ is  a locally compact group and $\La$ is a 
 lattice in $H;$ here the action of $G$ on $H/\La$ is given by left translations
 or by automorphisms, that is, by a homomorphism
 $G \to H$ or a homomorphism $G \to \Aut(X)$, where
 $\Aut(X)$ is the subgroup
of continuous automorphisms  $\varphi$ of $H$ such that $\varphi (\La) =\La.$
Apart from a few exceptions (such as Bernoulli actions), 
showing the Spectral Gap Property for these actions is usually a difficult 
problem. 
As a novelty, we establish by  completely elementary means
(using an idea from \cite{BeLu}) the Spectral Gap Property
for the action  of $G= PGL_2(\FF_q((t^{-1})))$
on  $X=PGL_2(\FF_q((t^{-1})))/PGL_2(\FF_q[t])$,
where $\FF_q$ is the finite field with $q$ elements and $\FF_q((t^{-1}))$  the local field  of  Laurent series
(see Section~\ref{S:PGL}). As a crucial tool in our approach, we prove and use a Cheeger type inequality for Markov chains
established in  \cite{LaSo} and  \cite{Sinclair}.

 We approach the  Spectral Gap Property mainly in terms of  averaging 
 operators, also known as Markov operators. Let $\mu$ be a probability
measure on $G$  and $\pi_X(\mu)$  the convolution operator  defined
on $L^2_0(X,m)$ by 
$$
\pi_X(\mu)f =\int_{G}  \pi_X(g)f d\mu(g) \tous f\in L^2_0(X,m).
$$
We have  $r_\spec(\pi_X(\mu)) \leq 1$  for the  spectral 
radius $r_\spec(\pi_X(\mu))$ of  $\pi_X(\mu)$.
Assume that  $\mu$ is  adapted,
that is, the support $\supp (\mu)$ of $\mu$ generates a dense subgroup of  $G$. Then the action of $G$ on $X$ has the Spectral Gap
Property if and only if $r_\spec(\pi_X(\mu))<1.$ 

The point of view of Markov operators  is relevant for applications in probability theory. 
Given an action $G\curvearrowright (X,m)$  and a probability measure
$\mu$ on $G,$ consider a sequence of independent $\mu$-distributed 
random variables $X_n$ with values in $G$
and the corresponding random products
$S_n=X_n\dots X_1$ for $n\in\NN.$ 
This defines a random walk on $X,$ given by the transitions probabilities
$$A\mapsto p(x, A)= \int_G  {\mathbf 1}_A(g^{-1} x)  d\mu(g) =( \pi_X(\mu)  {\mathbf 1}_A) (x),$$
for $x$ in $X$ and $A$ a measurable subset of $X.$

The  Spectral Gap Property  has several interesting applications in ergodic theory;
for instance, if $m$ is a probability measure, then, for every $f\in L^2(X,m),$ 
the sequence
of functions $x\mapsto \EE (f(S_n(x)))$ converges to $\int_X f dm$ in 
the $L^2$-norm,
with an exponentially fast rate of order $\lambda^n$ with $\lambda= \Vert \pi_X(\mu)\Vert.$ 
Other ergodic theoretic applications  to random walks  (see  \cite{CoGu2}, \cite{CoLe}, \cite{FuSh}, \cite{GoNe} and  \cite{GuivSpec})
include  the rate of convergence in the
random  ergodic theorem, pointwise ergodic theorems, analogues of the law of large numbers and of the central limit theorem, etc.
  Another application of the Spectral Gap Property is the uniqueness of $m$ as $G$-invariant mean
on $L^\infty(X,m);$ for this as well as for further
applications, see \cite{BHV}, \cite{BoGa},\cite{Lubotzky}, \cite{Popa}, \cite{Sar}.
To illustrate the use of the Spectral Gap Property in both situations (A) and (B),
we present  in the last section of this survey (Section~\ref{S:Appl}) two such applications:
  one to expanders graphs and one  to the escape rate of random matrix products.
 We also discuss (Section~\ref{S:QuantitativeSG}) the question  of quantifying the  Spectral Gap Property, that is,
 giving upper bounds for the norm or the spectral radius of the Markov operator  $\pi_X(\mu)$.

%%%%%%%%%%%%%%%%%%%%%%%%%%%%%%%%%%%%%%%%%%%%%%%%%%
\section[Measure]{Group actions on measure spaces and associated  representations}
\label{S:UnirepAction}
Let $G$ be  locally compact group, which we always assume to be second countable.
We will then say that $G$ is a \emph{separable locally compact group}. The group
$G$ has an (essentially unique) Haar measure $\lambda$, that is, a $\sigma$-finite  measure $\lambda$ on the Borel subsets of $G$
which is invariant under left translations.

Let $(X, m)$ be a measure space, where
 $m$ is a positive $\sigma$-finite measure on a fixed $\sigma$-algebra
 of subsets of $X.$  
We will only consider actions of $G$  of $X$ which are measurable, that is, 
 actions $G\curvearrowright X$  for which
$$G\times X\to X\quad (g,x)\mapsto gx$$
is measurable. 

We will always assume that the measure $m$ is \emph{quasi-invariant:} 
$m$ and  its image  $gm$ under $g$ are equivalent measures
(that is, $m$ and $gm$ have the same sets of measure $0$) for every $g\in G.$ 

Our actions will usually be ergodic group actions.
\begin{definition}
\label{Def-Ergodic}
The action $G\curvearrowright X$  on the measure space $(X,m)$
 is \emph{ergodic} if there are no nontrivial invariant subsets of $X$ in the following sense:
 if $A$ is a measurable subset of $X$ such that  $gA=A$ for all $g\in G,$
 then $m(A)=0$ or $m(X\setminus A)=0.$
   \end{definition}
   Ergodicity can be expressed in terms of functions on $X$
   as follows.
   We say that a measurable function $f:X\to \RRR$ is $G$-invariant
   if, for every $g\in G,$ we have $f(gx)=f(x)$ for $m$-almost
   every $x\in X.$ An action $G\curvearrowright X$  is ergodic if and only if 
   every $G$-invariant function is constant $m$-almost everywhere
   (see  Theorem~1.3 in \cite{BeMa}).

   \subsection{Examples of  group actions on measure spaces} 
\label{Exa-Actions}
We list some examples of  group actions, which will appear
throughout this survey.
\begin{enumerate}
\item  Let $\Ga$  be  a countable  group, 
$X= \{0, 1\}^\Ga,$ and  $m$ the probability measure
$m= \bigotimes_{\ga\in\Ga}\nu$ for
the measure $\nu$ on $\{0, 1\}$ given by $\nu(\{0\})= \nu(\{1\})=1/2.$ The
Bernoulli action of $\Ga$ is the measure preserving action
$\Ga \curvearrowright X$ defined by shifting coordinates:
$$
\ga (x_{\delta})_{\delta\in \Ga}= (x_{\ga^{-1}\delta})_{\delta\in \Ga}.
$$
\item Let $H$ be a separable locally compact group  and $L$ a closed subgroup of $H.$
The homogeneous space $X=H/L$  has a unique (up to equivalence) non-zero
$\sigma$-finite measure $m$ on its Borel subsets, which is quasi-invariant under the action of $H$
by left translations (see \cite[(2.59)]{Fol}).  
Every subgroup $G$ of $H$ acts 
on $X$ by left translations.

\item An important special case in the previous example arises when $L$ is 
a lattice in the locally compact group $H$, that is, 
$L $ is a discrete subgroup of $H$
and there exists a   $H$-invariant probability measure $m$ on the Borel subsets of
$H/L$. Every subgroup $G$ of $H$ acts in measure preserving way 
on $X=H/L$ by left translations.
Examples are given by $H= \RRR^n$ and $L = \ZZ^n,$ 
in which case  $X= \TT^n$
is the $n$-torus or by $H= SL_n(\RRR)$ and $L =SL_n(\ZZ),$ 
in which case  $X$ is the space of unimodular lattices in $\RRR^n.$

\item Let $H$,  $L$  and $m$ be as in Example 2.
Let  $\Aut(H)$ be the group of continuous automorphisms of 
  $H.$   The subgroup ${\Aut} (H/L)$ of $\Aut(H)$, defined by  
$${\Aut} (H/L) = \{ \vfi \in {\Aut}(H) \, : \,  \vfi (L) =L\},$$
 acts on $H/L$, leaving $m$ quasi-invariant.
 If $m$ is  $H$-invariant and finite, then $m$ is also invariant under ${\Aut} (H/L).$

In the example  given  by $H= \RRR^n$, $L= \ZZ^n$ 
and $X= \TT^n$, the group  ${\Aut} (H/L)$ can be identified with $GL_n(\ZZ).$
Other examples arise when $H$ is a  nilpotent Lie group  and $L$ is a lattice
in $H$ (see Section~\ref{S:Nil}).

\end{enumerate}

\subsection{Unitary representations associated to actions}
\label{SS:Rep}
Let $G$ a separable locally compact group and 
$G\curvearrowright X$ an action on a $\sigma$-finite measure space $(X, m)$,
with $m$ quasi-invariant.

For $g\in G,$   denote by $c(g,x) =\dfrac{dgm}{dm}(x)$  the Radon-Nikodym derivative of $gm$
with respect to $m.$  The mapping 
$\pi_X: G\to \B(L^2(X,m)),$ defined by 
 $$
\pi_{X}(g)f(x)=c(g^{-1},x)^{1/2}f(g^{-1}x) \tout f\in L^2(X,m), x\in X,
$$ 
is a continuous unitary   representation of $G$ on  $L^2(X, m),$
often called the \emph{Koopman representation} associated
to   $G\curvearrowright X$  (for more details,
see \cite[A.6]{BHV}).

Assume now that $m$ is a $G$-invariant probability measure.
Then 
 $$
 L^2_0(X)= \{f\in L^2(X)\ : \  \int_X fdm =0 \} = (\CCC {\Un}_X)^\perp
 $$
  is $G$-invariant.
 Denote again by $\pi_X$ the restriction of Koopman representation to   $L^2_0(X).$
 Then $G\curvearrowright X$ is ergodic if and only if  
$\pi_X$ has no non-zero invariant vectors:  $L^2_0(X)^G=\{0\}.$

Most of  the actions  we consider  will
be  mixing in the following sense. The probability measure preserving action $G\curvearrowright X$
 is called \emph{mixing} if  $\pi_X$ is a $C_0$ representation:  
  for all $f_1, f_2\in L^2_0(X),$ the matrix coefficient
$$C_{f_1, f_2}:  G\to \CCC, \qquad g\mapsto \langle \pi_X(g) f_1, f_2\rangle$$
belongs  to $C_0(G)$, that is, $\lim_{g\to \infty}  \langle \pi_X(g) f_1, f_2\rangle=0.$
Of course, mixing actions are ergodic.

\begin{example} 
%\begin{enumerate}
(i) The Bernoulli action of an infinite countable group $\Ga$ on 
$X= \{0, 1\}^\Ga$   is mixing.

\n
(ii) Let $\Ga$ be a subgroup of $GL_n(\ZZ).$
The action of $\Ga$ on the $n$-torus  $\TT^n =\RRR^n/\ZZ^n$
is ergodic if and only if  every $\Ga^t$-orbit in $\ZZ^n$ is infinite;
this action is mixing if and only if every point stabilizer for the
action of  $\Ga^t$ on  $\ZZ^n$
is finite. 
%\end{enumerate}
\end{example}
The claims in the previous examples follow from the next  proposition,
which illustrates the power  of the use of the Koopman representation.
Concerning  the first example,  observe that  
we can view    $\{0, 1\}^\Ga$ as the   compact abelian
group $A=\prod_{\ga\in \Ga} \ZZ/ 2\ZZ$ and $\Ga$ as a subgroup of 
$\Aut(A).$

\begin{proposition}
\label{Pro-Erg-Mixing}
Let $A$ be a compact abelian group with normalized Haar measure $m$ and $\Ga$ a subgroup
of $\Aut (A).$  Let $\widehat A$ be  the (discrete) dual group of $A$ (on which 
$\Aut(A)$ acts naturally). 
\begin{itemize}
\item The action
$\Ga\curvearrowright A $ is ergodic  if and only if  every $\Ga$-orbit in 
$\widehat A\setminus\{1\}$ is infinite;
\item the action $\Ga\curvearrowright A $ is mixing  if and only if   $\Ga$ acts properly on
$\widehat A\setminus\{1\}$ (that is, point stabilizers are finite).
\end{itemize}
\end{proposition}
\proof 
By Fourier transform, we have  isometric isomorphisms 
$L^2(A) \cong \ell^2(\widehat A)$ and
$L^2_0(A)\cong \ell^2(\widehat{A}\setminus\{1\})$ as $\Ga$-representations.

Now, $ \ell^2(\widehat{A}\setminus\{1\})^{\Ga}$ consists exactly of 
the $\Ga$-invariant $\ell^2$-functions 
on  $\widehat{A}\setminus\{1\}.$
Moreover, the $\Ga$-representation on  $\ell^2(\widehat{A}\setminus\{1\})$
is $C_0$ if and only if  all point stabilizers are finite.$\bsq$

The following result is the celebrated Howe-Moore theorem from \cite{HoMo},
which shows that, as a rule, ergodic actions of  simple  Lie groups 
are  mixing.
\begin{theorem}
\label{Theo-HoweMoore}
\textbf{(Howe-Moore Theorem)} 
Let $G$ be a connected simple Lie group with finite center
%the group of $\kk$-rational points
%of a connected simple algebraic group over a local field $\kk$ 
and $\Hpi$ a unitary representation of
$G$ on a Hilbert space $\H.$  Assume that $\H^G=\{0\}.$
Then $\pi$ is a $C_0$-representation: the matrix coefficients
$$C_{\xi, \eta}: G\to \CC, \quad g\mapsto \langle \pi(g) \xi, \eta\rangle$$
belong to  $C_0(G)$ for all $\xi, \eta\in\H.$
\end{theorem}

Here is one  striking consequence of the Howe-Moore theorem.
\begin{corollary}
\label{Cor-HoweMoore}
Let $G$ be as in Theorem~\ref{Theo-HoweMoore} and $\Lambda$ a lattice in $G.$
Let $H$ be a subgroup with a non-compact closure in $G.$ 
Then $H\curvearrowright  G/\Lambda$ is ergodic. Morevover,
$H\curvearrowright  G/\Lambda$ is even mixing 
if $H$ is closed.
\end{corollary}

\proof 
Set $X= G/\Lambda.$
Obviously, we have  $L^2_0(X)^G=\{0\}.$
Let $f\in L^2_0(X)^H$. Then 
$f\in L^2_0( X)^{\overline H}.$
Now, by Theorem~\ref{Theo-HoweMoore}, the matrix coefficient
$C_{f,f}$
belongs to $C_0(G)$.
Since $\overline H$ is not compact, it follows that $f=0.$ $\bsq$

%%%%%%%%%%%%%%%%%%%%%%%%%%%%%%%%%%%%%%%%%%%%%%%%%%%
\section{The Spectral Gap Property}
\label{S: SGP}
Let $G$ a separable locally compact group and 
$G\curvearrowright X$ an action on a measure space $(X, m)$,
with $m$ quasi-invariant.
%\subsection 
%{Actions with  spectral gap property}

We set $L^2_0(X)= (\CCC{\Un}_X)^\perp$ in case $m$ is a $G$-invariant probability measure,
and $L^2_0(X)= L^2(X,m)$ otherwise.
 The  corresponding unitary representation of $G$
 on $L^2_0(X)$ will always be denoted by $\pi_X.$
 
 \begin{definition}
 \label{Def-SGP}
  \textbf{(Actions with  the Spectral Gap Property)}
 The action of $G$ on $X$ has the \emph{Spectral Gap Property}, if there exists a compact set
 $Q$ of $G$ and $\eps>0$ such that
 $$
 \sup_{s\in Q} \Vert \pi_X(s) f-f\Vert \geq \eps \Vert f\Vert \tout f\in L^2_0(X).
 $$
 \end{definition}
 
 We can extend this definition to   arbitrary unitary representations.

 \begin{definition}
 \label{Def-RepSGP}
 \textbf{(Representations with  the Spectral Gap Property)}
A unitary representation $\Hpi$ of $G$
 has the Spectral Gap Property, if there exists a  compact subset
 $Q$ of $G$ and $\eps>0$ such that
 $$
 \sup_{s\in Q} \Vert \pi(s) \xi-\xi\Vert \geq \eps \Vert \xi\Vert \tout \xi\in \H.
 $$
   \end{definition}
 \n
 \begin{remark}
 \label{Rem-SGP}
 (i) \textbf{(Negation of the Spectral Gap Property)} The unitary representation $\Hpi$ of $G$ does not have the Spectral Gap Property if, for  
 every pair $(Q,\eps),$ where $Q$ is a compact subset of $G$ and $\eps>0,$
 there exists a unit vector  unit $\xi\in\H$ which is $(Q,\eps)$-invariant:
 $$\sup_{s\in Q} \Vert \pi(s) \xi-\xi\Vert < \eps.$$
 Since $G$ is $\sigma$-compact, observe that $\Hpi$  does not have the Spectral Gap Property
 if and only if there exists a sequence of unit vectors $(\xi_n)_n$ in $\H$ such that
 $$\lim_n \Vert \pi(g) \xi_n-\xi_n\Vert= 0 \tout g\in G.$$
 (The ``if" part of the previous statement follows from  a standard Baire category argument.)
 
\n
(ii) The negation of  the Spectral Gap Property may be formulated in terms of Fell's notion of weak containment: 
 $\pi$  does not have  the Spectral Gap Property  if and only if the trivial representation $1_G$
 is weakly contained in $\pi$.
 Recall that a  unitary representation $\rho$ is said to be \emph{weakly contained} in another
unitary representation $\sigma,$
if every diagonal matrix coefficient $C^{\rho}_{\xi,\xi}$ of $\rho$
can be approximated, uniformly on compact subsets of $G,$ by
convex combinations of diagonal matrix coefficients of $\sigma$
 (see Appendix F in \cite{BHV}).
 
   \end{remark}
  
 \subsection{ Spectral Gap Property in terms of averaging operators}
 Let $\Hpi$ be a unitary representation of a locally compact group $G$ and $\mu$ a 
probability measure on the Borel subsets of  $G.$ Define the averaging operator
$\pi(\mu)\in \B(\H)$ by 
$$\pi(\mu)\xi = \int_G \pi(g)\xi d\mu(g) \tout \xi\in \H.$$
Then clearly  $\Vert \pi(\mu)\Vert \leq 1$
and hence $r_\spec(\pi(\mu))\leq 1$
for the spectral radius  $r_{\spec}(\pi(\mu))$ 
of $\pi(\mu).$

We say that $\mu$ is \emph{adapted}
 if the subgroup generated
by the support $\supp (\mu)$ of $\mu$ is dense
in $G.$ 
We will also consider
the stronger condition
that $\supp (\mu)$ is not contained
in the coset of a proper closed subgroup
of $G$. In this case, we say that $\mu$ is \emph{strongly adapted}.

It is easy to see that  $\mu$ is strongly adapted if and only if 
the convolution product  $\check{\mu}\ast \mu$
is adapted, where  $\check{\mu}$ is defined by $d\check{\mu}(g)=d\mu(g^{-1}).$
We say that $\mu$ is absolutely continuous if it is
absolutely continuous with respect to a Haar measure on $G.$

\begin{proposition}
\label{SGP-Norm}
 Let $\Hpi$ be a unitary representation of the separable locally compact group $G.$
 The following statements are equivalent:
\begin{itemize}
\item[(i)] $\pi$ has the Spectral Gap Property;
\item[(ii)] $\Vert \pi(\mu)\Vert <1$ for any   (or for some) probability measure $\mu$ on  $G$ 
which is strongly adapted and absolutely continuous;
\item[(iii)]  $1$ is not a spectral value of $\pi(\mu)$ for any  (or for some) probability measure $\mu$ on  $G$ 
which is  adapted and absolutely continuous.
\end{itemize}

\end{proposition}

 \proof 
 We   just give the proof of the equivalence of (i) and (ii)  in the case where $G=\Ga$ discrete,
 and refer  to Proposition G.4.2 in \cite{BHV} for the complete proof.
 
 Assume that $\pi$ does not have the Spectral Gap Property; so,
 there exists a sequence $\xi_n$ of unit vectors in $\H$ such that
 $$ 
 \lim_n\Vert \pi(\ga) \xi_n-\xi_n\Vert=0 \tout \ga\in \Ga.
 $$
 Summing  against $\mu$ gives
\begin{align*}
\lim_n\Vert \pi(\mu)\xi_n-\xi_n\Vert &\leq \lim_n \sum_{\ga\in\Ga}\mu(\ga)  \Vert \pi(\ga) \xi_n-\xi_n\Vert=0.
\end{align*}
So, $1$ is in the spectrum of $\pi(\mu)$. Since $\Vert \pi(\mu)\Vert \leq 1,$
we have $\Vert \pi(\mu)\Vert =1.$

 Conversely, assume that $\Vert \pi(\mu)\Vert =1.$ 
  Since $\mu$ is strongly adapted, $\supp (\check\mu\ast \mu)$ generates $\Ga.$ 
 Now,
 $$
 \Vert\pi(\check\mu\ast \mu)\Vert= \Vert\pi(\mu)^*\pi( \mu)\Vert= \Vert \pi(\mu)\Vert^2=1
 $$
 and  $\pi(\check\mu\ast\mu)$ is a positive self-adjoint operator.
 Hence,  $1$ is a spectral value of $\pi(\check\mu\ast\mu)$ and there exists a sequence
 $\eta_n$ in $\H$ of approximate eigenvectors, that is, a sequence
 of unit vectors $\eta_n$ with $\lim_n\Vert\pi(\check\mu\ast\mu)\eta_n-\eta_n\Vert=0.$
  Then 
  $$
  \lim_n\Vert\pi(\ga)\eta_n-\eta_n\Vert=0
  $$
  for all $\ga\in \supp(\check\mu\ast\mu)$ and therefore for all $\ga\in\Ga. \bsq$

\subsection{Kazhdan's Property (T)}
 \label{SS:Kazhdan}
 We review a few basic facts on Kazhdan groups,
 with an emphasis on the existence of a spectral gap 
 for averaging operators. We  refer to 
 the monograph \cite{BHV} for missing  details.
 
 \begin{definition}
 \textbf{(\cite{Kazhdan})}
 A locally compact group  $G$ has \emph{Property (T)} or is a \emph{Kazhdan group},  if every unitary representation
 $\Hpi$ of $G$ with $\H^G=\{0\}$
 has the Spectral Gap Property.
 \end{definition}
 
 Recall that a local field $\kk$ is a locally compact non discrete field.
As is well-known,  a local field 
is isomorphic either to  $\RRR$, to $\CCC$, to a finite extension of the field of $p$-adic numbers
$\Q_p$, or to the field $k((X))$  of Laurent series over a finite field $k$.

 \begin{theorem}
 \label{Theo-T-Rank2}
 \textbf{(\cite{Kazhdan})}
 Let $\kk$ be a local field and $G$  the group of $\kk$-rational points of 
a simple algebraic group over $\kk$ with  $\kk$-rank at least two.
Then    $G$ has  Property (T).
 \end{theorem}
 
 \begin{example}
 The groups $SL_n(\RRR)$ for $n\geq 3$ and 
 $Sp_{2n}(\RRR)$ for $n\geq 2$ have Property (T).
 \end{example}

 When  $G$ has Property (T), we  can find a uniform pair $(Q, \eps)$ as in Definition~\ref{Def-RepSGP}
  for all representations without invariant vectors.
 \begin{proposition}
 \label{Pro-T-uniform}
 Let $G$ be a Kazhdan group.
Then there exists a pair $(Q,\eps)$ with $\eps>0$ and $Q$ a compact 
 subset of $G$ so that, for every  unitary representation
 $\Hpi$ of $G$ with $\H^G=\{0\},$ we have
 $$
 \sup_{s\in Q} \Vert \pi(s) \xi-\xi\Vert \geq \eps \Vert \xi\Vert \tout\xi\in \H.
 $$
 The pair $(Q,\eps)$ is called a \emph{Kazhdan pair}
 \end{proposition}
 \proof
 Assume, by contradiction, that no such a pair exists.  
 Denote by $I$ the set of  pairs $(Q,\eps)$
 with $\eps>0$ and $Q\subset G$  compact.
 So, for every $i=(Q,\eps) \in I,$ there exists a unitary representation
$(\pi_i, \H_i)$ of $G$ with $\H_i^G=\{0\}$ and a unit vector $\xi_i\in \H_i$
 such that
$$
 \sup_{s\in Q} \Vert \pi_i(s) \xi_i-\xi_i\Vert < \eps.
 $$
Then  the direct sum $\bigoplus_{i\in I} \pi_i$ is  a unitary representation of $G$
 in $\bigoplus_{i\in I} \H_i$ which has almost invariant vectors 
and has no non-trivial invariant vectors. This is a contradiction.$\bsq$

 Here is one important application of Kazhdan's Property (T).
 \begin{theorem}
 \label{Theo-T-CompactGen}
 \textbf{(\cite{Kazhdan})}
 If $G$ has Property (T), then $G$ is compactly generated.
 In particular, a discrete group with  Property (T) is finitely generated.
 \end{theorem}
 \proof
Let $\C$ be the family of open and compactly generated subgroups 
of $G.$ 
Consider the representation $\pi=\bigoplus_{H\in \C} \pi_{G/H},$
where $\pi_{G/H}$ is the regular representation on $\ell^2(G/H).$
Then $\pi$  does not have the Spectral Gap Property. So, $\pi$ has a non-zero 
invariant vector. This implies that $G/H$ is compact for some
$H\in \C$ and hence that $G$ is compactly generated. $\bsq$

%\subsection{Property (T) is inherited by lattices}
An important feature of Property (T) is that it is inherited by lattices.
The proof involves induced representations.

Let $G$ a separable locally compact group and $\Ga$ a lattice in $G.$
Let $\Hpi$ be a unitary representation of $\Ga.$
The induced representation $\widetilde\pi= \ind_\Ga^G\pi$ is the unitary representation
of $G$ which can be defined as follows.
Let $X\subset G$ be a Borel set which is a fundamental domain
for the  action of $\Ga,$ so that 
$$G= \coprod_{\ga\in \Ga} \ga X.$$
Given $g\in G$ and $x\in X,$ there are uniquely determined elements
$c(x,g)\in \Ga$ and $x\cdot g\in X$ such that
$$
xg= c(x,g) (x\cdot g) .
$$
Then $\widetilde\pi= \ind_\Ga^G\pi$ is
defined on the space $\widetilde \H$ 
of measurable maps $F\in L^2(X, \H) $
by 
$$
(\widetilde\pi(g) F)(x)= \pi(c(x,g))(F(x\cdot g)).
$$
 For more details on induced representations, see \cite{Mackey}.

\begin{proposition}
\label{Pro-T-Lattices}
Let $G$ a separable locally compact group and $\Ga$ a lattice in $G.$
Let $\Hpi$ be a unitary representation of $\Ga.$
Assume that $\pi$ does not have the Spectral Gap Property. Then $\ind_\Ga^G\pi$ does not have the Spectral Gap Property.
\end{proposition}
\proof
Since $\pi$ does not have the Spectral Gap Property, there exist 
a sequence of unit vectors $\xi_n\in \H$ with
$\lim_n\Vert \pi(\ga) \xi_n -\xi_n\Vert=0 $ for all $\ga\in \Ga$. 
Define $F_n\in \widetilde{\H}$ by 
$F_n(x)= \xi_n$ for all $x\in X.$
Then $\Vert F_n\Vert =1$ and 
\begin{align*}
\Vert \widetilde\pi(g) F_n -F_n \Vert^2&=  \int_X \Vert \pi(c(x,g)) \xi_n -\xi_n\Vert^2 dm(x),\\
\end{align*}
where $dm(x)$ denotes the Haar measure restricted to $X,$ normalized
 by $m(X)=1.$
Hence, $\lim_n\Vert \widetilde\pi(g) F_n -F_n \Vert=0.$
So, $\widetilde\pi= \ind_\Ga^G\pi$ does not have the Spectral Gap Property.$\bsq$
 
We deduce from the previous proposition the following important theorem,
a key result for producing examples of \emph{discrete}  Kazhdan groups.
\begin{theorem}
\label{Theo-T-Lattices}
\textbf{(\cite{Kazhdan})}
Let $G$ a locally compact group and $\Ga$ a lattice in $G.$
If $G$ has Property (T), then $\Ga$ has property (T).
\end{theorem}

\proof
Let $\Hpi$ be  a unitary representation of $\Ga$ without  the Spectral Gap Property.
Then $\widetilde\pi=\ind_\Ga^G\pi$ does not have the Spectral Gap Property, by the previous
proposition. So,   $\widetilde \H^G\neq \{0\}.$ On the other hand, it is immediate
from the definition of $\ind_\Ga^G\pi$  that 
$\widetilde \H^G$ consists exactly of the constant mappings   $F$ in $L ^2(X, \H)$ with 
$F(x)\in \H^\Ga$ for $m$-almost every $x\in X.$
Hence, $\H^\Ga\neq \{0\}. \bsq$

%\n
%\textbf{Example:}
\begin{example} 
The discrete groups  $SL_n(\ZZ)$ for $n\geq 3$ and $Sp_{2n}(\ZZ)$ for $n\geq 2$  have Property (T),
as they are   lattices in the Kazhdan groups $SL_n(\RRR)$  and $Sp_{2n}(\RRR),$
respectively.
\end{example}

\subsection{Uniform Spectral Gap Property for  Kazhdan groups}
\label{SS: SPG-T}
Let $G$ be a locally compact group with Property (T)
and $\pi$ a unitary representation of $G$ without non-zero invariant vectors.
Then, by  definition, $\pi$ has the Spectral Gap Property.
Hence, given  a  strongly adapted and absolutely continuous probability measure $\mu$ on  $G$,
we have $\Vert \pi(\mu)\Vert <1$, by Proposition~\ref{SGP-Norm}.
In fact, we can find a uniform bound  for $\Vert \pi(\mu)\Vert$,
independent of $\pi.$

\begin{theorem}
\label{Theo: SPG-T}
\textbf{(Uniform Spectral Gap Property)}
Let $G$ be a separable locally compact group with Kazhdan's Property (T) and
 $\mu$  a 
strongly adapted and absolutely continuous probability measure  on  $G$.
Then there exists $C<1$ such that 
$\Vert \pi(\mu)\Vert <C$ for every  unitary representation $\Hpi$ of $G$ without non-zero invariant vectors.
\label{Theorem}
\end{theorem}
\proof
 Assume, by contradiction, that  this is not the case. Then there exists
a sequence  $(\pi_n, \H_n)$ of  unitary representations
 of $G$ without non-zero invariant vectors
 such that $\lim_n \Vert \pi_n(\mu)\Vert=1.$
  Then  $\pi=\bigoplus_{n} \pi_n$ is  a unitary representation of $G$
 on $\H= \bigoplus_{n} \H_n$ which has no non-trivial invariant vectors
 and for which $\Vert \pi(\mu)\Vert =1$. This is a contradiction.
 $\bsq$
 
 \begin{remark}
  Let  $G$ be a  locally compact group with 
 the following property:  every  ergodic measure preserving action  $G\curvearrowright X$ 
 on a probability space $(X,m)$ has the Spectral Gap Property. 
 Of course, Kazhdan groups have this property. In fact, it was shown
 in \cite{ConnesWeiss} that this property  characterizes the class of Kazhdan groups.
 
  \end{remark}
 
\subsection{Amenable groups}
\label{AmenableGroups}
Amenability of locally compact groups may be expressed in 
several equivalent  ways. We give a brief review of 
a few number of these equivalent reformulations. 

\begin{definition}
 \label{Def-Amenable}
 Let $G$ be a locally compact group and denote by $m$ a  left Haar measure 
 on $G.$ The group $G$ is \emph{amenable}
 if  there exists a $G$-invariant mean 
on $L^\infty(G,m),$  that is,  
a positive linear functional $M$
on  $L^\infty(G,m)$ such 
that $M({\mathbf 1}_G)=1$ and $M({}_g\vfi) =M(\vfi)$
for all $g$ in $G$ and $\vfi$ in $L^\infty(G,m),$
where ${}_g\vfi(x)=\vfi (g^{-1}x).$
  \end{definition}

  We mention the following useful characterizations of amenable groups
  and refer to Appendix G in \cite{BHV} for proofs.
  
 \begin{proposition} 
 \label{Pro-Amenable}
 Let $G$ be a locally compact group.
 The following properties  are equivalent:
 \begin{itemize}
 \item[(i)] $G$  is amenable; 
 \item[(ii)] every continuous action $G\curvearrowright \C$ 
 by affine mappings on a  non-empty compact  convex subset $\C$ of a locally convex topological vector space
 has a fixed point;
\item[(iii)] for every continuous action $G\curvearrowright X$  on a compact non-empty 
 set $X,$ there exists a $G$-invariant probability  measure on the Borel subsets of $X.$  
 \end{itemize}
 \end{proposition}
 
 Examples of amenable groups  include  abelian groups 
 and more generally solvable groups. On the  other hand, 
 free non-abelian groups are not amenable.
 
We now rephrase amenability of a group in terms of 
 its regular representation.

 Let $G$ be a locally compact group with left Haar measure $m.$
 The left regular representation of $G$ is the unitary representation 
 $\pi_G$ which is defined by left translations on $L^2(G,m).$ 
  The following result is due to  Hulanicki (\cite{Hulanicki}) and Reiter (\cite{Reiter}).
  
 \begin{theorem}
 \label{HulaReiter}
 \textbf{(Hulanicki-Reiter's Theorem) }
 The following statements  are equivalent
 for a  locally compact group $G.$
 \begin{itemize}
 \item[(i)] The group  $G$  is amenable; 
 \item[(ii)] the regular representation $\pi_G$ does not have the Spectral Gap Property.
  \end{itemize}
\end{theorem}

We will give the proof of Hulanicki-Reiter's theorem  in the much wider  context of co-amenable
actions of groups on measure spaces in Section~\ref{S:Actions-Coam}
(see Remark~\ref{Remark-AmenablePairs}).

 F\o lner's theorem \cite{Foelner} is a refinement of the Hulanicki-Reiter Theorem. 
When $G$ is amenable, we  find  
a sequence of functions  $f_n\in L^2(G)$ with 
$\Vert f_n\Vert_2=1$ and $\lim_n\Vert \pi_G(g) f_n- f\Vert_2 =0$ for all 
$g\in G.$ 
One may ask whether the $f_n$'s can be chosen as normalized indicator functions
of Borel sets in $G.$ This is indeed the case.

\begin{theorem}
\label{Theo-Folner}
 \textbf{(F\o lner's theorem)}
 Let $G$ be an amenable locally compact group with left Haar measure $m.$
Then there exist so-called F\o lner sets: for every compact subset $Q$ of $G$ and
every $\varepsilon>0,$ there exists a Borel subset $U$ of $G$ 
with $0<m(U)<\infty$ such that
$$
\frac{m(xU\triangle U)}{m(U)}\leq \varepsilon \tout  x\in Q,
$$
where $\triangle$ denotes the symmetric difference.
\end{theorem}
We will generalize F\o lner's theorem to the context of co-amenable actions
 (Theorem~\ref{Folner-Coamenable}).
 
The following proposition on the relation between amenability and Property (T) 
is an obvious consequence of Hulanicki-Reiter's Theorem.
 \begin{proposition}
 \label{Pro-T-Ame}
 An amenable locally compact  group $G$ has Property (T) 
  if and only if $ G$ is compact.
 \end{proposition}
 
Amenability may be expressed in terms of averaging operators.
 Kesten proved  the following result in case
 $G$ is a discrete group and $\mu$ a symmetric probability measure
on $G$ (\cite{Kesten}). The general case is due to \cite{GuiDer}.
% see also \cite{Berg1}. 
 
 \begin{theorem}
 \label{Kesten}
 \textbf{(\cite{Kesten},  \cite{GuiDer})}
 Let $G$ be a locally compact group.
 and $\mu$ a strongly adapted probability measure on $G$. 
 The following statements are equivalent.
 \begin{enumerate}
 \item the group $G$ is amenable;
\item $\Vert \pi_G(\mu)\Vert =1.$
\end{enumerate}
 \end{theorem}
 
When $\mu$ is absolutely continuous, this is a consequence of  the Hulanicki-Reiter theorem
and the characterization of  the Spectral Gap Property for a unitary representation 
$\pi$ in terms of the averaging operator $\pi(\mu)$ from Proposition~\ref{SGP-Norm}.
 We will prove a  more general result in  the  context of co-amenable
actions  in Section~\ref{S:Actions-Coam}.

 \begin{remark}
\label{Exa-SpectralRadiusNorm}
\textbf{(Spectral radius versus operator norm)}
It may  happen that $\Vert\pi_{G}(\mu)\Vert=1$
for an adapted probability measure $\mu$ even when the group $G$
is non-amenable (of course, we then have  $r_{\spec}(\pi_G(\mu))<1$).
Indeed, let $G= \FF_2 $ be the free group on $2$ generators $a$ and $b$
and $\mu$ the uniform distribution on $\{a, b\}.$
Then $\Vert\pi_{G}(\mu)\Vert=1.$
We have 
\begin{align*}
\Vert\pi_{G}(\mu)\Vert&= \Vert\pi_{G}(\delta_{a^{-1}}\ast\mu)\Vert=
\Vert\pi_{G}(\nu)\Vert,
\end{align*}
where $\nu=\delta_{a^{-1}}\ast\mu$ is the uniform distribution on $\{e, a^{-1}b\}.$
Since the subgroup generated by $a^{-1}b$ is cyclic
and hence amenable, we have $\Vert\pi_{G}(\nu)\Vert=1.$
So,  $\Vert\pi_{G}(\mu)\Vert =1,$ although $\supp(\mu)$
generates $\FF_2 $ and $\FF_2 $ is not amenable.
\end{remark}

\section{Random walks and spectral radius of averaging operators}
\label{S:ApplicationSGP}
We introduce the spectral radius of a finitely generated group
(that is,  the spectral radius of the  simple random associated to a fixed generating set),
determine its exact value for free groups after Kesten and  relate
norms of averaging operators of Bernoulli actions  of groups
to  their spectral radius.
 
\subsection{Random walks on groups}
\label{ApplicationsSGP-RW}
Let $\Ga$ be a finitely generated group. Let $S$ be a finite generating set
of $\Ga$ with $S^{-1}=S.$
Let  $\G(\Ga, S)$ be the associated Cayley graph,
which is the graph defined as follows:
the vertex set  is $\Ga$ and  
$(x,y)\in \Ga\times \Ga$ is an edge if and only if $y=xs$ for some $s\in S.$

The simple random walk on $\G(\Ga, S)$ is the random walk $X_n$ defined as follows:
 every step consists in moving from a vertex $x$ to a neighbour $xs$ 
 with probability $1/{|S|}.$ 
The associated Markov operator 
 $$M=\pi_{\Gamma}(\mu_S)=\frac{1}{|S|}\sum_{s\in S}\pi_{\Gamma}(\delta_s)$$
acts on $\ell^2(\Gamma),$
where  $\mu_S$ is the uniform distribution on $S.$

Then
$$
\langle M^{n}\delta_e, \delta_e\rangle=\mu_S^{n}(e)
$$
is the probability 
${\mathbb P} (X_n=e| X_0= e)$ of return to the group unit at time $n,$
where $\mu_S^n$ denotes $n$-fold convolution.
Now, one can show (see Proposition~\ref{Pro-NormRegRep}) that 
$$
\lim_n \mu_S^{2n}(e)^{1/2n}=\Vert\pi_{\Gamma}(\mu_S)\Vert.
$$
The number 
$$\rho=\lim_n \mu_S^{2n}(e)^{1/2n}$$
is usually called the \emph{spectral radius} of the random walk
($\rho$ is indeed the spectral radius of the self-adjoint operator $M$).
So
$$
\rho= \lim_n {\mathbb P}( X_{2n}=e| X_0= e)^{1/2n}.
$$
In the case where $\Gamma$ is non-amenable,  $\rho<1$ and hence
${\mathbb P} (X_{2n}=e| X_0= e)$ decreases exponentially fast
as $n\to +\infty$.

\subsection{An example of the computation of the spectral radius}
\label{SS:KestenFree}
 Kesten determined in  \cite{Kesten-Random} the exact value  
for the spectral radius  of simple random walk on  a free group.

Let $\Ga= \FF_N $ be the free group on $N$ generators $a_1, \dots, a_N.$
Let $\mu$ be the uniform distribution on $\{a_1^{\pm 1}, \dots, a_N^{\pm 1}\}:$
$$\mu(a_i)= \mu(a_i^{-1})=\dfrac{1}{2N} \tout 1\leq i\leq N.$$
We claim that 
the spectral radius of the associated random
walk is
$$
\rho=\Vert\pi_{\Gamma}(\mu)\Vert=\frac{\sqrt{2N-1}}{N}.
$$
More generally, for $d\geq 2,$ let $T$ be the $d$-regular tree.
We consider the random walk on the vertices of $T$ with 
transition probability equal to $1/d$ to go from one vertex to one of its neighbours.
Let $M: \ell^2(T)\to \ell^2(T)$ be the associated Markov operator, defined by 
$$
M f(v)= \dfrac{1}{d} \sum_{w\sim v} f(w),
$$
where the sum is over the vertices $w$ which are neighbours of $v$.
(The Cayley graph of  $(\FF_N, \{a_1^{\pm 1}, \dots, a_N^{\pm 1}\})$
is the $2N$-regular tree, and $M$ can be identified with $\pi_{\FF_N}(\mu)$).
We  are going to show that 
$$\Vert M\Vert=\frac{2\sqrt{d-1}}{d}.$$

Apart from Kesten's original one, there are several proofs of this
formula  (see for instance Proposition 4.5.2 in \cite{Lubotzky}). 
We will follow a short  argument from \cite{Friedman}.

Since $T$ is normal and since $|\langle Mf, f\rangle|\leq \langle M|f|, |f|\rangle$ and  $\Vert f\Vert_2=\Vert |f|\Vert_2$ for $f\in\ell^2(T),$ we have 
$$\Vert M \Vert = \sup\left\{ \langle Mf, f\rangle \, : \, f\in \ell^2(T), f\geq 0,  \Vert f\Vert_2=1\right\}.$$

 Fix an origin $o\in T$ and denote by $\delta(v)$ the  graph distance of
 a vertex $v$ to $o.$
Observe that every vertex $v\neq 0$ has exactly $d-1$ neighbours $w$ with $\delta(w)= \delta(v)+1$
and one neighbour $w=w(v)$ with  $\delta(w)= \delta(v)-1$ and that $o$ has  $d$ neighbours $w$ all with
$\delta(w)=1.$

Let $f\in \ell^2(T)$ with $f\geq 0$ and $\Vert f\Vert_2=1.$ Then 
\begin{align*}
d\langle Mf, f\rangle 
&=\sum_{v\in T}\left( \sum_{w\sim v} f(v) f(w)\right)\\
&=\sum_{w:\delta(w)=1} f(o) f(w)+ 
\sum_{v\neq o}\left( \sum_{w: \delta(w)=\delta(v)+1} f(v)f(w)\right)
 +\sum_{v\neq  o} f(v) f(w(v))\\
 &= 2 \sum_{v\in T}\left( \sum_{w: \delta(w)=\delta(v)+1} f(v)f(w)\right)
 \end{align*}
Using the estimate
$$
f(v)f(w)\leq \dfrac{1}{2} \left(\dfrac{1}{\sqrt {d-1}} f(v)^2 + \sqrt{d-1} f(w)^2 \right),
$$
we obtain
\begin{align*}
d\langle Mf, f\rangle &\leq \sum_{v\neq o}  \left(\dfrac{d-1}{\sqrt{d-1}} f(v)^2 + \sqrt{d-1}f(v)^2\right)+
 \dfrac{d}{\sqrt {d-1}}f(o)^2 \\
 &\leq \sqrt{d-1} \sum_{v\in T} 2f(v)^2 = 2\sqrt{d-1}.
 \end{align*}
 This proves that $\Vert M\Vert \leq  \dfrac{2\sqrt{d-1}}{d}.$
 
To establish the lower bound, take an increasing sequence of real numbers  $\lambda_n<1$ with $\lim_n\lambda_n=1.$
Let $f_n$ be the radial function on $T$ defined by 
$f_n(v)= \left(\dfrac{\lambda}{\sqrt{d-1}}\right)^m$ if $\delta(v)=m.$
Then 
$$\Vert f_n\Vert_2^2=1+ \sum_{m\geq 1} \left(d(d-1)^{m-1} \left(\dfrac{\lambda_n}{\sqrt{d-1}}\right)^{2m}\right)= 1+ \dfrac{d}{d-1} \left(\dfrac{\lambda_n^2}{1-\lambda_n^2}\right)
$$
and
\begin{align*}
\langle Mf_n, f_n\rangle&= \dfrac{2}{d}\sum_{m\geq 0} \left( d(d-1)^{m-1} (d-1)\left(\dfrac{\lambda_n}{\sqrt{d-1}}\right)^{m+1}\left(\dfrac{\lambda_n}{\sqrt{d-1}}\right)^{m}\right)\\
&=\dfrac{2}{\sqrt{d-1}} \left(\dfrac{\lambda_n}{1-\lambda_n^2}\right).
\end{align*}
So, $\lim_n\dfrac{\langle Mf_n, f_n\rangle}{\Vert f_n\Vert_2^2} = \dfrac{2\sqrt{d-1}}{d}$
 and this proves the claim.

\begin{remark}
(i) With a little more effort, one can show that the spectrum of the operator $M$ as above is 
the whole interval $\left[-\dfrac{2\sqrt{d-1}}{d},\dfrac{2\sqrt{d-1}}{d}\right]$ (see \cite{Friedman} or \cite{Kesten-Random}).

\n
(ii) For $N\geq 2,$ let $\Ga$ be a group generated  by $N$ elements
$a_1,\ldots,a_N.$ For the spectral radius $\rho$
of the  random
walk on $\Ga$ defined by 
 the uniform distribution on $\{a_1^{\pm 1}, \dots, a_N^{\pm 1}\},$
one has $\rho \geq \frac{\sqrt{2N-1}}{N}.$
Indeed,  $\Ga$ is  a quotient of $\FF_{N}$  and the claim
follows  from Proposition~\ref{Pro-Shalom} below.
Kesten (see \cite{Kesten-Random}) proved that, if 
one has equality  $\rho= \frac{\sqrt{2N-1}}{N},$
then  $\Gamma$ is the free group on $a_1,\ldots,a_N$.

\n
(iii) Given a specific group $\Ga$ generated by a finite set $S$,
it is usually difficult  to compute or even to find bounds for the 
spectral radius  of the corresponding random walk.
For a recent result in the case where $\Ga$ is a surface group, 
see \cite{Gouezel}. 
The monography \cite{Woess} provides a comprehensive overview 
on results about  random walks on infinite groups
as well as on infinite graphs.
\end{remark}

\subsection{The norm of averaging operators for Bernoulli actions}
\label{Ex-Bernoulli}

Let $\Ga$ be an infinite countable group and $\Ga\curvearrowright X$ its Bernoulli action on
$X= \{0, 1\}^\Ga$ (see Section~\ref{Exa-Actions}).
Let $\mu$ be a symmetric and adapted probability measure on $\Ga$.
We claim that  
$$\Vert\pi_X(\mu)\Vert = \Vert\pi_{\Gamma}(\mu)\Vert.$$
In particular, it will follow from the Hulanicki-Reiter theorem that $\Ga\curvearrowright X$ has the
Spectral Gap Property if and only if $\Ga$ is not amenable.

To prove the claim, view $X$ as the  compact abelian group
$X= \prod_{\ga\in \Ga} \ZZ/2\ZZ.$ The dual
group is the discrete group $\widehat{X}= \bigoplus_{\ga\in \Ga} \ZZ/ 2\ZZ$.
 By Fourier transform, we have
$L^2_0(X)\cong \ell^2(\widehat{X}\setminus\{0\})$ as $\Ga$-representations.
Let $\Omega \subset \widehat{X}\setminus\{0\}$ be a set of representatives for
the $\Ga$-orbits in $\widehat{X}\setminus\{0\}.$
Then $\ell^2(\widehat{X}\setminus\{0\})$ decomposes as a direct sum
of $\Ga$-invariant subspaces:
$$
\ell^2(\widehat{X}\setminus\{0\})= \bigoplus_{x\in \Omega} \ell^2(\Ga x)
\cong \bigoplus_{x\in \Omega} \ell^2(\Ga /\Ga_x),
$$
where $\Ga_x$ is the stabilizer of $x.$ 
It is clear that
$\Ga_x$ is finite for every $x.$ Hence, we have $\ell^2(\Ga/\Ga_x)\subset \ell^2(\Ga)$
in an obvious sense. So,
$$\ell^2(\widehat{X}\setminus\{0\})\subset \bigoplus_{x\in \Omega} \ell^2(\Ga).$$
Therefore, $\Vert\pi_X(\mu)\Vert \leq \Vert\pi_{\Gamma}(\mu)\Vert.$
On the other hand,  we have, for every $x\in \Omega$ 
$$\Vert\pi_{\Gamma}(\mu)\Vert\leq  \Vert\pi_X|_{\ell^2(\Ga /\Ga_x)}(\mu)\Vert,$$
as follows from  Proposition~\ref{Pro-Shalom} below.

\section{Actions on homogeneous spaces with infinite measure}
\label{S:Actions-Coam}
In this section, we introduce and study a class of actions which we call co-amenable.
As  for amenability of groups,  these actions admit several characterizations,
the most notable one in our context being the absence of the Spectral Gap Property.

\subsection{Co-amenable actions}
\label{S:AmenableAct}
A useful generalization of amenable groups is the notion
of amenable homogeneous spaces in the sense of Eymard (\cite{Eymard}).
One can extend this notion to  actions on 
measure spaces as follows.
Let $G\curvearrowright X$ be an action of 
the separable locally compact group
$G$ on the measure space $(X,m),$ 
where as always $m$ is 
a $\sigma$-finite quasi-invariant measure.

\begin{definition} 
\label{Def-CoAmAction}
We say that the action of $G$ on $X$ 
is \emph{co-amenable} if 
there exists a $G$-invariant mean 
on $L^\infty(X,m),$ that is,
a positive linear functional $M$
on  $L^\infty(X,m)$ such that
that $M({\Un}_X)=1$ and $M({}_g\vfi) =M(\vfi)$
for all $g$ in $G$ and $\vfi$ in $L^\infty(X,m),$
where ${}_g\vfi(x)=\vfi (g^{-1}x).$
\end{definition}

\begin{remark}
\label{Remark-Co-Amenable}
(i) Consider the action of the locally compact group $G$ by left translation on $(G,m),$
where $m$ is a left Haar measure. This action is co-amenale
if and only if $G$ is amenable.

\n
(ii) Co-amenable actions $G\curvearrowright X$ 
were first considered by Greenleaf (\cite{Greenleaf}).
There were intensively studied by Eymard in the case 
of the action of $G$ on  a homogeneous space $X=G/H$ in \cite{Eymard}, 
where such an  action (or  space) is called  amenable. 
We prefer to call them  co-amenable in order to avoid 
confusion with the well-established notion of  an amenable action $G\curvearrowright X$ due to Zimmer (\cite{Zimmer}),
which in the case $X=G/H$ corresponds to the amenability of $H$;
a unification of both notions by means of an appropriate definition of amenable actions on pairs
of measure spaces is given in \cite{Zimmer2}.
For a further extension of these notions to the context on non-commutative measure spaces
(that is, to the context of von Neumann algebras), see \cite{Claire2}.

\n
(iii)  If $X$ is  a locally compact space,  one may define, as in \cite{Greenleaf}
or \cite{GuivAsym}, co-amenability  of the action $G\curvearrowright X$ 
through the existence of a $G$-invariant mean on the space $C^b(X)$ of 
continuous bounded functions on $X;$ as Example~\ref{CounterExa-AmenableSpace}
below shows,
this is in general  a weaker condition than co-amenability of the action of
$G$ on the measure space $(X,m)$, even in the case 
of a homogeneous space $X=\widetilde{G}/H$ 
(for a group $\widetilde{G}$ containing $G$), where a   natural quasi-invariant measure (class) $m$  
is given.

  \end{remark}

  Let  $G\curvearrowright (X,m)$ be an action 
 of $G$ on the measure space $m.$ 
  Given a measure space $(X,m),$ a mean $M$ on $L^\infty(X,m)$ defines 
  a \emph{finitely additive} probability measure $\mu_M$ on  the measurable subsets $A$ of $X,$ given
  by  $\mu_M(A)= M(\Un_A).$ 
  Such a finitely additive probability measure is absolutely continuous with respect to 
  $m,$ in the sense that if $m(A)=0$ then $\mu_M(A)=0.$
  Conversely,  a  finitely additive   probability measure $\mu$ on $X$ which is absolutely continuous with respect to 
  $m$ defines a unique mean $M$  on $L^\infty(X,m)$,
  given by $M(\vfi)=\sum_{i=1}^m \alpha_i \mu({A_i})$,
  if $\vfi=\sum_{i=1}^m \alpha_i \Un_{A_i}$
is a linear combination of indicator functions of measurable subsets $A_i$
of $X.$

 \begin{remark} 
    \label{Lem-AmCompactSpace} 
 (i)    Let $X$ be a \emph{compact} space.
    Then every  mean $M$ on  $L^\infty(X,m)$
  defines in a  natural way a genuine (that is,  $\sigma$-additive) probability measure
  on  $X.$
  Indeed,  let $\Phi: C(X)\to L^\infty (X,m)$ be the obvious mapping.
  Then $M\circ \Phi$ is a positive linear functional on $C(X)$
 with $M\circ \Phi(\Un_X)=1.$ By Riesz
 representation theorem, there exists a probability measure $\mu$ on $X$
such that 
$$M\circ \Phi(f) = \int_{X} f(x) d\mu(x) \tout f\in C(X).$$

Observe that, if a locally compact 
group $G$ acts on $(X,m)$ and if $M$ is $G$-invariant,
then the associated probability measure $\mu$ is $G$-invariant.

(ii) Recall that the space of probability 
measures $\P(X)$ on a compact space $X$ is compact
for the weak topology $\sigma(\P(X)), C(X));$ this is 
no longer true when $X$ is, say, a non compact
locally compact space.
However, for any measure space $(X,m)$,
 the space $\M(X)$ of means on $L^\infty (X,m)$
has an equally useful compactness property:
$\M(X)$ is compact
for the weak*-topology $\sigma(L^\infty (X,m), L^1 (X,m)).$
\end{remark}

  Means share with probability measures
  the property that they may be pushed forward through measurable mappings. 
  
  Let $(X,m)$ and $(Y,m')$ be $\sigma$-finite measure spaces
  and $\theta: X\to Y$ a measurable mapping
  such that  $m'$ is absolutely continuous with respect to
 the push-forward measure  $\theta_*(m).$
 Then 
  $$\Phi: L^\infty(Y,m')\to L^\infty(X, m), \quad \vfi \mapsto \vfi \circ \theta$$ 
  is a well-defined linear mapping which preserves positivity and maps $\Un_Y$ to $\Un_X.$
  So, if $M$ is a mean on $L^\infty(X, m),$ then  $\theta_*(M):= M\circ \Phi$ is 
 a mean on $L^\infty(Y, m'),$ which we call the \emph{push-forward mean} or the image
 of $M$ through $\theta.$

  An immediate  consequence of the consideration of push-forward
  means is the following useful fact.
  
   \begin{corollary}
   \label{Cor-AmenFactor}
   Let $G\curvearrowright X$ and $G\curvearrowright Y$ be   actions of the 
locally compact group $G$ on  
 measure spaces $(X,m)$ and $(Y,m')$. 
 Assume that there exists a measurable mapping $\theta: X\to Y$ which  intertwines the respective $G$-actions
 and  such that  $m'$ is absolutely continuous with respect to $\theta_*(m).$ 
 If  $G\curvearrowright X$ is co-amenable, then $G\curvearrowright Y$  is co-amenable.$\bsq$
 \end{corollary}
 \proof
 If $M$ is a $G$-invariant mean on $L^\infty(X, m),$ then  
its push-forward $\theta_*(M)$ is 
a $G$-invariant mean on $L^\infty(Y, m').$ $\bsq$

Given an action $G\curvearrowright (X,m)$,
it will be convenient in the sequel to consider $G$-invariant 
means  on  the subspace $L^\infty(X)_{G,u}$ 
of $G$-continuous functions   in $L^\infty(X),$ that is, the space 
$L^\infty(X)_{G,u}$ of all  $\vfi\in L^\infty(X)$ for which the mapping
$$G\to L^\infty (X), \qquad g\mapsto {}_g\vfi$$
 is norm-continuous. 

The space $L^\infty(X)_{G,u}$ is a large subspace, obtained by averaging
 functions from $L^\infty(X,m).$
 Indeed, it is easy to see 
that $f\ast \vfi \in L^\infty (X,m)_{G,u}$ for every $f\in L^1(G,\lambda)$ and $\vfi \in L^\infty(X),$
where $f\ast \vfi(x)= \int_G f(g)  {}_g\vfi(x) d\lambda(g)$ and $\lambda$ is Haar left measure on $G.$
(In fact, using Cohen's factorization theorem \cite{HewRos}, one can show that 
$L^\infty(X)_{G,u}=\{f\ast \vfi\ : \ f\in L^1(G,\lambda), \vfi \in L^\infty(X)\}.$)
It is obvious that $L^\infty(X)_{G,u}$ is a  $G$-invariant closed subspace of $ L^\infty(X)$
containing the constant functions.

Let $L^1(G)_{1,+}$ denote the
convex set of all $f\in L^1(G,\lambda)$ with $f\geq 0$ and $\Vert f\Vert_1=1.$
Observe that $L^1(G)_{1,+}$ is closed under convolution.

\begin{lemma}
\label{Lem-UnifContMean}
The following properties are equivalent.
\begin{itemize}
\item[(i)] There exists a $G$-invariant mean on $L^\infty(X).$ 
\item [(ii)] There exists a  $G$-invariant mean on $L^\infty(X)_{G,u}.$
\item[(iii)] There exists a \emph{topologically invariant}  mean on $L^\infty(X),$
that is, a mean $M$ such that  $M(f\ast \vfi)= M(\vfi)$ for 
all $\vfi\in L^\infty(X)$ and $f\in L^1(G)_{1,+}.$
\end{itemize}
\end{lemma}
\proof
The implication $(i)\Rightarrow (ii)$ is trivial.
To show that $(ii)\Rightarrow (iii),$  let $M$ be a $G$-invariant 
mean on $L^\infty(X)_{G,u}.$ Since
the mapping $G\to L^\infty(X)_{G,u}, \ g\mapsto {}_{g} \vfi$ is norm continuous for $\vfi\in L^\infty(X)_{G,u},$
one has 
$$
M(f\ast \vfi)=M(\vfi) \tout f\in L^1(G)_{1,+},\  \vfi\in L^\infty(X)_{G,u}.
$$

Let $(f_n)_n$ be an approximate identity in $L^1(G)_{1,+}$ for $L^1(G).$
Then, for each $\vfi\in L^{\infty}(X) $ and $f\in L^1(G)_{1,+},$ we have
$ \lim_n\Vert f\ast f_n \ast\vfi-f\ast \vfi\Vert_\infty=0$
and, hence, 
$$M(f\ast \vfi)=\lim_n M(f\ast f_n \ast\vfi)=\lim_n M( f_n \ast\vfi).$$
This shows that $M(f\ast \vfi)=M(f'\ast \vfi)$ for all  $f,f'\in L^1(G)_{1,+}$
and all $\vfi\in L^{\infty}(X) .$

Fix any $f_0\in L^1(G)_{1,+},$  and define a mean $\widetilde M$ 
on $L^{\infty}(X)$ by 
$$
\widetilde M (\vfi)=M(f_0\ast \vfi) \tout \vfi\in L^{\infty}(X).
$$ 
Then $\widetilde M$ is topologically invariant, since, for $f\in L^1(G)_{1,+}$ and $\vfi\in L^\infty(X),$ we have
\begin{align*}
\widetilde M  (f\ast\vfi)&=M(f_0\ast f\ast \vfi) = M(f_0 \ast \vfi) = \widetilde M(\vfi).
\end{align*}

To show that $(iii)\Rightarrow (i),$ let $M$ be a topologically invariant mean on $L^\infty(X).$
Then $M$ is $G$-invariant. Indeed,  fix
$f\in L^1(G)_{1,+}$. Then, for  
$\vfi\in L^\infty(X)$ and $g\in G$, we have  $f\ast {}_g\vfi= f_g \ast \vfi$ and hence
\begin{align*}
M({}_g\vfi)&=M(f\ast {}_g\vfi) = M(f_g \ast \vfi) =M(\vfi),
\end{align*}
for 
$f_g \in  L^1(G)_{1,+}$  defined by   $f_g(h)= \Delta(g^{-1}) f(h g^{-1}),$
where $\Delta$ is the modular function of $G.\bsq$

Recall that the   unitary representation $\pi_X$  of $G$ associated to the $G$-action on $X$
is defined  on $L^2(X,m)$ by
$$\pi_{X}(g) \xi(x)=\sqrt{\dfrac{dm(g^{-1}x)}{dm(x)}}\xi(g^{-1}x),
\qquad g\in G,\ x\in X,
\ \xi\in L^2(X).$$
(Observe that $G\curvearrowright X$ is trivially co-amenable if $m$ is a $G$-invariant
probability measure; so, our interest is in the case where $m$
is either infinite or not invariant.)

We will need the following general  fact concerning the spectral radius  of a
convolution operator acting on $L^2(X,m).$

\begin{lemma}
\label{Lem-SpctralRadius}
Let  $\mu$ be a probability
measure on $G$ with $r_\spec(\pi_X(\mu))=1$.
Then $1$ belongs to the spectrum $\sigma (\pi_X(\mu))$  of $\pi_X(\mu);$
more precisely, $1$ is an approximate eigenvalue of $\pi_X(\mu).$
\end{lemma}
\proof
Set $T:=\pi_X(\mu).$
Since  $r_\spec(T)=1,$ there exists $\la\in \sigma (T)$  
with $|\la|=1.$ 
 
 We claim that  $\la$ is an approximate eigenvalue of $T.$ Indeed, otherwise,
  $\Im(T-\la I)$ would be  a proper closed subspace of $L^2(X)$. So, we would have
 $\ker (T^*-\overline{\la}I) \neq \{0\}.$  However, since $T$ is a
 contraction, the equality case of Cauchy-Schwarz inequality shows that  $\ker (T^*-\overline{\la}I)= \ker (T-{\la}I)$.
 Hence, $\la$ would be an eigenvalue of $T$ and this would be a contradiction.
 
 So, there exists a sequence  $(\xi_n)_n$ in $L^2(X)$ with $\Vert\xi_n\Vert=1$ such that
$\lim_n \Vert T\xi_n-\la \xi_n \Vert =0$ or, equivalently,
$$
\lim_n\int_G  \langle \pi_X(g) \xi_n,\xi_n\rangle d\mu(g)=\lim_n \langle T\xi_n,\xi_n\rangle=\la.
$$
In particular, we have 
$$
\lim_n\left|\int_G \langle\pi_X(g) \xi_n,\xi_n\rangle d\mu(g)\right|=1.
$$
Since
\begin{align*}
1&=\int_G \Vert\pi_X(g)\xi_n\Vert\Vert \xi_n\Vert d\mu(g) \geq\int_G  \langle \pi_X(g)|\xi_n|,|\xi_n|\rangle d\mu(g)\\
&\geq\int_G  |\langle \pi_X(g) \xi_n,\xi_n\rangle| d\mu(g)\geq \left|\int_G \langle \pi_X(g)\xi_n,\xi_n\rangle d\mu(g)\right|,
\end{align*}
it follows that 
$$
\lim_n \langle T|\xi_n|,|\xi_n|\rangle=
\lim_n\int_G \langle \pi_X(g)|\xi_n|,|\xi_n|\rangle d \mu(g)= 1,
$$
that is, $\lim_n \Vert T|\xi_n|- |\xi_n|\Vert=0.$
Hence, $1$ is  an approximate eigenvalue of $T.$
$\,\bsq$

The following result, which is the main result of this section,
 was obtained by Guivarc'h in \cite[Proposition~1]{GuivAsym}.
The equivalence $(i)\Leftrightarrow (ii)$ is due to Eymard (\cite{Eymard})
in the case of  an action $G\curvearrowright G/H$ for a closed
subgroup $H$ of $G$. The equivalence 
$(ii) \Leftrightarrow (iv)$ (or $(ii) \Leftrightarrow (v)$) has been
extended to a class of more general unitary representations  in \cite{BeGui06}.

\begin{theorem}\textbf{(\cite{GuivAsym})} 
\label{Theo-Guiv}
Let $G\curvearrowright X$ be an action of 
the separable locally compact group
$G$ on the measure space $(X,m),$ where $m$
is a $\sigma$-finite quasi-invariant measure on $X.$ The following properties are equivalent:

\n
(i)  The action $G\curvearrowright X$ is co-amenable;

\n 
(ii) the representation $\pi_X$ of $G$ does not have the Spectral Gap Property;

\n
(iii) we have  $r_{\spec} (\pi_X(\mu))=1$ for  every  probability  measure $\mu$ on $G;$

\n
(iv)   we have  $\Vert \pi_X(\mu)\Vert=1$ for  some   strongly adapted probability  measure $\mu$ 
on $G;$ 

\n
(v) we have  $r_{\spec} (\pi_X(\mu))=1$ for  some  adapted probability  measure $\mu$  on $G.$ 

\end{theorem}

\proof 
$(i)\Rightarrow (ii):$ 
The set $\M$ of all means on $L^{\infty}(X)$ 
is a weak* closed (and hence compact) convex subset of the
unit ball of $L^{\infty}(X)^*.$
We can  view  the set $L^1(X)_{1,+}$ of densities  as a  subset of  $\M,$
since every $\xi$ in $L^1(X)_{1,+}$  defines 
an element in $\M,$ via integration against $\xi.$ Hahn-Banach's theorem
shows that  $L^1(X)_{1,+}$ is weak* dense in $\M.$

The group $G$ acts by isometries
 on $L^1(X,m),$ through the formula
 $$
 \pi_X(g) \xi(x) = \dfrac{dm(g^{-1}x)}{dm(x)}\xi(g^{-1}x) \tout g\in G, \ \xi\in L^1(X,m),\ x\in X,
 $$
and so $L^1(X)$ is a continuous $L^1(G)$-module, via 
$$f\ast \xi (x):=\int_G f(g) (\pi_X(g) \xi)(x) d\lambda(g) \tout f\in L^1(G), \ \xi\in L^1(X).$$
One checks that, for $f\in L^1(G), \xi \in L^1(X)$ and $\vfi \in L^\infty(X)$, 
one has 
 $$
 \int_X (f\ast \xi)(x) \vfi(x) dm(x)= \int_X  \xi(x) (\check{f}\ast  \vfi)(x) dm(x),
 $$
 where $\check{f}\in L^1(G)$ is defined by $\check{f}(g)=\Delta(g^{-1}) f(g^{-1}).$

Let $\L$ be  a countable dense subset of $L^1(G)_{1,+}$
for the $L^1$-norm.
For each $f\in \L,$ take a copy of $L^1(X)$ 
and consider the product space $
E=\prod_{f\in \L} L^1(X),$ 
 with the product of the norm topologies.
Then $E$ is a locally convex space, and the weak
topology on $E$ is the product of the
weak topologies. Consider the convex set 
$$
\Sigma=\{(f\ast \xi-\xi)_{f\in \L} :\ \xi\in L^1(X)_{1,+}\}\subset E.
$$ 

Since $G\curvearrowright X$ is co-amenable, there exists a  $G$-invariant mean $M$ on $L^{\infty}(X),$
which is even topologically invariant (see Lemma~\ref{Lem-UnifContMean}). 
 Take a net $(\xi_i)_i$ in $ L^1(X)_{1,+}$ with
$$M(\vfi) =\lim_i\int_X \xi_i(x)\vfi(x) dm(x) \tout \vfi\in L^{\infty}(X).$$
Then, since $M$  is topologically invariant, we have
 $\lim_i f\ast \xi_i- \xi_i=0$ in the weak topology of $L^1(X),$
 for  every $f\in L^1(G)_{1,+}.$
 So, $0$ belongs to the  closure 
 of $\Sigma$ for the product of  the weak topologies and hence  for the product of the norm-topologies,
 since $\Sigma$ is convex.
 Therefore,   we can find
a sequence  $(\xi_n)_n$
in $L^1(X)_{1,+}$ such that,
for every $f\in\L$, we have
$$\lim_{n}\Vert f\ast \xi_n-\xi_n\Vert_1= 0.$$
It follows that 
$\lim_n\Vert f\ast \xi_n-\xi_n \Vert_1=0$ for every $f\in L^1(G)_{1,+},$ 
by density of $\L$.
Since 
$$
\Vert \pi_X(g)(f\ast \xi_n)- f\ast \xi_n \Vert_1\leq \Vert ({}_{g}f)\ast \xi_n-  \xi_n \Vert_1 +\Vert f\ast \xi_n-  \xi_n\Vert_1,
$$
we  have therefore 
$$\lim_n\Vert \pi_X(g)(f\ast \xi_n)- f\ast \xi_n \Vert_1=0 \tout g\in G.$$

Choose $f\in \L$ and set $\eta_n = \sqrt{f\ast \xi_n}$. 
Then $\Vert \eta_n\Vert_2^2= \Vert f\ast\xi_n\Vert_1=1;$
moreover, for every $g\in G,$ we have 
$$\lim_n\Vert  \pi_X(g)\eta_n-\eta_n \Vert_2=0,$$ since
$$
\Vert \pi_X(g)\eta_n-\eta_n \Vert_2^2\leq \Vert \pi_X(g)(f\ast \xi_n)- (f\ast\xi_n) \Vert_1,
$$
using the elementary inequality $|a-b|^2\leq |a^2-b^2|$ for all 
non-negative real numbers $a$ and $b.$
It follows then that 
$\lim_n\Vert  \pi_X(g)\eta_n-\eta_n \Vert_2=0$ uniformly on compact subsets of $G$
(see Remark~\ref{Rem-SGP}.i).
This proves $(ii).$

\n 
$(ii)\Rightarrow (iii):$ 
Assume that $\pi_X$ has does not have Spectral Gap Property.
 So, there exists a sequence $\eta_n \in L^2(X)$ with
 $\Vert \eta_n\Vert_2=1$ and $\lim_n\Vert \pi_X(g) \eta_n-\eta_n \Vert_2=0$
 for all $g\in G.$
 Set $\xi_n = |\eta_n|^2.$ Then $\xi_n \in L^1(X)_{1, +}$ and
 $\lim_n\Vert \pi_X(g)\xi_n-\xi_n \Vert_1 =0$ for all $g\in G,$
 since
 $$
 \Vert \pi_X(g) \xi_n-\xi_n \Vert_1 \leq 2 \Vert  \pi_X(g) \eta_n-\eta_n \Vert_2,
 $$
 by Cauchy-Schwarz inequality.
 
 Let $\mu$ be any probability measure on $G.$
 Since
 $$
\Vert \pi_X(\mu) \eta_n -\eta_n\Vert_2 \leq \int_G\Vert \pi_X(g) \eta_n -\eta_n\Vert_2 d\mu(g),
$$
it follows from Lebesgue convergence theorem that $\lim_n\Vert \pi_X(\mu) \eta_n -\eta_n\Vert_2=0$.
Hence, $r_{\spec} (\pi_X(\mu))=1.$

  \n 
$(iii)\Rightarrow (iv)$ is obvious.

\n
$(iv)\Rightarrow (v)$ Let  $\mu$ be a strongly adapted probability
measure on $G$ with $\Vert\pi_X(\mu)\Vert=1.$
Then $\check{\mu}\ast \mu$ is adapted and $\pi_X(\check{\mu}\ast \mu)$ is a self-adjoint operator on $L^2(X)$
with  $\Vert\pi_X(\check{\mu}\ast \mu)\Vert=\Vert\pi_X(\mu)\Vert^2=1.$
Hence, $r_{\spec} (\pi_X(\check{\mu}\ast \mu))=1.$

  \n 
$(v)\Rightarrow (i):$
Let  $\mu$ be an adapted probability
measure on $G$ with  $r_{\spec} (\pi_X(\mu))=1.$

By Lemma~\ref{Lem-SpctralRadius}, 
$1$ belongs to  $\sigma (\pi_X(\mu))$ and is  an approximate eigenvalue.
Hence, there exists a sequence of unit vectors $\eta_n$ in $L^2(X)$ with 
$$\lim_n\Vert \pi_X(\mu) \eta_n -\eta_n\Vert_2=0.$$
 So,
$\lim_n \langle  \pi_X(\mu) \eta_n, \eta_n\rangle =1,$ that is, 
$$
\lim_n \int_G \langle  \pi_X(g) \eta_n, \eta_n \rangle d\mu(g)=1.
$$
It follows that there exists a subsequence, again denoted by $\eta_{n},$ such that
$$
\lim_n  \langle  \pi_X(g) \eta_n, \eta_n \rangle =1.
$$
for $\mu$-almost every $g\in G.$
Since
$$
\Vert \pi_X(g) |\eta_n| -|\eta_n|\Vert_2 \leq \Vert \pi_X(g) \eta_n -\eta_n\Vert_2,
$$
we can assume that $\eta_n \geq 0.$
Set 
$$
H=\{ g\in G \ : \ \lim_n  \langle  \pi_X(g) \eta_n, \eta_n \rangle =1\} =\{ g\in G \ : \ \lim_n \Vert \pi_X(g) \eta_n -\eta_n\Vert_2  =0\}.
$$
Then $H$ is a measurable subgroup of $G$ with $\mu(H)=1.$
Hence, $H$ contains  $\supp(\mu)$ and it follows that $H$ is dense in $G.$

Set $\xi_n=\sqrt{ \eta_n}\in L^1(X)_{1,+}.$ 
 Let now $M$ be a mean on $L^\infty (X)_{G,u}$ which is a limit
 of $(\xi_n)_n$ in the weak-*-topology.
 Then $M$ is $H$-invariant.  Since, $g\mapsto {}_g\vfi$ is 
 norm-continuous for $\vfi\in L^\infty (X)_{G,u},$ it follows that 
 $M$ is $G$-invariant. Hence, $G\curvearrowright X$ is co-amenable
 by Lemma~\ref{Lem-UnifContMean}. $\bsq$

 The following consequence of  the equivalence between (i),  (iii) and (v)
 in the previous theorem is worth mentioning.
 \begin{corollary}
 \label{Cor-Co-amenable}
 Let $G\curvearrowright X$ be an action of 
the separable locally compact group
$G$ on the measure space $(X,m).$ 
Let $H$ be a separable locally compact group,
$j:H\to G$ a continuous homomorphism 
and $H\curvearrowright X$ the corresponding action.

\n
(i)  If  $G\curvearrowright X$ is co-amenable, then $H\curvearrowright X$ is co-amenable.
In particular,  $H\curvearrowright X$ is co-amenable if $H$ is a closed
subgroup or a countable dense subgroup of $G.$

 \n
 (ii)  Assume that $G\curvearrowright X$ is not co-amenable and that $j(H)$ is dense in $G$.
 Then $H\curvearrowright X$ is not co-amenable.
  \end{corollary}

  \begin{remark}
\label{Remark-AmenablePairs}
(i) Considering the action $G\curvearrowright G$ given
by left translation, we see that  Theorem~\ref{Kesten} 
is a special case of  Theorem~\ref{Theo-Guiv}. 

\n
(ii)
For an extension of Theorem~\ref{Theo-Guiv} to amenable pair of actions
(including actions on von Neumann algebras), see \cite{Claire1} and \cite{Claire2}).
\end{remark}

\begin{example}
\label{CounterExa-AmenableSpace}
Identify the group  $G=SL_2(\RRR)$ with the subgroup $
\left( \begin{array}{ccc}
1&0&0\\
0&\ast &\ast\\
0&\ast &\ast\\
\end{array}\right)
$
of  $H=SL_3(\RRR);$
the standard action of $G$ on   $X=\RRR^3\setminus\{0\}$ 
(which is a homogeneous space of the form $H/L$),
preserves the Lebesgue measure $m$ on $X,$ and fixes the 
first unit vector $e_1\in \RRR^3.$
%\left( \begin{array}{c}
%1\\
%0\\
%0\\
%\end{array}\right).
%$
The Dirac measure $\delta_{e_1}$ is a  $G$-invariant mean on the space $C^b(X);$
however, there exists no $G$-invariant mean on $L^\infty(X,m).$
Indeed, otherwise, the corresponding unitary representation $\pi_X$ of $G$ on $L^\infty(X)=L^2(\RRR^3)$
would have almost invariant vectors (by Theorem~\ref{Theo-Guiv}).
On the other hand, upon neglecting a set of $m$-measure $0,$   the $G$-orbits in $X$ are the sets  
%$X_x=\{ xe_1+ye_2+z e_3\ \ : \ y,z\in \RRR\},$
$$X_x=\left\{ \left(\begin{array}{c}
x\\
y\\
z\\
\end{array}\right) \: |\ y,z\in \RRR\right\}$$ for $x\in \RRR,$
 which all are isomorphic to $G/N$ where $N$ is the subgroup of unipotent upper-triangular matrices. So, $\pi_X$ is 
weakly equivalent to  the quasi-regular representation  $\pi_{G/N}$ in $G/N.$ 
Since, $N$ is amenable, $\pi_{G/N}$ is weakly contained in  
the regular representation $\pi_G$ (see Theorem F.3.5 in \cite{BHV}).
Therefore, $\pi_G$
would have almost invariant vectors and this impossible, since $G=SL_2(\RRR)$ is not amenable.
\end{example}

\subsection{F\o lner sequences}
\label{Folner-Coam}
We now prove the existence of  F\o lner sequences in $X$
for co-amenable actions $G\curvearrowright (X,m),$  provided the measure $m$
is \emph{$G$-invariant}. More precisely,
the following extension of Theorem~\ref{Theo-Folner} holds
(see also \cite{Greenleaf} and \cite{Iozzi-Nevo}).
 
\begin{theorem}
\label{Folner-Coamenable}
 \textbf{(Existence of F\o lner  sequences)}
 Let $G$ be a separable locally compact group
 and let $G\curvearrowright X$ be a co-amenable action of 
$G$ on the measure space $(X,m).$
Assume that $m$ is $G$-invariant.
Then, for every compact subset $Q$ of $G$ and
every $\varepsilon>0,$ there exists a measurable subset $U$ of $X$ 
with $0<m(U)<\infty$ such that
$$
\frac{m(gU\triangle U)}{m(U)}\leq \varepsilon \qquad\tout  g\in Q.
$$
\end{theorem}
\proof 
We follow the proof  given in    \cite[Theorem G.5.1]{BHV} for 
the group case (Theorem~\ref{Theo-Folner}),
which carries  over \emph{mutatis mutandis}  in our situation.

Let $Q$ be a compact subset of $G$ containing $e,$ and  let $\varepsilon>0.$ 
Set $K: =Q^2.$ 
Since $\pi_X$ does not have the Spectral Gap Property,
we can find $\xi\in L^2(X)$ with $\Vert \xi\Vert_2=1$ such that 
$$
\sup_{g\in K}\Vert \pi_X(g) \xi-\xi\Vert_2 \leq \frac{\varepsilon |Q|}{2|K|},
$$
where we denote by $|A|$ the Haar measure of a subset $A$ of $G.$

Set $f= |\xi|^2.$ Then  $f\in L^1(X)_{1,+}$ and
$$
\sup_{g\in K}\Vert{}_{g^{-1}} f-f\Vert_1 \leq \frac{\varepsilon |Q|}{2|K|}.
$$
For $t\geq 0,$ let  $E_t:=\{y\in X\ :\ f(y)\geq t\}.$ 
Then 
$$
gE_t=\{ y\in X\ :\ {}_{g^{-1}} f(y)\geq t\}.
$$
By the   lemma below, we have
$$
\Vert{}_{g^{-1}} f-f\Vert_1=\int_{0}^{\infty}m(gE_t\triangle E_t)dt .
$$
for every $g\in G.$ Hence, for every $g\in K,$ we have
\begin{align*}
\int_{0}^{\infty}m(E_t)\left(\int_K\frac{m(gE_t\triangle E_t)}{m(E_t)}dg\right) dt
&=\int_K\Vert{}_{g^{-1}} f-f\Vert_1 dg \leq \frac{\varepsilon |Q|}{2}.
\end{align*}
Since $\int_{0}^{\infty}m(E_t)dt=\Vert f\Vert_1=1,$
 it follows that there
exists $t$ with $0<m(E_t)<\infty$ and such that 
$$
\int_K\frac{m(gE_t\triangle E_t)}{m(E_t)}dg\leq 
\frac{\varepsilon |Q|}{2}.
$$
Set 
$$
A=\{g\in K:\ \frac{m(gE_t\triangle E_t)}{m(E_t)}\leq\varepsilon\}.
$$
Then  $|K\setminus A|<|Q|/2 $ and we claim that $Q\subset AA^{-1}.$

To show this, let $g\in Q.$ Then $|gK\cap K|\geq |gQ|=|Q|$ and, hence,
\begin{align*}
|Q|\leq |gK\cap K|&\leq |gA\cap A|+|K\setminus A|+|g(K\setminus A)|=|gA\cap A|+2|K\setminus A|\\
&<|gA\cap A|+|Q|. 
\end{align*} 
Therefore, $|gA\cap A|>0,$ and this implies that $g\in AA^{-1},$ as claimed.

Now, for $g_1, g_2\in A,$ we have
$$
g_1g_2^{-1}E_t\triangle E_t\subset (g_1g_2^{-1}E_t\triangle g_1E_t)\cup (g_1E_t\triangle E_t),
$$
and hence
\begin{align*}
m(g_1g_2^{-1}E_t\triangle E_t)
&\leq m(g_2^{-1}E_t\triangle E_t)+m(g_1E_t\triangle E_t)\\
&=m(g_2 E_t\triangle E_t)+m(g_1 E_t\triangle E_t)\leq 2\varepsilon m(E_t).\bsq
\end{align*} 

The following formulas, which are  versions of the area and co-area
formulas from  Lemma~\ref{Lem2-Cheeger}, have a
similar   elementary proof (see Lemma G.5.2 in \cite{BHV}).
\begin{lemma}
\label{Lem-Folner}
Let $(X,m)$ be a measure space.
Let $f,f'$ be non-negative functions in $L^1(X).$
For every $t\geq 0,$ let
$E_t=\{x\in X:\ f(x)\geq t\}$ and $E'_t=\{x\in X:\ f'(x)\geq t\}.$
Then 
$$
\Vert f-f'\Vert_1=\int_{0}^{\infty} m(E_t\triangle E'_t)dt.
$$
In particular, $\Vert f\Vert_1= \int_{0}^{\infty}m(E_t)dt.$
\end{lemma}

\begin{remark}
\label{Rem-Foelner-Coamenable}
As noticed in \cite{Greenleaf}, the previous theorem may fail if  one drops the assumption that the measure $m$ on $X$
is $G$-invariant. Indeed, consider for example the action of $G=\ZZ$
on $X=\ZZ$ by translations, where $X$ is equipped with the measure defined by 
$m(\{n\})= 2^{|n|}$ for all $n\in \ZZ.$ It is easy to see that, for $Q=\{\pm 1\}$ and
for every finite set $U$ of $\ZZ,$ we have
$\sup_{g\in Q}\frac{m(gU\triangle U)}{m(U)}\geq 1.$
\end{remark}

\subsection{Norm of averaging operators under the regular representation}
\label{SS:RegularRep}
Let $G$ be a separable locally compact group,
and denote by $m$ a left Haar measure on $G.$
We consider the action of $G$ on itself by left translation.

The associated unitary representation $\pi_G$ is the 
left regular representation of $G$ on on $L^2(G,m).$
It is a remarkable fact that, given a probability measure $\mu$
on $G,$  the  norm of $\pi_G(\mu)$ under the regular representation
gives a lower bound for  $\pi_X(\mu)$ for \emph{every} action 
$G\curvearrowright X$ on a measure space $X$.
This fact, which may be viewed as a version of  Herz's  majoration principle from \cite{Herz},
 was proved by Shalom (see Lemma 2.3 in \cite{Shalom-AnnMath}).

We first  give a formula, due to Kesten 
in the discrete case and to Berg and Christensen in general,
  for the norm of the convolution operators defined by probability measures on $G.$

\begin{proposition}
\label{Pro-NormRegRep}
\textbf{(\cite{Kesten-Random}, \cite{Berg-Christensen})}
Let $\mu$ be a probability measure on the separable locally compact group $G.$
Then, for every compact neighbourhood $U$ of the group unit $e,$ we have
$$
\Vert \pi_G(\mu) \Vert =  \lim_{n\to +\infty}  \langle \pi_G(\check{\mu}\ast\mu)^{n}{\mathbf 1}_{U},{\mathbf 1}_{U}\rangle^{1/2n}= 
\lim_{n\to +\infty} \left((\check{\mu}\ast\mu)^{n}(U)\right)^{1/2n}.
$$
In particular, if $G$ is discrete, we have 
$$
\Vert \pi_G(\mu) \Vert = \lim_{n\to +\infty}(\check{\mu}\ast\mu)^{n} (e)^{1/2n} =  \lim_{n\to +\infty} \langle \pi_G(\check{\mu}\ast\mu)^{n}) \delta_e, \delta_e\rangle^{1/2n}.
$$
\end{proposition}
\proof
We give the proof of the formula  in the case where $G$ is discrete.

First, we observe that the sequence $(\check{\mu}\ast\mu)^{n} (e)^{1/2n}$ converges.
Indeed, set $a_n= (\check{\mu}\ast\mu)^{n} (e).$ Let $T$ be the square root of 
the positive selfadjoint operator $\pi_{G}(\check{\mu}\ast\mu).$ 
Then 
$$a_n= \Vert T^n \delta_e\Vert^2= \langle T^n \delta_e, T^n \delta_e\rangle $$ and, using Cauchy-Schwarz inequality,
we see that $a_{n}^2 \leq a_{n+1} a_{n-1}.$ Hence
$({a_{n+1}}/{a_n})_n$ is increasing and therefore convergent. It follows that $(a_n^{1/2n})_n$ converges (to the same limit).

Next, since $\pi_{G}(\check{\mu}\ast\mu)$ is selfadjoint, we have  
$$\Vert\pi_{G}(\check{\mu}\ast\mu)\Vert=
\lim_n \Vert\pi_{G}(\check{\mu}\ast\mu)^{n}\Vert^{1/n},$$
so that 
\begin{align*}
\Vert\pi_{G}(\mu)\Vert^2
=\sup_{f\in \CCC[G]}\left(\limsup_n
\langle\pi_{G}(\check{\mu}\ast\mu)^{n}f,f\rangle^{1/n}\right),
\end{align*}
where $\CCC[G]$ is the group algebra of $G,$ that is, the linear span of 
${\{\delta_x:\ x\in G\}.}$ Now
\begin{align*}
\limsup_n \langle\pi_{G}(\check{\mu}\ast\mu)^{n}\sum_{i=1}^k c_i\delta_{x_i}, 
\sum_{i=1}^k
c_i\delta_{x_i}\rangle^{1/n}\leq \max_{i=1}^k\left(\limsup_n \left\langle\pi_{G}(\check{\mu}\ast\mu)^{n}\delta_{x_i},
\delta_{x_i}\right\rangle^{1/n}\right),
\end{align*}
for all $x_1,\ldots, x_k\in G$ and $c_1,\ldots,c_k$ in $\CCC.$
Since
$$
\langle\pi_{G}(\check{\mu}\ast\mu)^{n}\delta_{x_i},
\delta_{x_i}\rangle= (\check{\mu}\ast\mu)^n(e),
$$
this proves the claim. $\bsq$

\begin{proposition}
\label{Pro-Shalom}
\textbf{(\cite{Shalom-AnnMath})}
Let $G\curvearrowright X$ be an action of 
the separable locally compact group
$G$ on the measure space $(X,m),$ where $m$
is a $\sigma$-finite quasi-invariant measure on $X.$ 
For  every probability measure $\mu$ on  $G,$ we have 
$$\Vert \pi_X(\mu) \Vert \geq  \Vert \pi_G(\mu) \Vert.$$
\end{proposition}

\proof
Set $\nu= \check{\mu}\ast \mu$ and 
let $\xi$ be a unit vector in $L^2(X)$ with $\xi\geq 0.$  The non-negative function
$$G\mapsto \RRR, \quad g\mapsto \langle \pi_X(g) \xi, \xi\rangle$$
is continuous and takes the value $1$ at $e.$ 
Hence, there exists a compact neighbourhood $U$ of $e$ such 
that
$$
\langle \pi_X(g) \xi, \xi\rangle \geq \dfrac{1}{2} \tout g\in U.
$$
For every $n\geq 1,$ we have
\begin{align*}
\langle \pi_X(\nu)^{n} \xi, \xi\rangle &=\langle \pi_X(\nu^{n}) \xi, \xi\rangle= \int_G \langle \pi_X(g) \xi, \xi\rangle d\nu^{n}(g)\\
&\geq \dfrac{1}{2} \nu^{n}(U).
\end{align*}
Since 
$$\Vert \pi_X(\mu)\Vert \geq  \langle \pi_X(\mu)^{n} \xi, \pi_X(\mu)^{n}\xi\rangle^{1/2n}= \langle \pi_X(\nu)^{n} \xi, \xi\rangle^{1/2n},
$$
it follows from  Proposition~\ref{Pro-NormRegRep} that 
$$
\Vert \pi_X(\mu)\Vert \geq  \Vert \pi_G(\mu) \Vert.\bsq
$$

As a consequence of the previous proposition and
the spectral characterization of amenability (Theorem~\ref{Kesten}),
 we see that  an action of  an amenable group
on a measure space \emph{never} has the Spectral Gap Property.

\begin{corollary}
\label{Pro-AmenableCoam}
Let $G$ be an amenable separable locally compact group. Then every action
$G\curvearrowright X$ on a measure space $(X,m)$
is co-amenable. $\bsq$
\end{corollary}
 Alternatively, the previous corollary follows
also from  the fixed point property of amenable groups
(Proposition~\ref{Pro-Amenable}.ii), applied to 
the convex set of means on $(X,m).$ 

In contrast to this, the co-amenable actions of  a  Kazhdan group 
 are the actions which are co-amenable  
for trivial reasons. More precisely, the following result holds.

\begin{corollary}
\label{Pro-KazhdanCoam}
Let $G$ be  locally compact group with Kazhdan's Property (T) and 
$G\curvearrowright X$ a co-amenable ergodic action on a measure space $(X,m).$
Then $m$ is equivalent to a $G$-invariant probability measure on $X.$
\end{corollary}
\proof
Indeed, the assumptions imply that there exists a function
$f\in L^1(X,m)_{1,+}$ which is $G$-invariant, that is,
$$\frac{dgm}{dm}(x)f(gx)= f(x) \tout g\in G, \ x\in X.$$
One checks that the density  $f(x)dm(x)$ is a $G$-invariant
measure on $X.$
Moreover, since the action the action is ergodic,  $f>0$
almost everywhere.$\bsq$

\subsection{Linear actions with the Spectral Gap Property}
\label{SS:RegularRep}
Given a separable locally compact group $G$ and an action $G\curvearrowright X$ on a  measure space $(X,m),$
there are only few general results ensuring  the Spectral Gap Property for this action.

We will consider linear actions on a vector spaces $V=\kk^d$ over an arbitrary 
 local field $\kk$.
  Let $G$ be a subgroup
 of $GL(V)=GL_d(\kk).$
 We have the following 
sufficient condition for the Spectral Gap Property
for the linear action $G\curvearrowright V\setminus \{0\},$ 
where $V$ is equipped with a
 translation invariant measure $m.$
This condition  involves the
induced action of $G$ on   the projective space  $\PP(V)$ of $V.$
Recall that $\PP(V)$ is compact for the quotient topology induced
from that of $V.$

\begin{lemma}
Let $G$ is a locally compact group which embeds continuously
in $PGL(V)=PGL_d(\kk).$
Assume that there exists no $G$-invariant probability measure
on $\PP(V)$. Then the action $G\curvearrowright V\setminus \{0\}$ 
has the Spectral Gap Property.
\end{lemma}
\proof
Indeed, assume, by contraposition, that the action $G\curvearrowright V\setminus \{0\}$ 
is co-amenable. 
The canonical
mapping $p: V\setminus \{0\} \to \PP(V)$ is $G$-equivariant 
and  the Lebesgue measure $m'$ on $\PP(V)$ is  absolutely continuous with respect to
$p_*(m)$. Hence, $G\curvearrowright  \PP(V)$ is co-amenable (see 
Corollary~\ref{Cor-AmenFactor}). 
So, there exists a $G$-invariant mean $M$ on $L^\infty(\PP(V), m').$
Since $ \PP(V)$ is a {compact} space, $M$ defines a 
a probability measure on $\PP(V)$,
which is $G$-invariant  (see Remark~\ref{Lem-AmCompactSpace}.i ).
$\bsq$

The following result is a consequence of Furstenberg's celebrated lemma
(see \cite{Furst} or  \cite[Corollary 3.2.2]{Zimmer})
on stabilizers of probability measures on projective spaces.
A  subgroup $G$ of $GL(V)$  is said to be \emph{totally irreducible}
if every finite index subgroup of $G$ acts irreducibly on $V.$

\begin{theorem}
\label{Theo-FurstenbergLemma}
Let $G$ be a totally irreducible subgroup in $SL(V).$  Assume that
$G$  is not relatively compact in $SL(V).$
Then $G\curvearrowright V\setminus \{0\}$ has
the Spectral Gap Property.
\end{theorem}
\proof
By the previous lemma, it suffices to show that there exists no $G$-invariant probability measure
on $\PP(V).$

Since  $G$ is not relatively compact, there exists a sequence $(g_n)_n$ in $G$
with $\lim_n\Vert g_n\Vert=\infty$ (for any norm on $\End (V)$). Set 
$u_n=\dfrac{g_n}{\Vert g_n\Vert}.$ Then $\Vert u_n\Vert=1$ and $\lim_n \det (u_n)= \lim_n {1}/{\Vert g_n\Vert^d}=0$.
Upon passing to a subsequence, we  can therefore assume that  $u=\lim_n u_n$ exists
in $\End (V)$ with $\Vert u\Vert=1$ and $\det (u)=0.$
So,  $\ker(u)$ and $u(V)$ are proper non-zero subspaces of $V.$
Denote by  $X_1$ and $X_2$ their  images  in $\PP(V)$.

Assume now, by contradiction, that there exists a  $G$-invariant probability measure $\nu$
on $\PP(V).$ 
Write 
$\nu=\nu_1 +\nu_2,$ where $\nu_1$ and 
$\nu_2$ are positive measures on  $\PP(V)$ with 
$$\nu_1(\PP(V)\setminus X_1)= 0 \qquad \text{and} \qquad \nu_2( X_1)= 0.$$
By compactness of the set of bounded positive measures on  $\PP(V)$, we can assume that
$\lim_n  g_n \nu_1 =m_1$ and $\lim_n  g_n \nu_2 =m_2$ exist. 
Then 
$$\nu= \lim_n g_n \nu= m_1+ m_2.$$
 Since $\nu_2$ is supported on $\PP(V)\setminus X_1$ and since 
$ \lim_n g_n x= u x \in X_2$ for every $x\in \PP(V)\setminus X_1,$
the measure $m_2$ is supported on $X_2.$ 
We can also assume that $X'= \lim_n g_n X_1$ exists, where
$X'$ is the image in $\PP(V)$  of a subspace of $V$ of the same dimension
as $\ker(u).$ 
So, $\nu$ is supported on the union $X' \cup X_2$  of two  proper projective subspaces.
This contradicts the fact that $G$ is totally irreducible.
Indeed, let $F$ be a projective subspace of minimal dimension 
with $\nu(F) >0.$ If  $\nu( gF \cap F)>0,$
then $gF=F$ by the minimality of dimension of $F$. 
Since $\nu$ is a probability measure, it follows that 
$\{ gF \ : \ g\in G\}$ is a finite set of projective subspaces.
Therefore, the proper projective subspace $F$ is fixed by a 
subgroup of $G$ of finite index. $\bsq$

\begin{remark}
\label{Rem-Furstenberg}
The proof of Theorem~\ref{Theo-FurstenbergLemma} 
shows that actually the following stronger statement holds
(see Lemma 3 in \cite{Furst}).
Let $G$ be  locally compact group which is
minimally almost periodic (that is, 
every continuous homomorphism of $G$ into
a compact group is trivial) and let  $\pi: G\mapsto GL(V)$  be a representation
of $G$ on a finite dimensional vector space $V$
over a local field. Consider the action of $G$ on $\PP (V)$
associated to $\pi.$
 If  $G$ preserves a probability  measure  $\nu$ on $\PP (V),$
 then $\nu$ is concentrated on the $G$-fixed points in $\PP(V).$
\end{remark}
\begin{remark}
\label{Rem-Shalom}
Shalom obtained in \cite{Shalom-TAMS}
the following  nice characterization of  co-amenable actions of 
an algebraic group $G$
over a local field $\kk$ on an arbitrary algebraic $\kk$-variety: there exists 
an amenable quotient of $G$ through which every such action
factorizes.
 \end{remark}

\subsection{Co-amenable subgroups of semisimple groups}
\label{SS:CoAmSS}
Given a  locally compact group $G,$
one may ask which closed subgroups $H$  are co-amenable in $G,$
that is, for which closed subgroups $H$  is $G\curvearrowright G/H$ co-amenable.
Observe that, by Corollary~\ref{Pro-KazhdanCoam},
 if $G$ has Kazhdan property, then $G\curvearrowright G/H$
is co-amenable  if and only $H$ has finite covolume in $G$
(that is, $G/H$ has a non trivial finite $G$-invariant Borel measure).
So, the only   simple real Lie groups interesting  for this question
are the isometry groups of real or complex hyperbolic spaces.
As we will see (Theorem~\ref{Theo-CoAm-CritExp}), in this case, the co-amenability of  a subgroup  can be read off
from its critical exponent.

It is a remarkable fact that co-amenable subgroups of
semi-simple algebraic groups are Zariski-dense, a result 
established in   \cite[Corollaire, p. 192]{GuivAsym} for the  case of real algebraic groups
and in  \cite[Corollary 1.7]{Shalom-TAMS} in general (this generalizes a result from \cite{Stuck};
see also \cite{Iozzi-Nevo}).

\begin{theorem}
\label{Theo-CoAmZariski}
\textbf{(\cite{GuivAsym},\cite{Shalom-TAMS})}
Let  $G=\mathbf G(\kk)$ be the group of $\kk$-rational points
of  a semi-simple connected algebraic group over the local field $\kk$
without compact factors.
Let $H$ be a closed co-amenable subgroup of $G$. Then 
$H$ is Zariski dense in $G.$
\end{theorem}

\proof
Let $L\subset G$ be the Zariski-closure  of $H$ in $G$. 
By a well-known result of Chevalley (see Theorem 5.1 in \cite{Borel}), there exist a 
$\kk$-rational representation $ G \to GL(V)$ 
over a finite dimensional $\kk$-vector space $V$ and a line
$\ell\subset V$ such that $L$ is   the stabilizer of $\ell$
in $G.$ The next proposition shows that $\ell$ is $G$-fixed and so $L=G.\bsq$

\begin{proposition}
\label{Pro-CoAmZariski}
Let $G$ be as in  Theorem~\ref{Theo-CoAmZariski}
and let $\pi:G\to GL(V)$  be a continuous representation
of $G$ on a finite dimensional vector space $V$
over $\kk.$
Let $H$ be a closed co-amenable subgroup of $G$ and let $\ell\subset V$
be a line which is $H$-invariant. Then $\ell$ is $G$-invariant (and hence
$G$ fixes every point in $\ell$).
\end{proposition}

\proof
We proceed by induction on $\dim V.$ The case
$\dim V=1$ being trivial, assume that $\dim V\geq 2$.

Let $x\in \PP(V)$ be the image
of $\ell$ in   $\PP(V)$.
The orbital mapping 
$$G\to \PP(V), \quad g\mapsto gx$$
induces  a $G$-equivariant 
continuous mapping ${f}: G/H\to \PP(V).$ 
If $M$ be an invariant mean on $L^\infty(G/H),$
then  $f_*(M)$ defines
a $G$-invariant probability measure $\nu$
on $\PP(V)$. 
Since $G$ is minimally almost periodic, $\nu$ is supported
by the $G$-fixed points in $\PP(V)$; see Remark~\ref{Rem-Furstenberg}
above. So, $G$ fixes a line $\ell'$ in $V.$ 

Consider now the representation $\overline{\pi}:G\to GL(V/\ell')$ associated to $\pi.$
Since $\dim (V/\ell') < \dim V,$ the induction hypothesis shows that 
the image of $\ell$ in $V/\ell'$ is $\overline{\pi}(G)$-invariant. 

Let $W= \ell +\ell'$. Then $W$ is a $G$-invariant subspace of $V.$
Moreover, the corresponding homomorphism $G\to GL(W)$ 
has its image contained in a solvable group. Hence,
$G$ acts trivially on $W$ and the claim is proved.$\bsq$

\begin{corollary}
\label{Cor-CoAM-Z}
Let $G$ be a connected simple  non compact real Lie  group
with trivial center and let $H$ be a  proper closed co-amenable subgroup of $G$.
Then $H$ is a  discrete subgroup of $G.$
\end{corollary} 
\proof
Let  $\mathfrak g$ be the Lie algebra of $G$  and 
${\mathrm Ad}: G\to GL(\mathfrak g)$ the adjoint representation.
Then $G$ can be identified with the connected component
(in the Hausdorff topology) of the real points of the   simple
algebraic group $\Aut(\mathfrak g)$ over $\RRR.$

The Lie subalgebra $\mathfrak h$ corresponding to $H$ is ${\mathrm Ad}(H)$-invariant and hence 
${\mathrm Ad}(G)$-invariant, since $H$ is Zariski dense in $G.$  
So, $\mathfrak h$ is an ideal in $\mathfrak g.$ Since $\mathfrak g$ is a simple
Lie algebra, it follows that either $\mathfrak h= \mathfrak  g$ 
or $\mathfrak h=\{0\}.$ 
The first case does not occur 
as $H$ would be an open proper subgroup  of the connected group $G.$
So, $\mathfrak h=\{0\}$ and this means that $H$ is discrete.$\bsq$

We give below  examples of co-amenable subgroups $H$
 in $SL_2(\RRR)$ which are not lattices.
These examples are based on the following proposition about ``transitivity" of co-amenability.
\begin{proposition}
\label{Pro-CoAmTrans}
Let $G$ be a locally compact group and $H\subset L$
be closed subgroups of $G.$
If $L$ is co-amenable in $G$ and $H$ is co-amenable in $L$,
then $H$ is co-amenable  in $G.$
\end{proposition}

\proof
Since $H$ is co-amenable in $L,$ the regular representation
$\pi_{L/H}$ of $L$ in $ L^2(L/H)$  weakly contains the trivial representation 
$1_L$. By continuity of induction (see Theorem F.3.5 in  \cite{BHV}),
the induced representation $\Ind_L^G( \pi_{L/H})$ weakly contains
$\Ind_L^G 1_L$.  Now, $\Ind_L^G( \pi_{L/H})$ is equivalent to 
the regular representation $\pi_{G/H}$  in $ L^2(G/H)$
and $\Ind_L^G 1_L$ to the regular representation $\pi_{G/L}$  in $ L^2(G/L)$.
It follows that $1_G$ is weakly contained
in $\pi_{G/H}$, since $1_G$ is weakly contained in  $\pi_{G/L}$
by the co-amenability of $L$  in $G.$ So, $H$ is co-amenable  in $G.\bsq$

\begin{example}
\label{Exa-CoAmHom}
Let $G$ be a locally compact group and $\Gamma$ a lattice in $G.$
Let $\vfi: \Ga \to S$ be a surjective homomorphism into an amenable discrete
group $S$. 
The previous proposition shows that $H= \ker(\vfi)$ is a co-amenable subgroup
in $G.$ 

For instance, if $G=PSL_2(\RRR)$ and $\Ga=F_2$ is the free group
on two generators realized as lattice in $G,$ then 
$H=[\Ga, \Ga]$ is co-amenable in $G$ and is not a lattice.
 This is also true for every co-compact lattice $\Gamma$
 in  $G=PSL_2(\RRR).$
 For other examples,  see \cite{Iozzi-Nevo} and \cite{Stuck}.

\end{example}
 
 As a consequence of  results from
\cite{Corlette} and \cite{Sullivan}, 
we now characterize the co-amenability of a (discrete) subgroup of  $SO(n,1)$  or $SU(n,1)$   in terms of the value of its 
critical exponent. 
Fix a $G$-invariant Riemannian metric $d$ on  the 
hyperbolic space $\HH^n(\RRR)$ in case $G=SO(n,1)$
or   $\HH^n(\CC)$  in case $G=SU(n,1),$   normalized to have constant
(respectively, maximal)  sectional curvature equal to $-1$ in the real case
(respectively, in the complex case).
  
Recall that the \emph{critical exponent} $\delta(\Ga)$ of   a discrete subgroup $\Gamma$ of  $G=SO(n,1)$  or $G=SU(n,1),$  
is defined as
$$\delta(\Ga):=   \inf \left\{s\in \RRR\ : \  \sum_{\ga\in \Ga} e^{-sd(\ga x, x)} <\infty \right\},$$
where $x$ is  an arbitrary  point  in $\HH^n(\RRR)$ or   $\HH^n(\CC).$ 
If $\Ga$ is a lattice in $G,$ then  $\delta(\Ga) =n-1$  in case $G=SO(n,1)$ 
and  $\delta(\Ga) =2n$  in case $G=SU(n,1);$  these are the maximal
possible values of $\delta(\Ga)$  for a discrete subgroup $\Ga.$

\begin{theorem}
\label{Theo-CoAm-CritExp}
Let $G=SO(n,1)$ for $n\geq 2$  or $G=SU(n,1)$ for $n\geq 1.$
The closed proper subgroups  of $G$ which are co-amenable 
are exactly the discrete subgroups $\Ga$ with maximal critical exponent, that is, such that
$\delta(\Ga)=n-1$ in case $G=SO(n,1)$ 
and  $\delta(\Ga) =2n$  in case $G=SU(n,1).$
\end{theorem}
\proof
Let $\Ga$ be a closed co-amenable proper  subgroup of $G;$ in view of Corollary~\ref{Cor-CoAM-Z},
we  can assume  that $\Ga$ is discrete, 

Let $\lambda_0\geq 0$ be the bottom of the
spectrum of the Laplace-Beltrami operator on the locally symmetric space 
$\Ga\backslash X,$ where $X=\HH^n(\RRR)$ or   $\HH^n(\CC).$ 
It is well-known and easy to prove that $1_G$ is not 
weakly contained in the regular representation $\pi_{\Ga\backslash G}$
on $L^2(\Ga\backslash G)$ if and only if $\lambda_0>0.$
Now, it is shown in \cite[\textsection 4]{Corlette} (see also  \cite[2.17]{Sullivan}) that
$\lambda_0=d^2/4$ if  $\delta(\Ga) \leq 1/2d$ and $\lambda_0=\delta(\Ga)(d-\delta(\Ga))$ if $\delta(\Ga) \geq 1/2d$,
where $d=n-1$ in case $G=SO(n,1)$ 
and  $d=2n$  in case $G=SU(n,1).$ So, $\lambda_0=0$
if and only if $\delta(\Ga)=d. \bsq$

\section{Quantifying the Spectral Gap Property}
\label{S:QuantitativeSG}
Let $G$ be a separable locally compact group
and  $G\curvearrowright X$ an action on a measure space $(X,m).$
Assume that the corresponding representation  
$\pi_X$ has the Spectral Gap Property.
(Recall that in case $m$ is  $G$-invariant probability measure, $\pi_X$ denotes the corresponding representation 
 in $L^2_0(X)$.)
  So, given a strongly adapted  probability measure $\mu$ on $G$,
we have $\Vert \pi_X(\mu)\Vert<1.$ 
 For various applications (see Section~\ref{S:Appl}), it is important to have
an upper bound for  $\Vert \pi_X(\mu)\Vert.$ 
Such a bound may involve  the norm of $\mu$   under some
``known" representations of the group $G,$  such as  the regular representation $\pi_G$.
As we now see, estimates of this kind are available when $G$ is a subgroup of a simple Lie group,
as a consequence of the strong decay of the matrix coefficients of their unitary representations.

Let $(\pi,\H)$ be a unitary representation of the 
locally compact group $G.$   
For a real number  $p$ with $1\leq p <+\infty,$ the representation $\pi$ is said to be 
\emph{strongly $L^p$},
if there is dense subspace $D\subset \H.$ 
such that, for every $\xi,\eta\in D,$  the matrix coefficient 
$$C^{\pi}_{\xi,\eta}: G\to \CCC, \quad g\mapsto \langle \pi(g)\xi,\eta\rangle$$
belongs to $L^p(G).$ Observe that then $\pi$ is 
strongly $L^q$ for any $q>p,$
since $C^{\pi}_{\xi,\eta}$ is bounded.

Strongly $L^p$-representations for $p=2$ or $p\geq 2+\eps$ are closely tied to the 
regular representation $\pi_G$.

\begin{proposition}
\label{Pro-TemperedRep}
Let $(\pi,\H)$ be a unitary representation of the 
locally compact group $G.$  
\begin{itemize}
\item[(i)] If  $\pi$ is strongly $L^2$, then  $\pi$ is contained in
 an infinite multiple of the regular representation $\pi_G.$
\item[(ii)] If $\pi$ is strongly $L^p$ for every $p>2,$ then  $\pi$ is weakly contained in
  the regular representation $\pi_G.$
\end{itemize}
\end{proposition}

Concerning the proofs, see Proposition~1.2.3 in Chapter V of \cite{HoTa} for  (i)
and  Theorem~1 in  \cite{CHH} for (ii).
Representations which  are strongly $L^p$ for every $p>2$ also 
often called tempered representations.

A crucial fact for the sequel is the following theorem;
part (i) is due to Cowling (Theorem 2.4.2 in \cite{Michael}) and
part (ii) to Moore (Proposition~3.6 in \cite{Moore}).

\begin{theorem}
\label{Theo-CowlingMoore}
Let $G$ be a simple real Lie group with finite center.
\begin{itemize}
\item[(i)] \textbf{(\cite{Michael})}
If  the real rank of $G$ is at least two, then there exists
$p(G)$ in $(2,+\infty)$ such that every unitary representation of $G$
with no no-zero invariant vectors  is strongly $L^p$ for every $p>p(G).$
\item[(ii)] \textbf{(\cite{Moore})} If  the real rank of $G$ is one, then  every unitary representation 
of $G$ with the Spectral Gap Property is strongly $L^p$ for some $p\in [2,+\infty)$
\end{itemize}
\end{theorem}

\begin{remark}
\label{Rem-CowlingMoore}
\n 
Part (i) of the previous theorem holds more generally when $G$ is the group of $\kk$-rational points
of a simple algebraic group over a local field $\kk$, with $\kk$-rank at least two
(see Theorem 5.6 in  \cite{HoMo}). 
Estimates for the optimal bounds $p(G)$ in (i) have been given 
in \cite{Howe}, \cite{Li1}, \cite{Li2} and \cite{Oh}.
For instance, one has  $p(SL_n(\RRR)) = 2n-2$ for $n\geq 3.$
Cowling's result also covers the case where  $G$ has $\RRR$-rank one and is a  Kazhdan group
(that is, when $G$ is locally isomorphic to  $Sp(1, n)$ or  $F_4^{-20})$).

To our knowledge, it is not known whether (ii) is true for
all local fields $\kk$  and all simple algebraic groups 
 with $\kk$-rank 1. We suspect that it is indeed the case.
 
  The results above are derived using pointwise 
 estimates for matrix coefficients of the involved unitary representation,
 and this provides a more precise and often useful information about them (see \cite{Michael}, \cite{HoTa}).
\end{remark}

Combining Theorem~\ref{Theo-CowlingMoore} with Proposition~\ref{Pro-TemperedRep},
we obtain the following remarkable result.
Recall that, for an integer $k\geq 1,$   the $k$-fold tensor product $\pi^{\otimes k}$
of a unitary representation $\Hpi$  is the unitary representation of $G$ 
defined  on the tensor product Hilbert space ${\H}^{\otimes k}$
by 
$$\pi^{\otimes k}(g) (\xi_1\otimes \cdots \otimes \xi_k)= \pi(g)\xi_1\otimes \cdots \otimes\pi(g) \xi_k \tout g\in G,\, \xi_1, \dots, \xi_k\in \H.
$$

\begin{corollary}
\label{Cor-Tensor}
Let $G$ be a simple real Lie group with finite center and
$\pi$ a unitary representation of $G$ with the Spectral Gap Property.
Then there exists $N \geq 1$ such that 
$\pi^{\otimes N}$ is contained in an infinite  multiple of the
regular representation $\pi_G$.
Moreover,
in case $G$ has Property (T), the integer $N$
can be chosen independently of $\pi.$

\end{corollary}
\proof
Let $1\leq p<\infty $ be such that  $\pi$ is strongly $L^p$.
Let $N$ be an integer with $N\geq p/2.$ 
Then $\pi^{\otimes N}$ is strongly $L^2$ and the claim
follows. $\bsq$

Given a probability measure, we now deduce from the previous corollary estimates
for the norm $\pi(\mu)$ in terms of $\pi_G(\mu)$,
For this, we will use in a crucial way the following estimate which
 appears in the  proof of Theorem~1 in \cite{Nevo}.
 Recall that, for a unitary representation $(\pi,\H)$ of $G$, the contragredient (or conjugate) representation
 $\overline\pi$ acts on the conjugate Hilbert space $\overline \H$.
 
\begin{proposition}
\label{Pro-Nevo}
\textbf{(\cite{Nevo})}
 Let $\mu$  be a  probability measure  on 
the locally compact group $G.$
 Let  $(\pi,\H)$ be a unitary representation  of $G$.
For every integer  $k\geq 1,$ we have 
$$
\Vert \pi(\mu)\Vert \leq \Vert \left( \pi\otimes\overline{\pi}\right)^{\otimes k}(\mu)\Vert^{1/2k},
$$
\end{proposition}
\proof
 Using Jensen's inequality, we have  for every   vector $\xi\in\H,$
\begin{align*}
\Vert\pi(\mu)\xi\Vert^{4k} 
&=\vert\langle\pi(\check{\mu}\ast\mu)\xi,\xi \rangle\vert^{2k}=\left\vert\int_{ G} \langle \pi(g)\xi,\xi\rangle
d(\check{\mu}\ast\mu)(g)\right\vert^{2k}\\
&\leq\int_G\langle \vert\pi(g)\xi,\xi \rangle\vert^{2k}
d(\check{\mu}\ast\mu)(g)\\
&=\int_G\vert\langle (\pi\otimes\overline{\pi})(g)(\xi\otimes \xi),\xi\otimes \xi\rangle\vert^k d(\check{\mu}\ast\mu)(g)\\
&=\int_G\langle (\pi\otimes\overline{\pi})^{\otimes k}(g)(\xi\otimes \xi)^{\otimes k},(\xi\otimes \xi)^{\otimes k} \rangle d(\check{\mu}\ast\mu)(g)\\
&=\vert\langle(\pi\otimes\overline{\pi})^{\otimes k}(\check{\mu}\ast\mu)(\xi\otimes \xi)^{\otimes k}, (\xi\otimes \xi)^{\otimes k}\rangle\vert\\
&= \Vert (\pi\otimes\overline{\pi})^{\otimes k}(\mu)(\xi\otimes \xi)^{\otimes k}\Vert^{2}.
\end{align*}
and the claim follows. $\bsq$

  \begin{corollary}
\label{Cor-Bound}
Let $G$ be a simple real Lie group with finite center and
$\pi$ a unitary representation of $G$ with the Spectral Gap Property.
Then there exists $N \geq 1$ such that, for every 
 probability measure $\mu$ on $G$, we have
$$
\Vert \pi(\mu)\Vert \leq \Vert \pi_G(\mu)\Vert^{1/N}.
$$
Moreover, in case $G$ has Property (T), the integer $N$
can be chosen to be independent of $\pi.$ 
\end{corollary}

\proof
 Let $1\leq p<\infty $ be such that  $\pi$ is strongly $L^p$.
Let $k$ be an integer with $k\geq p/4.$ Then 
$\left(\pi\otimes\overline{\pi}\right)^{\otimes k}$ is strongly
 $L^2$ and is hence is contained in  
 a  multiple of $\pi_G.$ 
 The claim follows now from Proposition~\ref{Pro-Nevo}. $\bsq$

\begin{remark}
\label{Rem-Gorodnik-Nevo}
Results about decay of matrix coefficients of semisimple algebraic  groups $G$
as  described in Theorem~\ref{Theo-CowlingMoore} and Corollary~\ref{Cor-Bound}
have been used
in  the monograph  \cite{GoNe} in order to prove impressive
quantitative ergodic theorems for  families of averaging operators on $G$
or on a lattice  in $G$ as well as  results on counting lattice points (with explicit error term),
with various applications to counting problems  from number theory.
\end{remark}

\section{Probability measure preserving actions }
From now on, we will deal only with group actions which
preserve a probability measure. 
So, let $G$ be a locally compact group 
acting on a  measure space $(X,m),$ 
where $m$ is a $G$-invariant probability measure. 
Recall that $G\curvearrowright X$ has  the Spectral Gap Property
if the Koopman representation $\pi_X$
on $L^2_0(X)= (\CCC{\Un}_X)^\perp$
does not have almost invariant vectors.

\subsection{Amenability as obstruction to the Spectral Gap Property}
 The following result shows that amenability is  
 an obstruction to the  Spectral Gap Property also in the probability measure
 preserving case, at least for discrete groups.
 \begin{theorem}
 \label{Theo-Rosenblatt}
 \textbf{(\cite{JuRo})}
 Let $\Ga$ be a countable amenable group
 with a measure preserving action on a \emph{non atomic}
 probability space  $(X,m)$.
 Then  $\Ga\curvearrowright X$ does not have the Spectral Gap Property.
 \end{theorem}
 
 \proof
 Let $S$ be a finite symmetric subset of $\Ga$
 and $\eps>0.$ We want to construct
 a function $f\in L^2_0(X)$ which is $(S, \eps)$-invariant.
 Since $\Ga$ is amenable, 
there exists   a finite subset $F$ of $\Ga$
with  $|F s \triangle F| \leq \eps |F|$ for all $s\in S,$ 
by   F\o lner's Theorem~\ref{Theo-Folner}.

Since $m$ is not atomic, there exists  a measurable subset $A$
of $X$ with $m(A)= {1}/{2|F|}.$
Consider the function $\xi: X\to \NN$ 
defined by $\xi =\sum_{\ga\in F} \pi_X(\ga^{-1})({\bf 1}_A):$
$$ \xi(x) =  \sum_{\ga\in F} {\bf 1}_{\ga^{-1} A}(x) \tout x\in X.
$$
 We have 
 \begin{align*}
 \Vert \xi\Vert_1&= \sum_{\ga\in F} \int_X \pi_X(\ga^{-1})({\bf 1}_A) dm(x) =  \sum_{\ga\in F} m(A) = |F| m(A) = \dfrac{1}{2}.
 \end{align*}
 Let $s\in S$ be fixed. Then 
\begin{align*} 
  \Vert \pi_X(s^{-1})\xi- \xi\Vert_1&= \left\Vert \sum_{\ga\in F}  {\bf 1}_{(\ga s)^{-1}A} - {\bf 1}_{\ga^{-1}A} \right\Vert_1 \\
 &= \left \Vert\sum_{\ga\in Fs}  {\bf 1}_{\ga^{-1}A}- \sum_{\ga\in F}{\bf 1}_{\ga^{-1}A}\right \Vert_1\\
 &=  \left\Vert\sum_{\ga\in Fs\setminus F}   {\bf 1}_{\ga^{-1}A}-  \sum_{\ga\in F\setminus Fs}{\bf 1}_{\ga^{-1}A} \right\Vert_1\\
 &\leq  |F s \triangle F| m(A) \leq \eps |F| m(A) = \eps/2.
 \end{align*}

  Let  $g= \sqrt{\xi}.$ Then $\Vert g\Vert_2^2= \Vert \xi\Vert_1= 1/2$ 
 and 
 $$
 \Vert \pi_X(s^{-1})g- g\Vert_2^2 \leq \Vert \pi_X(s^{-1})\xi- \xi\Vert_1 \leq  \eps/2.
 $$
 for all $s\in S.$
  Since ${\xi}$ takes its values in $\NN,$ we have $g= \sqrt{\xi} \leq \xi.$
 Hence,
 \begin{align*}
( \int_X g dm)^2 &\leq ( \int_X \xi dm)^2 = \Vert \xi\Vert_1^2= 1/4
 \end{align*}
 and so 
 $\Vert g - \int_X g dm {\bf 1}_X\Vert^2_2 \geq \dfrac{1}{2}- \dfrac{1}{4}= \dfrac{1}{4}.$
 Set now 
 $$f:= \dfrac{g-  \int_X g dm {\bf 1}_X}{\Vert g-  \int_X g dm {\bf 1}_X\Vert_2}.$$
 Then  $\Vert f\Vert_2=1$ and
 $$
 \Vert \pi_X(s^{-1})f- f\Vert_2^2\leq \dfrac{4\eps}{2}= 2\eps. \bsq
 $$
\begin{remark}
\label{Rem-JuRo}
Theorem~\ref{Theo-Rosenblatt} may fail when $G$ is amenable
and not discrete.
Indeed, when $G$ is a compact infinite group,
the action $G$ on itself by left translation has the Spectral Gap
Property (as probability measure preserving action).
\end{remark}

 Let $\Ga\curvearrowright X$ be  a measure preserving action on a 
 probability space $(X,m)$. A \emph{factor} of the triple
 $(X, m, \Ga)$ is a triple $(Y, m', \Ga')$ 
 where $\Ga'$  is a group acting 
 on a  probability space $(Y,m')$
 such that there exists a homomorphism $\Ga\to \Ga'$  and a
  measurable map $\Phi : X\to Y$ with  $\Phi_*(m)= m'$
 and  intertwining the two actions. 
  Here is a corollary of Theorem~\ref{Theo-Rosenblatt}.
 \begin{corollary}
 \label{Cor-JuRo}
 Let  $(Y, m', \Ga')$  be a factor of  $(X, m, \Ga)$.
 Assume that $m'$ is non atomic and that
 $\Ga'$ is amenable. 
 Then  $\Ga\curvearrowright X$ does not have the Spectral Gap Property.
 \end{corollary}
\proof
 Observe that $\{f \circ \Phi: f\in L^2_0(Y,m')\}$
is  a $\Ga$-invariant closed subspace of 
$L^2_0(X,m).\bsq$

\begin{remark}
The previous corollary shows that, given a discrete group $\Ga$
and an action $\Ga\curvearrowright X$ on a 
 probability space $(X,m),$ the existence
 of a factor $(Y, m', \Ga')$   of  $(X, m, \Ga)$,
 with $\Ga'$ amenable and $m'$ non atomic,
 is an obstruction to the Spectral Gap Property
 for $\Ga\curvearrowright X$. In some cases, such 
 as $\Ga \subset SL_n(\ZZ)$ acting by automorphisms
 on $X= \TT^n,$ this is the only obstruction
 (see Theorem~\ref{Theo-SGP-AutoTorus}).
 \end{remark}

\subsection{Spectral Gap Property  and Orbit Equivalence}
\label{SS:OE}
 Orbit Equivalence is a subject where one
 aims to study  groups $G$ and group actions $G\curvearrowright X$
 on measure spaces through the equivalence relation 
 on $X$ defined by the partition   into $G$-orbits
  (for a recent survey on this theme, see \cite{Damien}).

 \begin{definition}
\label{Def-OE}
Two  measure preserving ergodic actions of two countable groups $\Ga_1$ and $\Ga_2$ 
on  probability spaces $(X_1,m_1)$ and $(X_2, m_2)$ are
orbit equivalent 
 if there exist  measurable subsets 
$X_1'$ and $X_2'$ with measure  $1$ in  $X_1$
and $X_2$ and a Borel isomorphism 
$f:X_1'\to X_2'$
 with $f_*(m_1)=m_2$ such that,
for $m_1$-almost every $x\in X_1',$
we have  $$f(\Gamma_1 x)=\Gamma_2 f(x).$$
\end{definition}

A natural question is whether the Spectral Gap Property
is an invariant  of orbit equivalence.
This is the not case,
as is  shown 
by the following result of Hjorth and Kechris
(Theorem A. 3.2 in \cite{Hjorth-Kechris}).
Recall that, if $\Ga$ is a non amenable group,
its Bernoulli action  $\Ga\curvearrowright\{0,1\}^\Ga$
 has the Spectral Gap Property
(see Example~\ref{Ex-Bernoulli}).

\begin{theorem}
\label{Theo-Hjorth-Kechris}
\textbf{(\cite{Hjorth-Kechris})}
The Bernoulli action of the free group $\Ga=F_2$ on 2 generators
on $X=\{0,1\}^\Ga$  is orbit equivalent to another action of
$\Ga$ on $X$ which does not have the Spectral Gap Property.
\end{theorem}

 There is however a  property related to the Spectral Gap Property
  which is an  invariant  of orbit equivalence, namely
 strong ergodicity in the sense
 of  Schmidt (see \cite{Schmidt1}, \cite{Schmidt2}),
 
 Let $\Ga\curvearrowright X$ be a measure preserving action 
 of a  countable group $\Ga$ on
the probability space  $(X,m).$ If this action is ergodic,
there is no non-trivial invariant measurable subset of 
$X.$ Nevertheless, there might exist non-trivial
asymptotically invariant  subsets in the following sense.
 \begin{definition}
 \label{Def-StrongErgodic}
 A sequence of measurable subsets $(A_n)_n$ of $X$
is said to be {asymptotically invariant} if 
$$
\lim_n  m(\ga A_n\triangle A_n)=0  \tout  \ga\in \Ga.
$$
An asymptotically invariant sequence $(A_n)_n$
is said to be {non-trivial} if
$$
\inf_n m(A_n)(1-m(A_n))>0.
$$
The action of $\Ga$ on $X$ is 
\emph{strongly ergodic} if there exists
no non-trivial  asymptotically invariant sequence in $X.$
\end{definition}

Strong ergodicity
is an invariant of orbit equivalence.
Indeed, let 
$[\Ga]$ denote the \emph{full group} of $\Ga,$ that is, the group 
of all measure preserving measurable bijections
 $\varphi: X\to X$ with $\varphi(x) \in \Ga x$ for 
 $m$-almost every $x\in X.$ The following
 lemma shows that the collection of asymptotically invariant sequences 
 for $\Ga$ and hence strong ergodicity only depends
 on the equivalence relation on $X$ defined by the action of
 $\Ga.$

 \begin{lemma}
 \label{Lem-AIS}
 Let $(A_n)_n$ be an asymptotically invariant sequence 
 for $\Ga.$ Then 
 $\lim_n  m(\varphi (A_n)\triangle A_n)=0$,
 for every $\varphi\in [\Ga].$
 \end{lemma}
 \proof
 Let $\eps>0.$  There exist  Borel subsets $X_1, \dots, X_k$ of $X$
 with 
 $$m(X\setminus \bigcup_{i=1}^k X_i)\leq \eps/4$$
  and elements $\ga_1, \dots, \ga_k$ in $\Ga$ such that 
 $\vfi(x)= \ga_i x$ for all $x\in X_i$ and $i=1,\dots,k.$
 Let $N\in \NN$ be such that 
 $$m(\ga_i A_n\triangle A_n)\leq \eps/4k \tout n\geq N,\,  i= 1,\dots,k.$$
 Since
 $$
 \varphi(A_n)\setminus A_n \subset   \bigcup_{i=1}^k (\ga_i A_n\setminus A_n) \cup \vfi(X\setminus  \bigcup_{i=1}^k X_i),
 $$
 it follows that
 $$m(\varphi(A_n)\setminus A_n)\leq \eps/2 \tout n\geq N,\,  i= 1,\dots,k.$$
 As
 $m(A_n\setminus \varphi(A_n))= m(\varphi^{-1}(A_n)\setminus A_n)$
 and, similarly,
$$
 \varphi^{-1}(A_n)\setminus A_n \subset   \bigcup_{i=1}^k(\ga_i^{-1} A_n\setminus A_n )\cup  \vfi^{-1}(X\setminus  \bigcup_{i=1}^k \ga_iX_i),
 $$
 we obtain $ m(\varphi (A_n)\triangle A_n)\leq \eps$ for all $n\geq N.\bsq$

The following elementary proposition shows
that the  Spectral Gap Property implies strong ergodicity
and this is often the way strong ergodicity is established
in specific examples.

\begin{proposition} 
\label{Pro-AsympInv}
Let $\Ga\curvearrowright X$ be a measure preserving action 
 of a  countable group $\Ga$ on
the probability space  $(X,m).$ 
If  $\Ga\curvearrowright X$ has the  Spectral Gap Property,
then $\Ga\curvearrowright X$ is strongly ergodic.
\end{proposition}

\proof
Assume, by contradiction, that there exists a 
 non-trivial asymptotically invariant sequence $(A_n)_n$
of measurable subsets of $X.$   Set 
$$f_n={\mathbf 1}_{A_n}-m(A_n){\mathbf 1}_{X}.$$
 Then $f_n\in L^2_0(X,m)$ and, for every $\ga\in \Ga,$
 we have
$$
\Vert f_n\Vert^2=m(A_n)(1-m(A_n))\qquad\text{and}\qquad
\Vert\pi_{X}(\ga)f_n-f_n\Vert^2=m(\ga A_n\triangle A_n).
$$
So, $\pi_{X}^0$  has almost invariant vectors and this is a contradiction.
$\bsq$

The following example, due to K. Schmidt, shows that  the converse does not hold in the previous 
proposition (that is, strong ergodicity and the Spectral Gap Property do not coincide).
Of course, since strong ergodicity is an invariant of  orbit equivalence,
this also follows from Theorem~\ref{Theo-Hjorth-Kechris} above.

\begin{example}
\label{Exa-Schmidt}
\textbf{(\cite{Schmidt2})}
Let $\Ga=F_3$ be the free group on the three generators 
$a,b,c$. Let $\Ga_2$ be the free subgroup of  $\Ga$ generated by $a,b$ and 
let $Y={\mathbf T}^2$ be the $2$-torus endowed with Lebesgue measure $\la.$
Then $\Ga_2$  acts as a group of   automorphisms on $Y$
through a surjective homomorphism $\Ga_2\to SL_2(\ZZ)$.
The action $\Ga_2\curvearrowright Y$ has the Spectral Gap Property
(see Example~\ref{Exa-Torus}) and is hence strongly ergodic by the
previous proposition.

Set $X= Y\times \NN$ and define a finite positive measure $m$ on $X$
by 
$$
m(A)= \sum_{n=1}^{\infty} \frac{1}{n^2} \la\left(\{ y\in Y \ : \ (y,n) \in A\}\right)
$$
for every Borel subset $A$ of $X.$ 
We define a measure preserving action of $\Ga_2$
on $(X,m)$ by
$$
(\ga, (y,n))\mapsto (\ga y, n) \tout \ga\in \Ga_2, \, (y,n)\in X.
$$

Let $T:X\to X$ be a measure preserving  bijection 
with $T(X_n) \subset X_{n-1}$ for all $n\geq 2,$
where $X_n= Y\times \{n\}.$
Such a mapping $T$ can be constructed as follows.
For every $n\geq 2,$ choose a Borel subset $Z_{n-1}$ of $X_{n-1}$ 
with $m(Z_{n-1})= m(X_{n})= 1/n^2.$
There exist a  measure preserving  bijection $T_{n}: X_{n}\to Z_{n-1}$ for every $n\geq 2$
as well as a  measure preserving  bijection $T_1: X_1\to X\setminus \cup_{n\geq 2} Z_{n-1}.$
This gives rise to   a  measure preserving  bijection $T: X\to X$
 defined by $T|_{X_n}= T_n$ for all $n\geq 1.$

We extend the  action of $\Ga_2$ on $X$ to a measure preserving action of $\Ga$ by 
letting $c$ act as $T.$
 As we are going to show, the action $\Ga\curvearrowright X$ is strongly ergodic and does not have the Spectral Gap Property.

Let $(A_k)_k$ be an asymptotically invariant sequence 
 for $\Ga$. For every $n\geq 1,$ let $A_{k,n}= A_k\cap X_n.$
 Since $\Ga_2$ leaves $X_n$ invariant,
  $(A_{k,n})_k$ is an asymptotically invariant sequence 
 for $\Ga_2$. Hence, 
 $$\lim_k m(A_{k,n})(m(X_n)-m(A_{k,n}))=0$$ 
and therefore, for every $n\geq 1,$ we have either $\lim_k m(A_{k,n})=0$ or  $\lim_k m(A_{k,n})=m(X_n).$
 We claim that we have either $\lim_k m(A_{k,n})=0$ for all $n$ or 
 $\lim_k m(A_{k,n})=m(X_n)$ for all $n.$ Once proved, this will
 imply that  $\lim_k m(A_{k})=0$ or $\lim_k m(A_{k})=m(X)$
 
  To prove the claim, assume
   that 
  $\lim_k m(A_{k,n})=m(X_n)$ for some $n.$ 
  Then
  $\lim_k m(A_{k,n-1})=m(X_{n-1})$ and  $\lim_k m(A_{k,n+1})=m(X_{n+1})$
  in case $n\geq 2$ and $\lim_k m(A_{k,l})=m(X_{l})$ for all $l\geq 1,$ 
  in case $n=1.$
  
  Indeed, let  $n\geq 2.$
  Since $T(A_{k,n})\subset X_{n-1}$ and $T^{-1}(A_{k,n})\subset X_{n+1}$
  and since $\lim_k m(T^{\pm 1} A_k\triangle A_k)=0,$
  it follows that   $\lim_k m(A_{k,n-1})\neq 0$ and $\lim_k m(A_{k,n+1})\neq 0$
  and so $\lim_k m(A_{k,n-1})=m(X_{n-1})$ and  $\lim_k m(A_{k,n+1})=m(X_{n+1})$.
  The case $n=1$ is treated similarly.
     This proves that  the action  $\Ga\curvearrowright X$ is strongly ergodic.
   
  For every $n\geq 2,$ let $C_n= \cup_{k\geq n} X_k.$ 
  Then $C_n$ is $\Ga_2$-invariant, 
  $$m(T(C_n) \triangle C_n)=\dfrac{1}{(n-1)^2},$$
  and 
  $$m(C_n)= \sum_{k\geq n} \dfrac{1}{k^2}={\rm O}\left(\dfrac{1}{n}\right).$$
  Hence, $a(C_n) \triangle C_n =\emptyset, b(C_n) \triangle C_n =\emptyset,$ and
  $$\lim_n \dfrac{m(T(C_n) \triangle C_n)}{m(C_n)}= 0.$$
  This implies that 
  $$\lim_n \dfrac{m(\ga C_n \triangle C_n)}{m(C_n)}= 0 \tout \ga\in \Ga.$$
 With
$f_n:={\mathbf 1}_{C_n}-m(C_n){\mathbf 1}_{X} \in  L^2_0(X,m),$
we have,  for every $\ga\in \Ga,$
$$
\lim_n\dfrac{\Vert\pi_{X}(\ga)f_n-f_n\Vert^2}{\Vert f_n\Vert^2}=\lim_n\dfrac{m(\ga C_n\triangle C_n)}{m(C_n)(1-m(C_n))}=0,
$$
and this shows that  $\Ga\curvearrowright X$ does not have the Spectral Gap Property.
\end{example}

 Schmidt proved in \cite[(2.4) Theorem]{Schmidt2}  the following  stronger version of Theorem~\ref{Theo-Rosenblatt}.
 \begin{theorem}
\label{Theo-Schmidt}
\textbf{(\cite{Schmidt2})}
Let $\Ga$ be a countable amenable group.
Then no measure preserving action $\Ga\curvearrowright X$ 
on  a non-atomic  probability space $(X,m)$  is strongly ergodic.
\end{theorem}

\section{Spectral gap property for homogeneous  spaces}
\label{S:Homogneous}
In this section, we will deal with the class
of probability measure preserving actions arising from
homogeneous spaces associated to lattices.
 The setting we consider is as follows.
 
 Let $G$ be a locally compact group (for instance, a connected Lie group)
 and $\Ga$ a lattice in $G.$ Let $G/\Ga$ be equipped
 with its unique $G$-invariant probability measure $m.$
 We consider the action of $G\curvearrowright G/\Ga$
 and the corresponding Koopman representation $\pi_{G/\Ga}$
 of $G$ on $L^2_0(G/\Ga).$ 
 We then ask: does  the action $G\curvearrowright G/\Ga$ have the Spectral Gap Property?

\subsection{The case of a cocompact lattice}
\label{SS:Uniform}
% For notational reasons, we prefer here to consider the space $\Ga\backslash G$
% instead of $G/\Ga$
As we now see, the question above has a positive answer for \emph{uniform}  lattices.
 \begin{proposition}
 \label{Pro-Uniform}
 Let $G$ be a locally compact group and  $\Ga$ a
 cocompact lattice in $G.$
 Then $G\curvearrowright G/\Ga$ has the Spectral Gap Property. 
 \end{proposition}
\proof
We denote by $\pi$ the Koopman representation on $L^2(G/\Ga)$
(and \emph{not} on $L^2_0(G/\Ga)$).

We first check the crucial fact that, for every $f\in C_c(G),$ 
 the convolution operator $\pi(f): L^2(G/\Ga)\to L^2(G/\Ga)$
 is a compact operator.
  Indeed,  let $X\subset G$ be a compact  fundamental domain for the action of $\Ga$ on $G.$
  Thus, $X$ is a compact subset of $G$ such that 
  $$G=\coprod_{\ga\in \Ga} X\ga.$$
  We denote by  $m$  a Haar measure on $G$.
  (Observe that $G$ is unimodular, since it contains a lattice and so $m$ is right and left invariant;
  see Proposition B.2.2 in \cite{BHV}.)
  
 View $\xi\in L^2(G/\Ga)$ as a function on $G$ which is $\Ga$-invariant on the right; then,
  for every $x\in X\cong G/\Ga,$
  \begin{align*}
  \pi(f) \xi(x)&= \int_G f(g) \xi(g^{-1} x)dm(g) \\
  &= \int_G f(xg^{-1}) \xi(g)dm(g)\\
  &=\int_X \sum_{\ga\in\Ga} f(x\ga^{-1}y^{-1}) \xi(y\ga)dm(y)\\
   &=\int_X K(x,y) \xi(y)dm(y),
  \end{align*}
  where $K(x,y)= \sum_{\ga\in\Ga} f(x\ga^{-1}y^{-1}).$ 
  Observe that  there are only finitely many 
  $\ga$ in  $\Ga$ for which $x\ga^{-1}y^{-1}$ is in the compact support of $f$
  for some $(x,y)\in X\times X$. So, $K$ is  continuous on $X\times X$ and hence
  $\pi(f)$ is an integral operator with continuous kernel. Since, $X$ is compact,
  $\pi(f)$ is therefore a Hilbert-Schmidt operator on $L^2(X)=L^2(G/\Ga).$

 Now,  let $f\in C_c(G)$  with  $f\geq 0, \int_G f(g)dm(g) =1$ 
 and $\check f= f$ and  such that 
 $\supp (f)$ generates $G$. Then $\pi(f)$
 is a compact self-adjoint operator. Hence,
 $1$ is an isolated spectral value of $\pi(f)$ and so 
 $G\curvearrowright G/\Ga$ has the Spectral Gap Property,
 by Proposition~\ref{SGP-Norm}.$\bsq$
\begin{remark}
\label{Rem-Uniform}
In fact, pushing the analysis a bit further, one can show that 
$L^2(G/\Ga)$ decomposes as a Hilbert space direct sum 
$L^2(G/\Ga)= \bigoplus_{i} m_i\H_i$ of 
irreducible $G$-invariant subspaces $\H_i$, each of which occuring with finite multiplicity $m_i$
(see Theorem in Chap. I, Section 2.3 in \cite{GGPS}).
\end{remark}

  \subsection{The case of a non cocompact lattice}
 \label{SS:Nonuniform}.
 The problem of establishing the Spectral Gap Property 
 for  $G\curvearrowright   G/\Ga$ is much harder in the case 
 of a non cocompact lattice $\Ga.$ Only partial results are known.
 
Part (i) of the following theorem  has been conjectured in \cite[Chapter III.  Remark 1.12]{Margulis}
and  proved in  \cite{Bekka-Cornulier}; part (ii) is from  \cite{BeLu}.
\begin{theorem}
\label{Theo-LieGroups} 
\textbf{(\cite{Bekka-Cornulier}, \cite{BeLu})}
Let $G$ be a locally compact group and $\Ga$ a lattice in $G$.
Then  $G\curvearrowright   G/\Ga$ has the Spectral Gap Property in the following cases:
\begin{itemize}
 \item [(i)] $G$  is a real Lie group;
 \item[(i)] $G=\mathbf G(\mathbf k)$ is the group of $\mathbf k$-rational points
  of  a simple algebraic  group $\mathbf G$ over  a local field $\mathbf k.$
\end{itemize}
 \end{theorem}

Concerning part (ii) of the theorem, observe that when $\mathbf k$ is non-archimedean 
with characteristic $0$, every lattice $\Ga$ in $\mathbf G(\mathbf k)$
is cocompact (see \cite[p.84]{Serre}) and  the result
follows from Proposition~\ref{Pro-Uniform}. By way of contrast, $G$ has many 
non uniform lattices when the characteristic of $\mathbf k$ is non zero
(see \cite{Serre} and \cite{Alex-rank1}).
So, for the proof, it suffices to consider the case where the characteristic of $\mathbf k$ is non-zero
and where $\mathbf k-{\rm rank} (\mathbf G)=1.$ (Recall that $\mathbf G(\mathbf k)$
has Kazhdan's Property (T) when $\mathbf k-{\rm rank} (\mathbf G)\geq 2$; see Theorem~\ref{Theo-T-Rank2}.)
In this situation, it is known that $G=\mathbf G(\mathbf k)$ acts by automorphisms on the associated Bruhat-Tits tree $T$ (see \cite{Serre}),
which is a regular or a bi-partite bi-regular tree. 
The proof then is reduced to showing that the projection on the quotient graph $X=\Ga\backslash T$ of the standard random walk on $T$
  has a spectral gap.

\subsection{An example: the case of $PGL_2(\FF_q((t^{-1})))/PGL_2(\FF_q[t])$}
\label{S:PGL}
We will give in this section a complete proof of part (ii) of  Theorem~\ref{Theo-LieGroups}   for the special case $G=PGL_2(\FF_q((t^{-1})))$
 and $\Ga=PGL_2(\FF_q[t]),$ where $\FF_q((t^{-1}))$
 is the local field of formal Laurent series with coefficients in the finite field with $q$ elements.
Following an idea from \cite{BeLu}, the proof uses  a version of Cheeger's inequality
for Markov chains on a countable state space for which we give a full proof.

\subsubsection{A Cheeger inequality for Markov chains on a countable state space}
\label{SS-Cheeger}
Let $X$ be a countable set and let $\mu$ be a Markov kernel on $X,$
that is, a mapping $\mu: X\times X\to \RRR^+$
such that $\sum_{y \in X} \mu(x,y) = 1$ for all $x \in X$.
Such a kernel defines a Markov chain $(Z_n)_{n\geq 0}$ on $X$, with transition probabilities
$$\PP(Z_{n+1}=y | Z_{n}= x)= \mu(x,y).$$

We assume that   $\mu$ is irreducible, that is,
given any pair $(x,y)$ of distinct points in $X$,
there exist an integer $n \ge 1$ and a sequence 
$x = x_0, x_1, \hdots, x_n = y$ in $X$
such that $\mu(x_{j-1},x_j) > 0$ for any $j \in \{1,\hdots,n\}$.
We also assume that $\mu$ is reversible:  there exists a stationary measure $m$ for  $\mu$,
that is a function $m: X \to \RRR_+^*$
such that 
$$
m(x)\mu(x,y)  = m(y)\mu(y,x)
\quad \text{for all} \quad
x,y \in X .
$$
The corresponding  Markov operator $M_\mu$  on $\ell^2(X,m)$  is defined by
$$
M_\mu f(x)= \sum_{(x,y)\in X^2 } \mu(x,y)f(y) \tout f\in \ell^2(X,m).
$$
Since $m$ is stationary  measure for $\mu,$ one checks that the operator $M_\mu$ is self-adjoint
with $\Vert M_\mu\Vert \leq 1.$

From now on, we assume that the stationary measure is a \emph{finite measure}
and, hence without loss of generality, that $m$  is a \emph{probability measure.}
Then $1$ is an eigenvalue of  $M_{\mu},$ with ${\mathbf 1}_X$ as eigenfunction.
We will be concerned with finding upper bounds for the spectrum $\sigma (M_{\mu})$ 
of $M_\mu$ restricted  to  $\ell^2_0(X,m)= (\CCC {\mathbf 1}_X)^{\perp}$.

It is convenient to consider the Laplacian  $\Delta_\mu= I-M_{\mu},$
which is a nonnegative operator on $\ell^2(X,m)$ with $\Vert \Delta_\mu\Vert \leq 2.$
So, we seek a lower bound for 
$$\lambda_1= \inf \sigma \left(\Delta_{\mu}|_{\ell^2_0(X,m))}\right).$$
Such a bound
will be given in terms of the \emph{Cheeger constant} $h(X)$ of the random walk $\mu$
which is defined as follows.

Let $\widetilde{\mu}$ be the (symmetric) measure on $X\times X$ defined by 
$$
\widetilde{\mu}(x,y) =m(x) \mu(x,y)\quad \text{for all} \quad (x,y)\in X\times X.
$$
Set 
$$h(X):= \inf \frac{\widetilde{\mu}(S \times S^c)}{m(S)m(S^c)},$$
where the infimum is taken over all non empty subsets $S$ of $X$ and where $S^c=X\setminus S.$

The Cheeger inequality is an isoperimetric inequality originally proved for the
Laplacian acting on the $L^2$-space of a compact  Riemannian manifold (see \cite{Cheeger}, \cite{Chavel})
and carried over to the setting of Markov chains  in \cite{LaSo},\cite{Sinclair}.
Versions of Cheeger's inequality for  weighted graphs were considered by several authors (see, for instance, \cite{Alon}, \cite{Diaconis},
\cite{Dodziuk}, \cite{Mokhtari}, \cite{Morgenstern}, \cite{Bauer}).

\begin{theorem}
\label{Theo-Cheeger}
\textbf{(\cite{LaSo},\cite{Sinclair})}
We have 
$$\frac{h(X)^2}{8} \leq \lambda_1.$$
Consequently, if $h(X)>0$ then $\lambda_1>0.$
\end{theorem}

We now proceed with the proof of Cheeger's inequality.
The first ingredient is the following lemma, which is straightforward to check.
\begin{lemma}
\label{Lem1-Cheeger}
\n
(i) $\lambda_1$ is the infimum of the Rayleigh quotients  $ \frac{\langle \Delta_{\mu} f  ,  f \rangle}{\Vert f\Vert^2} $
over all $f\in \ell^2_0(X), f\neq 0$ and $f$ real valued.

\n
(ii) For every $f\in \ell^2(X,m),$ we have 
$$ \langle \Delta_{\mu} f  ,  f \rangle= \frac{1}{2}\sum_{(x,y) \in X^2} |f(y) - f(x)|^2 m(x) \mu(x,y)= 
\frac{1}{2}\sum_{(x,y) \in X^2} |f(y) - f(x)|^2 \widetilde{\mu}(x,y). \bsq$$ 
\end{lemma}

The next ingredient in the proof are the so-called  area and co-area formulas.

\begin{lemma}
\label{Lem2-Cheeger}
Let $u: X\to \RRR^+$  be in $\ell^1(X,m)$ and set $S_t= \{ x\in X\ : \ u(x) >t\}$ for $t\geq 0.$
Then the following formulas holds:

\n
(i) \textbf{(Area formula)} 
$$\sum_{x\in X} u(x)m(x) =\int_{0}^{\infty}m(S_t) dt.$$

\n
(ii) \textbf{(Co-area formula)} 
$$
\dfrac{1}{2}\sum_{(x,y) \in X^2} |u(y) - u(x)| \widetilde{\mu}(x,y) =\int_{0}^{\infty} \widetilde{\mu}(S_t \times S_t^c) dt.
$$
\end{lemma}

\proof
(i) For $x\in X,$ we have $x\in S_t$ if and  only if $\mathbf{1}_{(t,\infty)} (u(x))=1$ and 
hence 
\begin{align*}
\int_{0}^{\infty}m(S_t) dt &= \int_{0}^{\infty} \left(\sum_{x\in X} m(x)\mathbf{1}_{(t,\infty)} (u(x)) \right) dt= 
\sum_{x\in X} m(x) \int_{0}^{\infty} \mathbf{1}_{(t,\infty)} (u(x)) dt\\
&=\sum_{x\in X} m(x) \int_{0}^{\infty} \mathbf{1}_{[0,u(x))} (t) dt= \sum_{x\in X} m(x) u(x).
\end{align*}

\n (ii)  
For $(x,y)\in X\times X,$ 
denote by $I_{x,y} \subset \RRR^+$ the interval between $u(x)$ and $u(y).$ So, $|I_{x,y}|=|u(y) - u(x)|.$ 
We have  $(x,y)\in (S_t\times S_t^c)\cup (S_t^c\times S_t)$ if and  only if $t\in I_{x,y},$
that is if and only if $\mathbf{1}_{I_{x,y}} (t) =1.$
Since $\widetilde{\mu}$ is symmetric, we have 
$$\widetilde{\mu}\left((S_t\times S_t^c)\cup (S_t^c\times S_t)\right)= 2\widetilde{\mu}(S_t \times S_t^c)$$
 and hence 
\begin{align*}
2\int_{0}^{\infty} \widetilde{\mu}(S_t \times S_t^c) dt &=
\int_{0}^{\infty} \left(\sum_{x,y} \widetilde{\mu}(x,y)\mathbf{1}_{I_{x,y}} (t) \right) dt\\
&=\sum_{x,y}\widetilde{\mu}(x,y)\int_{0}^{\infty} \mathbf{1}_{I_{x,y}} (t) dt \\
&=\sum_{x,y} \widetilde{\mu}(x,y)|u(y) - u(x)|. \ \bsq
\end{align*}

\n
\textbf{Proof of  Cheeger's inequality (Theorem~\ref{Theo-Cheeger})}
Let $f\in \ell^2_0(X,m)$ with real values. 
In view of Lemma~\ref{Lem1-Cheeger}, we have to prove the following inequality:
$$
\dfrac{h(X)^2}{4} \Vert f\Vert^2\leq  \sum_{(x,y) \in X^2} |f(y) - f(x)|^2 \widetilde{\mu}(x,y). \eqno{(*)}
$$
Let $c\in \RRR$ be such that $m(\{ x\in X\ : \ (f+c)(x) >0\}) \leq 1/2.$
Since $\sum_{x\in X} f(x) m(x) =0,$ one checks that 
$\Vert f+c\Vert^2 \geq \Vert f\Vert^2.$ So, upon replacing $f$ by $f+c$, 
we can assume that $m(\{ x\in X\ : \ f(x) >0\})\leq 1/2$
(observe that the right
hand side of $(*)$ does not change when $f$ is replaced by $f+c$).

Let $f_+$ and $f_-$ be the positive and negative parts of $f,$ so that
$f= f_+- f_-$ and $ \Vert f\Vert^2= \Vert f_+\Vert^2 + \Vert f_-\Vert^2 $.
Set $u=f_+^2$ or $u=f_-^2.$ Writing $S_t=\{x\in X\ : \ u(x) >t\},$
observe that $m(S_t) \leq 1/2$ for all $t>0$ (and hence $m(S_t^c)\geq 1/2$).
Using first the area formula,  then
the definition of $h=h(X)$ and finally the co-area formula, we have
\begin{align*}
h\sum_{x\in X} u(x)m(x) &= h \int_{0}^{\infty}m(S_t) dt\\
&\leq  2\int_{0}^{\infty} \widetilde{\mu}(S_t \times S_t^c) dt =
\sum_{x,y}|u(y) - u(x)|  \widetilde{\mu}(x,y).
\end{align*}
So, 
$$
h\Vert f\Vert ^2\leq \sum_{x,y}|f_+^2(y) - f_+^2(x)|  \widetilde{\mu}(x,y) + \sum_{x,y}|f_-^2(y) - f_-^2(x)|  \widetilde{\mu}(x,y)
$$
By Cauchy-Schwarz inequality, we have
\begin{align*}
&\sum_{x,y}|f_{\pm}^2(y) - f_{\pm}^2(x)|  \widetilde{\mu}(x,y)\\
&\leq \left( \sum_{x,y}|f_{\pm}(y) - f_{\pm}(x)|^2 \widetilde{\mu}(x,y)\right)^{1/2}
\left(\sum_{x,y}|f_{\pm}(y) + f_{\pm}(x)|^2 \widetilde{\mu}(x,y)\right)^{1/2}.
\end{align*}
Now, using  the symmetry of $\widetilde \mu$, the fact that 
 $\sum_y \widetilde{\mu}(x,y) = m(x)$  and again  Cauchy-Schwarz inequality,
we have
\begin{align*}
\sum_{x,y}|f_{\pm}(y) + f_{\pm}(x)|^2 \widetilde{\mu}(x,y)&=\sum_{x,y}( f_{\pm}^2(y) + 2  f_{\pm}(x) f_{\pm}(y)+ f_{\pm}^2(x))  \widetilde{\mu}(x,y)\\
&= 2\sum_{x,y} f_{\pm}^2(x)  \widetilde{\mu}(x,y)  + 2 \sum_{x,y} f_{\pm}(x) f_{\pm}(y)\widetilde{\mu}(x,y)  \\
&\leq 4\Vert  f_{\pm}\Vert^2 
\end{align*}
Hence,
\begin{align*}
h\Vert f\Vert^2 &\leq 2\Vert  f_{+}\Vert \left(\sum_{x,y}|f_{+}(y) - f_{+}(x)|^2 \widetilde{\mu}(x,y)\right)^{1/2} +\\
&+2\Vert  f_{-}\Vert\left( \sum_{x,y}|f_{-}(y) - f_{-}(x)|^2 \widetilde{\mu}(x,y)\right)^{1/2}
\end{align*}
Now, it is straightforward to check that 
$$
\sum_{x,y}|f_{\pm}(y) - f_{\pm}(x)|^2 \widetilde{\mu}(x,y) \leq \sum_{x,y}|f(y) - f(x)|^2 \widetilde{\mu}(x,y). 
$$
Hence,
\begin{align*}
h^2\Vert f\Vert ^4&\leq 4\Vert  f\Vert^2  \sum_{x,y}|f(y) - f(x)|^2 \widetilde{\mu}(x,y)
\end{align*}
and inequality $(*)$ follows.$\bsq$

\begin{remark}
\label{Rem-Cheeger}
There is also an upper bound which is almost trivial, namely  
$$\lambda_1\leq 2 h(X).$$ 
Indeed,  it is straightforward
to check that, for  a non empty subset $S$ of $X,$ we have
$$
  \frac{\langle \Delta_{\mu} f  ,  f \rangle}{\Vert f\Vert^2} =\frac{\widetilde{\mu}(S \times S^c)}{m(S)m(S^c)}
 $$
 for the function
 $f=\mathbf{1}_S -m(S)\mathbf{1}_X$ in $\ell^2_0(X)$ and the inequality follows from Lemma~\ref{Lem1-Cheeger}.
\end{remark}

\subsubsection{Spectral Gap Property for $PGL_2(\FF_q((t^{-1})))/PGL_2(\FF_q[t])$}
\label{SS:SPG-PGL}
Let $\FF_q$ be the finite field with $q$ elements and $\FF_q((t^{-1}))$  the local field  of  Laurent series,
which is the the completion of $\FF_q(t)$ with respect to the valuation at infinity
defined by $v(a/b)=\deg (b)-\deg (a)$ for $a,b\in \FF_q[t], b\neq 0.$
The corresponding  compact subring 
is the ring of formal series $\FF_q[[t^{-1}]].$

Let $G= PGL_2(\FF_q((t^{-1}))).$ The subgroup  $K= PGL_2(\FF_q[[t^{-1}]])$
is compact and open in $G.$ 
As described in  II.1.1 and II.1.6  of \cite{Serre} (see also \cite{Efrat}), $T= G/K$ can be endowed 
with the structure of a $q+1$-regular tree: one can take as set of representatives for the cosets in $G/K$
the set of matrices
$$
\left(\begin{array}{cc}
t^{n}&\alpha\\
0&1\end{array}\right),
$$
with $n\in \ZZ$ and $\alpha$ from a set of representatives of  $\FF_q[[t^{-1}]]/(t^n);$ 
 the neighbours of the vertex $\left(\begin{array}{cc}
t^{n}&\alpha\\
0&1
\end{array}\right)$ are the $q+1$ vertices 
$$\left(\begin{array}{cc}
t^{n+1}&\alpha\\
0&1 
\end{array}\right), \ 
 \left(\begin{array}{cc}
t^{n-1}&\alpha +\beta t^n\\
0&1
\end{array}\right), \qquad \beta\in \FF_q.
$$
The group $G$ acts on $T$ by isometries 
on the left. 

Let $\Ga = PGL_2(\FF_q[t]),$ which is a discrete subgroup of $G.$
The quotient graph $X:=\Ga\backslash T\cong \Ga\backslash G/K$ is a half-line tree
$$
\underset{x_0}\circ \rule{1cm}{0,2mm}\underset{x_1}\circ\rule{1cm}{0,2mm}\underset{x_2}\circ\rule{1cm}{0,2mm}\underset{x_3}\circ\rule{0,5cm}{0,2mm}\dots$$
 given by (the
cosets of) the elements 
$$x_n=\left(\begin{array}{cc}
t^{n}&0\\
0&1
\end{array}\right)  \qquad n\geq 0.
$$
The  $q+1$ edges $(x_0,y) \in E(T)$ are  mapped to the edge
$(x_0, x_1)\in E(X)$; for $n\geq 1,$ the $q$-edges 
$(x_n, y)$ with $y=
 \left(\begin{array}{cc}
t^{n-1}&\beta t^n\\
0&1
\end{array}\right)$ for $\beta\in \FF_q$ are mapped to the edge $(x_n,x_{n-1})\in E(X)$
and the edge  $(x_n, y)$ with
$y=\left(\begin{array}{cc}
t^{n+1}&0\\
0&1 
\end{array}\right)$ is mapped to $(x_n,x_{n+1})\in E(X).$
Let $\la$ be the Haar measure on $G$ normalized by $\la(K)=1.$
We have  a decomposition of $\Ga\backslash G$ as disjoint union
$$\Ga\backslash G = \coprod_{n\geq 0} x_n K,$$
where $x_n K\subset \Ga\backslash G$ is the $K$-orbit  of $x_n$
in $\Ga\backslash G.$ So, $x_n K\cong ( x_n^{-1} \Ga_n x_n) \backslash K$,
where $\Ga_n= \Ga\cap x_n K x_n ^{-1}$ is the  stabilizer in $\Ga$ of $x_n$
for the action $\Ga\curvearrowright T$. One checks that 
$|\Ga_0|= q(q^2-1)$ and $|\Ga_n|= q^{n+1}(q-1)$ for $n\geq 1.$
Let $\overline{m}$ be the measure on the set $X$ defined by 
$$
\overline{m}(x_0)= \dfrac{1}{|\Ga_0|}=  \frac{1}{q(q^2-1)}, 
 \qquad \overline{m}(x_n)= \frac{1}{|\Ga_n|}= \dfrac{1}{q^{n+1}(q-1)} \quad\text{for}\qquad n\geq 1.
$$
Since $\sum_{n=1}^{\infty}  \dfrac{1}{q^{n+1}(q-1)} <\infty,$ 
$X$ has finite measure and hence we have
$$\vol(\Ga\backslash G) = \la\left(\coprod_{n\geq 0} x_n K \right)= \overline{m}(X)<\infty,$$
 showing that  $\Ga$ is a (non-uniform) lattice in $G$.

The simple random walk on $T,$ which is given by 
the  transition probabilities $$
\mu(x,y)=
\begin{cases} 
\dfrac{1}{q+1}& \text{if $(x,y)\in E(T)$}\\
0&\text{otherwise,}
\end{cases}
$$
 is reversible, with  stationary  measure $m: x \longmapsto 1.$
The associated projected random on $X$ is given 
by 
the  transition probabilities  $\overline{\mu}(x_0,x_1)=1$ and
$$
\overline{\mu}(x_n,x_{n+1})= \dfrac{1}{q+1}, \ \overline{\mu}(x_n,x_{n-1})=  \dfrac{q}{q+1} \ \  \text{for} \quad n\geq 1.
$$
As is easily  checked, one has
$$
\overline{m}(x_n)  \overline{\mu}(x_n,x_{n-1}) = \overline{m}(x_{n-1})  \overline{\mu}(x_{n-1},x_{n})
$$
for all $n\geq 1$, which means that 
$\overline m$ is a stationary measure for $\overline{\mu}$.
The   Markov  operator  $M_T$ on $\ell^2(T,m)$ for the random walk on $T$, which is defined by
$$
M_Tf(x)= 
\frac{1}{q+1}\sum_{(x,y)\in E(T)} f(y),\tout f\in \ell^2(T),
$$
 commutes with the $\Ga$-action;  it 
induces the Markov operator $M_X$  on $\ell^2(X, \overline m)$
corresponding to the projected random walk on $X$ and is
given by 
$$
M_Xf(x_n)= 
\begin{cases} 
  \dfrac{q}{q+1}f(x_{n-1}) +\dfrac{1}{q+1}f(x_{n+1}) & \text{for $n\geq 1$}\\
f(x_1)&\text{for $n=0$}
\end{cases}
\qquad f\in \ell^2(X,\overline m).
$$
The operator $M_X$ is self-adjoint; the corresponding Laplacian operator 
$\Delta_X= {\rm Id}_X -M_X$ is a non-negative operator  with spectrum contained 
 $[0,2].$
It is easy to show  that $G\curvearrowright \Ga\backslash G$ has the Spectral Gap
Property if and only if $1$ does not belong to the spectrum of the restriction of $M_X$ to $\ell^2_0(X),$
the orthogonal space to the constants (see Proposition~6 in \cite{BeLu}; here the compactness
of $K$ is crucial).
So, by Theorem~\ref{Theo-Cheeger}, it suffices to show that the Cheeger constant $h(X)$ of the random walk on $X$ is strictly positive.
This is indeed the case. 
\begin{proposition} 
\label{Pro-SPG}
We have $h(X)\geq \min\left\{\dfrac{q-1}{q+1}, \dfrac{4q^2}{(q+1)(q^2-1)}\right\} >0.$
Hence, the action $G\curvearrowright \Ga\backslash G$ has the Spectral Gap
Property.
\end{proposition}
\proof
Recall that $$h(X):= \inf_{S\subset X, S\neq \emptyset} \frac{\widetilde{\mu}(S \times S^c)}{\overline{m}(S)\overline{m}(S^c)},$$
where $\widetilde{\mu}(x,y)=\overline{m}(x) \overline{\mu}(x,y).$
In order to simplify the computations, we rescale $\overline{m}$ 
 by
$$\overline{m}(x_0)= \frac{1}{q+1},  \qquad \overline{m}(x_n)=  \dfrac{1}{q^{n}} \tout n\geq 1.$$
Let $S$ be  a non empty subset of $X;$  replacing $S$ by $S^c$ if necessary, we can assume
that $\overline{m}(S)\leq m(X)/2.$

One checks that $\overline{m}(X)= \dfrac{2q}{q^2-1}$ and that
$$\overline{m}(x_0)+\overline{m}(x_1)=  \dfrac{2q+1}{q(q+1)}> \overline{m}(X)/2,$$
 since $q\geq 2.$
So, two cases may occur.

\noindent
$\bullet$\emph{First case:} $x_0\in S.$ 
Then $x_1\notin S$ and hence $(x_0,x_1)\in S\times S^c$.
Therefore,
\begin{align*}
\frac{\widetilde{\mu}(S \times S^c)}{\overline{m}(S)\overline{m}(S^c)} \geq 2 \dfrac{\widetilde{\mu}(x_0,x_1)}{\overline{m}(X)^2}= \dfrac{2}{(q+1)\overline{m}(X)^2}.
\end{align*}
\noindent
$\bullet$\emph{Second case:} $x_0\notin S.$  Let $n\geq 0$ be minimal with the
property that $x_{n+1}\in S.$ Then $(x_{n+1},x_{n})\in S\times S^c.$
Hence, 
$$\overline{m}(S) \leq \sum_{k\geq n+1}\overline{m}(x_k)= \dfrac{1}{q^n(q-1)}$$ and so
\begin{align*}
\frac{\widetilde{\mu}(S \times S^c)}{\overline{m}(S)\overline{m}(S^c)} &\geq \dfrac{\widetilde{\mu}(x_{n+1}, x_{n})}{\overline{m}(X)\sum_{k\geq n+1}\overline{m}(x_k)}=\dfrac{{1}/{q^n(q+1)}}{{\overline{m}(X)}/{q^n(q-1)}}\\
& = \dfrac{q-1}{(q+1) \overline{m}(X)}.
\end{align*}
Normalizing $\overline{m}$ to a probability measure, we obtain
$$h(X)\geq \min\{\dfrac{q-1}{q+1}, \dfrac{2}{(q+1) \overline{m}(X)}\}=\min\{\dfrac{q-1}{q+1}, \dfrac{4q^2}{(q+1)(q^2-1)}\} .\bsq$$

\begin{remark}
\label{Rem-Efrat}
The precise spectral decomposition of $M_X$ (or $\Delta_X$) acting on $\ell^2(X)$ is determined  in \cite{Efrat}.
In particular, it is shown there that the spectrum of $M_X$ is $\left[{-2\sqrt{q}}/{(q+1)},{2\sqrt{q}}/{(q+1)} \right]\cup\{\pm 1\}.$
(Observe that $-1$ is indeed an eigenvalue of $M_X$ with eigenfunction $f: X\to \RRR$ defined by
$f(x_n)= (-1)^n$ for all $n\geq 0;$ this is related to the fact that $G=PGL_2(\FF_p((t^{-1}))$ has a 
a one dimensional character $\neq 1_G$, as its abelianization $PGL_2(\FF_p((t^{-1})))/PSL_2(\FF_p((t^{-1})))$
is a group of order $2$.)

\end{remark}
\subsection{Lattices without the Spectral Gap Property}
\label{SS:LatticesWithout}
 In view of Theorem~\ref{Theo-LieGroups}, one might think that $G\curvearrowright  G/\Ga$
 has the Spectral Gap Property for any   locally compact group $G$ and any lattice $\Ga$  in $G.$
 This is however not the case, as the following result from  \cite{BeLu} shows.
\begin{theorem}
\label{Theo-NonSGP} 
\textbf{(\cite{BeLu})}
For an integer $k>2,$ let $T_k$ be the $k$--regular tree and  $G={\rm Aut} (T_k).$
Then $G$ contains a lattice $\Ga$ such that $G\curvearrowright  G/\Ga$ does not have the Spectral Gap Property.
\end{theorem}
The proof here is based on  the  (easy) inequality $\lambda_1 \leq 2 h(X)$ between 
the Cheeger constant and the bottom of the spectrum of the Laplacian $\lambda_1$ 
of an appropriate random walk which is determined as follows.
We can find a reversible random walk on a countable graph $X$ with transition probability
 $\overline{\mu}: X\times X\to \RRR$ and stationary measure $\overline m$ with the following properties:
$\overline{m}(X)<\infty,$ the Cheeger constant $h(X)$ of $X$ is $0$ and, moreover,
 $k \overline{\mu}(x,y)$ and ${1}/{\overline{m}(x)}$ are  integers for all $x,y.$
By the ``inverse    Bass--Serre theory'' of groups acting on trees (see \cite{BaLu}),
there exists  a lattice $\Ga$ in  $G={\rm Aut} (T_k)$
such that the projection of the standard random walk  
on $\Ga\backslash T_k$ can be identified with $(X,\overline{\mu},\overline{m} ).$
For more details,  see \cite{BeLu}.

\section{Actions on tori and nilmanifolds with the Spectral Gap Property }
\label{S:Nil}
We review the main results from \cite{BeGui-AutoNil}, concerning actions
of a (countable) group $\Ga$ by automorphisms (or affine transformation)
on a compact nilmanifold.
First, we state the result for actions on the torus.

Let  $T= \RRR^d/\ZZ^d$ be the $d$-dimensional torus.
Observe that   $\Aut (T)$ can be identified with $GL_d(\ZZ)$.
Set  $V=\RRR^d$  and $\La=\ZZ^d.$ If $W$ is a
rational linear subspace of $V,$
then  $S= W/ (W\cap \La)$ is a subtorus  of $T$ 
and we have a torus factor  $\overline T= T/S.$ 
 Let $\Ga$ be a subgroup
of $\Aut (T)$ and assume that $W$ is $\Ga$-invariant.
Then  $\Ga$  leaves  $S$  invariant and 
the induced action of $\Ga$ on $\overline T$ is a  factor of 
the action of $\Ga$ on $T.$ 
We will say that $\overline T$ is a 
$\Ga$-invariant  torus factor of $T.$   Here is our main result
for actions by torus automorphisms.
\begin{theorem}
\label{Theo-SGP-AutoTorus}
Let $\Ga$ be a subgroup of $GL_d(\ZZ)$.
The following properties are equivalent.
\begin{itemize}
 \item [(i)] The action of $\Ga \curvearrowright T=\RRR^d/\ZZ^d$ has  the  Spectral Gap
 Property.
\item [(ii)] There exists no   non-trivial $\Ga$-invariant  torus factor $\overline T$
such that the projection of $\Ga$ on $\Aut (\overline T)$  is 
amenable.
 \item [(ii)]  There exists no  non-trivial  $\Ga$-invariant  torus  factor $\overline T$ of $T$    such that
the projection of  $\Ga$ on $\Aut(\overline T)$ is a virtually abelian group (that is,  it  contains an abelian subgroup of finite index).
\end{itemize}
\end{theorem}

The following corollary gives
a large class of examples of groups of automorphisms of the torus
with  the  Spectral Gap Property.

\begin{corollary}
\label{Cor-SGP-AutoTorus}
Let $\Ga$ be a subgroup of $GL_d(\ZZ)$.
Assume that is $\Ga$ is not virtually abelian and that 
$\Ga$ acts   $\QQ$-irreducibly  on $\RRR^d$
(that is, there are no non trivial $\QQ$-rational 
subspaces which are invariant under $\Ga$).
Then  $\Ga \curvearrowright T=\RRR^d/\ZZ^d$ has  the  Spectral Gap
 Property.
 \end{corollary}
 
\begin{remark}
\label{Rem-SGP-AutoTorus}
The previous corollary was obtained in \cite{FuSh}
under the stronger assumption that $\Ga$ acts   $\RRR$-irreducibly  on $\RRR^d.$
\end{remark}

\begin{example}
\label{Exa-Torus}
Let $\Ga$ be a subgroup of $GL_d(\ZZ)=\Aut(T)$.
We identity the dual group of $T=\RRR^d/\ZZ^d$ 
with $\ZZ^d$ in the usual way.
As in the proof of Proposition~\ref{Pro-Erg-Mixing},
the Fourier transform sets up a unitary equivalence between
the Koopmann representation $\pi_T$ on $L^2_0(T)$ 
  and the  natural representation   of $\Ga$ on
$\ell^2\left(\ZZ^d\setminus\{0\}\right)$ 
defined by the  dual action of $\Ga$ on $\ZZ^d,$
which is given by $(\ga,x)\mapsto (\ga^t)^{-1}x.$
So, $\Ga\curvearrowright T$ has the  Spectral Gap
 Property if and only if  the action $\Ga\curvearrowright \ZZ^d\setminus\{0\}$
 is not co-amenable.
 
 Choose a set of representatives $S$ for the $\Ga$-orbits in
$\ZZ^d\setminus\{0\};$ we have 
a  direct sum decomposition 
$$\ell^2\left(\ZZ^d\setminus\{0\}\right)=
\bigoplus_{s\in S}  \ell^2(\Ga/\Ga_s),
$$
into $\Ga$-invariant subspaces, where $\Ga_s$ is the stabilizer of $s.$ 
Therefore,  
$$\Vert\pi_T(\mu)\Vert=\sup_{s\in S}\Vert\pi_{\Ga/\Ga_s}(\mu)\Vert $$
for every  probability measure $\mu$ on $\Ga.$

Let $d=2.$ Then every subgroup $\Ga_s$ is amenable, as it is   conjugate inside $GL_2(\RRR)$ to a subgroup 
of $N=\left\{
\left(\begin{array}{cc}
1&x\\
0&1
\end{array}\right)
 \ : \ x\in \RRR
\right\}.$ Hence, $\Vert\pi_{\Ga/\Ga_s}(\mu)\Vert =\Vert\pi_{\Ga}(\mu)\Vert $
and so $\Vert\pi_{T}(\mu)\Vert =\Vert\pi_{\Ga}(\mu)\Vert.$
In particular, $\Ga\curvearrowright T$ has the  Spectral Gap
 Property if and only if $\Ga$ is not amenable.

\end{example}

Actions by automorphisms on nilmanifolds are a natural generalization
of actions by torus automorphisms. The setting is as follows. 

Let $N$ be a connected and simply connected nilpotent Lie group.
Let $\La$ be a lattice in $N;$ the associated nilmanifold
 $N/\La$ is known to be compact.
 Observe that not every  nilpotent Lie group has a lattice: a necessary and sufficient condition  for
 this to happen is that  $N$ is an algebraic group defined over $\QQ$ (see \cite{Raghunathan}).
 
Let $\Aut(N)$  be the group  of  continuous automorphisms of $N$
and denote by  $\Aut (N/\La)$ the subgroup
of continuous automorphisms  $g$ of $N$ such that $g(\La) =\La.$
Every $g\in \Aut (N/\La)$  
preserves the translation invariant probability
measure $m$ on $N/\La$ induced by a Haar measure
on $N.$  
The nilsystem   $(N/\La, m)$ has a  natural maximal  torus
factor $(T, m');$  every automorphism $g\in \Aut (N/\La)$ induces a torus
automorphism  $\overline{g}\in \Aut(T)$ and the mapping
$g\mapsto \overline g$ is a homomorphism $ \Aut (N/\La)\to \Aut(T)$;
see  \cite{Parry} for this, as well
as for other results on the ergodic theory of automorphisms of nilmanifolds.

The following theorem reduces the question of the  Spectral Gap Property 
for groups of automorphism of nilmanifolds to  the   same question for  groups acting on tori.

\begin{theorem}
\label{Theo-SGP-AutoNil}
Let $N/\La$ be a compact  nilmanifold,  
with associated maximal torus factor $T.$ 
Let  $\Ga$ be a  subgroup of $\Aut(N/\La)$.
 The following properties are equivalent.
\begin{itemize}
 \item [(i)] The action of $\Ga \curvearrowright N/\La$  has the  Spectral Gap
 Property.
 \item[(ii)] The action of  $\Ga \curvearrowright T$  has the  Spectral Gap
 Property.
 \end{itemize}
\end{theorem}

\begin{example}
\label{Exa-Heisenberg}
Let $N={H_{3}(\RRR)}$ be the
$3$--dimensional real Heisenberg group.
Recall that $N$ can be realized as the group 
with underlying set $\RRR^3=\RRR^{2}\times \RRR$ and product
$$((x,y),s)((x',y'),t)=\left((x+x',y+y'),s + t +  xy'-x'y))\right).$$
The Lebesgue measure $m$ on $\RRR^3$  is 
a (left and right) Haar measure on $N.$
The group $N$ is  a two-step nilpotent Lie group;  
its  centre $Z$ coincides with its commutator subgroup 
and is given by $Z=\{((0,0),s)\ :\ s\in \RRR\}.$
The  group $SL_{2}(\RRR)$
acts by automorphisms on $N,$ via
$$
 g((x,y),t)= (g(x,y),t) \tout g\in SL_{2}(\RRR),\ (x,y)\in \RRR^{2},\ t\in \RRR.
$$
The discrete subgroup 
$
\La= \{((x,y),s)\ :\ x,y\in \ZZ^2, s\in \ZZ\}
$
is a  lattice in $N$. The group $SL_2(\ZZ)\subset \Aut(N)$ preserves
$\La$ and acts therefore on 
the Heisenberg nilmanifold $X:=N/\La.$

The maximal torus factor $T$ can be identified with $\RRR^2/\ZZ^2,$ via the 
$SL_2(\ZZ)$-equivariant mapping
$X\to \RRR^2/\ZZ^2, [((x,y),s)]\mapsto [(x,y)].$
Identifying $L^2(T)$  with  a closed $\pi_X(\Ga)$-invariant subspace of $L^2(X)$,
we  have a decomposition 
$$L^2(X)=L^2(T)\oplus \H$$
into $\pi_X(\Ga)$-invariant subspaces,
 where $\H$ the orthogonal complement of
 $L^2(T)$ in $L^2(X).$
 
Let $\mu$ be a symmetric probability measure  on $SL_2(\ZZ)$
and $\Ga= \Ga(\mu)$ the subgroup  generated 
by the support of $\mu.$ 
By Example~\ref{Exa-Torus},  the restriction of 
$\pi_X(\mu)$ to $L_0^2(T)$ has norm 
$$ \Vert \pi_X(\mu)\Vert_{L_0^2(T)}=\Vert \pi_{\Ga}(\mu)\Vert.$$
On the other hand, it can be shown that the restriction of 
$\pi_X(\mu)$ to $\H$ has norm
$$
\Vert \pi_X(\mu)\Vert_{\H}\leq \Vert \pi_{\Ga}(\mu)\Vert^{1/4}.
$$
The proof of this inequality
involves the consideration of the so-called Weil representation 
of (a two-fold cover of)  the simple  Lie group $SL_2(\RRR)$ and 
estimates of its matrix coefficients, much in the spirit of Section~\ref{S:QuantitativeSG}.
Summarizing, we have
$$
\Vert \pi_X(\mu)\Vert \leq \Vert \pi_{\Ga}(\mu)\Vert^{1/4}.
$$
In particular,  
$\Ga \curvearrowright X$  has the Spectral Gap Property if and only if $\Ga$
is non-amenable. 
For more details on this example, see \cite{BeHe}.
\end{example}

 \section{Some applications of  the Spectral Gap Property}
 \label{S:Appl}
 We give two  applications of the Spectral Gap Property of group actions,
 one to the construction of expander graphs and the other 
 to the escape rate of random walks on linear groups.
\subsection{Expander graphs}
 \label{Appl-Expanders}
Let ${\G}=(X,\EE)$ be a   finite $k$-regular graph, where $X$
is the set of vertices and $\EE\subset X\times X$ the set of edges 
of ${\G}$. 
We consider the simple random walk on $X$
defined by the transition probabilities
$p:X\times X \to \RRR$ given by $p(x,y)=\frac{1}{k}$
if $(x,y)\in \EE$ and $p(x,y)=0$ otherwise.
The map  $m:X\to \RRR, x\mapsto 1/|X|$ is
a stationary measure for $\mu.$  
The Cheeger constant (also known as expanding constant) of $X$ or $\G$, as defined in Section~\ref{SS-Cheeger}, is 
\begin{align*}
h(X)&=\dfrac{|X|}{k}\min\left\{\frac{|\partial S|}{|S||S|^c}\ :  S\subset X, S\neq \emptyset\right\},
\end{align*}
where the boundary $\partial S$ of a subset $S$ of $X$  is the set of edges
$(x,y)$ with $x\in S$ and $y\notin S.$
A more commonly used constant is the so-called expanding or isoperimetric constant of the graph,
defined by 
$$
\widetilde{h}(X)=\min\left\{\frac{|\partial S|}{|S|}\ :  0<|S|<|X|/2\right\};
$$
 $\widetilde{h}(X)$ and $h(X)$ are related by ${\widetilde{h}(X)}/k\leq   h(X)\leq {2\widetilde{h}(X)}/{k}.$

Expander graphs are families of graphs which are both
 sparse  and strongly  connected. More precisely, 
a  sequence   
 of finite $k$-regular graphs ${\G}_n=(X_n,\EE_n)$
 with $\lim_{n\to\infty}|X_n|=\infty$ 
is  a family of \emph{expanders} if
$\inf_{n}h(X_n)>0.$
A constant $\eps>0$ with $\inf_{n}h(X_n)\geq \eps$ is called an \emph{expanding constant} for 
the sequence of expanders $({\G}_n)_n.$

Let ${\G}_n=(X_n,\EE_n)$ be a family of finite $k$-regular graphs 
 with $\lim_{n\to\infty}|X_n|=\infty.$ 
 Let $\Delta_n$ be the corresponding Laplacian
on $X_n$; recall that $\Delta_n$  is defined on  $\ell^2(X_n)$ by
$$\Delta_n f(x)= f(x)-\dfrac{1}{k}\sum_{(x,y)\in \EE_n} f(y).$$
 Let  $\lambda_1^{(n)}$ denote
the smallest eigenvalue $\neq 0$
of  $\Delta_n.$
 In view Cheeger's inequalities (Theorem~\ref{Theo-Cheeger} and Remark~\ref{Rem-Cheeger}),
we see that $(X_n)_n$ is a family of expanders if and only 
if $\inf_{n}\lambda_1^{(n)}>0.$

The existence of  {expander graphs}  is settled
by elementary counting arguments (see \cite[Proposition 1.2.1]{Lubotzky}).
However, their constructions seem to require sophisticated
mathematical tools. 
We will  give the   explicit construction of a family of expander graphs
using Kazhdan's Property (T), following the original  idea of  Margulis;
recently, families of expanders have been found using the so-called
zig-zag construction (see \cite{zigzag}, \cite{AlLuW--01}).

Let $G$ be  a finitely generated group, with a fixed
 finite generating set $S$ with $S^{-1}=S.$
The Cayley graph $\G(G,S)$  is a connected $k$-regular graph for $k=|S|$
(see Section~\ref{ApplicationsSGP-RW}).

We assume now that $G$ is a \emph{finite} group.  
Let $\pi_G$ be the right regular representation of $G$ on
$\ell^2(G).$ Denote by $\pi_G^0$  the corresponding representation of
$G$ on the $G$-invariant subspace $\ell_0^2(G)=\{{\Un}_G\}^{\perp}.$  

Let $\mu$ be the probability measure $\mu=\frac{1}{|S|} \sum_{s\in S} \delta_s$ on $G.$
The following crucial lemma establishes a   link between the norm of 
$\pi_G^0(\mu)$ and the
smallest eigenvalue $\lambda_1$ in $\ell_0^2(G)$
of the Laplace operator $\Delta$ associated
to the simple random walk on $\G(G, S).$
\begin{lemma}
\label{Lem-ExpConstT} 
We have  $ \lambda_1\geq   \dfrac{1}{2}(1-\Vert \pi_G^0(\mu)\Vert)^2.$
\end{lemma}
\proof
Let $f\in \ell^2_0(G)$ be an eigenfunction
of $\Delta$ with $\Vert f\Vert =1,$ for the eigenvalue $\lambda_1.$
Denoting by   $\EE$  the set  of edges of $\G(G, S),$ 
and using Lemma~\ref{Lem1-Cheeger} as well as Jensen's inequality, we have
\begin{align*}
2\lambda_1&=2{\langle \Delta f  ,  f \rangle} =\dfrac{1}{k} \sum_{(x,y) \in \EE} \vert{f(y) - f(x)}\vert^2
=\dfrac{1}{k}\sum_{s\in S} \sum_{x\in G} \vert{f(xs) - f(x)}\vert^2\\
&=\dfrac{1}{k}\sum_{s\in S}  \Vert{\pi_G^0(s)f - f}\Vert^2
 \geq \left(\dfrac{1}{k}\sum_{s\in S}  \Vert{\pi_G^0(s)f - f}\Vert\right)^2 \\
 &\geq\left\Vert \dfrac{1}{k}\sum_{s\in S} (\pi_G^0(s)f - f)\right\Vert^2 =\Vert\pi_G^0(\mu)f - f \Vert^2\geq (1-\Vert \pi_G^0(\mu)\Vert)^2.\bsq
\end{align*}

As a consequence of the previous lemma, we obtain a construction
scheme for expanders.
\begin{theorem}
\label{Theo-ExpGrT}
Let $\Ga$ be a group with Kazhdan's Property (T),
and $S$ a finite generating  set  of $\Ga$ with $S^{-1}=S.$
 Let $(N_n)_n$ be a sequence of normal subgroups of $\Ga$ of
finite index with $\lim_n |\Ga/N_n|=\infty.$
Then  the sequence  of the Cayley  graphs $\G(\Ga/N_n, \vfi_n(S))$
is an expander  family, where 
$\vfi_n(S)$ is the image of $S$ under the canonical projection 
$\vfi_n: \Ga\to\Ga/N_n.$
\end{theorem}
\proof
Let $\mu=\frac{1}{|S|} \sum_{s\in S} \delta_s.$
By Proposition~\ref{Theo: SPG-T},
that there exists a constant $\delta>0$ such that 
$1-\Vert \pi(\mu)\Vert \geq \delta$  for every
  unitary representation $\pi$ of $\Ga$ without non-zero invariant vectors.
 This holds in particular for the representations $ \pi_{\Ga/N_n}^0 \circ \vfi_n$
 acting on $\ell^2_0(\Ga/N_n).$
 The result follows now from the previous lemma.$\bsq$
 
 \begin{remark}
 With the notation as in the previous theorem, 
 let $\lambda_1^{(n)}$ denote the first non-zero eigenvalue
 of the Laplacian on $X_n=\Ga/N_n.$
 Since $h(\Ga/N_n) \geq \dfrac{\lambda_1^{(n)}}{2},$ 
 we see that  an expanding constant for the family  
 $\G(\Ga/N_n, p_n(S))$ is 
 $\eps= {\delta^2}/{2},$ where
 $\delta =\inf_{\pi} (1-\Vert \pi(\mu)\Vert)$  with $\pi$ running over all
  unitary representations of $\Ga$ without non-zero invariant vectors.
 \end{remark}
\begin{example}
\label{Exa1-GrExp}
Let $\Ga=SL_3(\ZZ).$ 
Then $S=\{ E_{ij}^{\pm 1}\ :\ 1\leq i, j\leq 3,\ i\neq j \}$
is a generating  set of $\Ga,$ where 
the $E_{ij}$'s  are the usual elementary matrices.
For every prime integer $p,$ let 
$\Ga(p)$ be the so-called
principal congruence subgroup, that is,
$$\Ga(p)=\{A\in \Ga\ :\ A\equiv I \mod p\}$$ 
is the kernel of the 
surjective homomorphism
$\vfi_p:SL_3(\ZZ)\to SL_3(\ZZ/p\ZZ)$ given by 
reduction modulo $p.$ 
Since $\Ga/\Ga(p)\cong SL_3(\ZZ/p\ZZ),$
the subgroup $\Ga(p)$ has finite index
$p^3(p^3-1)(p^2-1)\approx p^8.$
The family of Cayley graphs
$(\G(\Ga/\Ga(p), \vfi_p(S))_p$ is a family of 
$k$-regular expanders with 
$k=12.$
It can be shown that $\eps\approx \dfrac{10^{-6}}{4}$ is an expanding constant for this family
(see Example 6.1.11 in \cite{BHV}).
\end{example}

For a comprehensive account on expander graphs and their applications, see \cite{HLW}.
An overview of recent developments in this subject is given in \cite{Breuillard}.

\subsection{Growth of products of random matrices}
We now give an application of the Spectral Gap Property 
to a result of Furstenberg from \cite{Furst-NRP} about random walks on linear groups.
The setting is as follows.

Let $\mu$ be a probability measure on the special linear group $G=SL_d(\RRR).$
We set $V=\RRR^d.$  We will consider the 
  operator norm on ${\mathrm End}(V)= M_d(\RRR)$ 
 associated to the   Euclidean norm   $\Vert \cdot\Vert$ on $V.$

Let $S_n(\omega)= X_n(\omega) \cdots X_1(\omega)$
be a sequence of random products, where $(X_n)_{n\geq 1}$ is a sequence
of independent random variables defined on a common probability
space $(\Omega, \PP),$  with values in $G$ and identically distributed
according to $\mu.$
One is interested in a  non-commutative analogue of the \emph{Law of Large Numbers}
describing the top Lyapunov exponent $\lambda_1(\mu)$
of the random matrix products which is defined as follows. 

Assume that  $\mu$ has finite first moment, that is,$
\int_G \log\Vert g\Vert d\mu(g) <\infty.
$ 
It follows from Kingman's subadditive ergodic theorem (see Theorem 10.1 in \cite{Walters}) that 
$$\lim_{n\to \infty} \dfrac{1}{n}\log\Vert S_n(\omega)\Vert = \lim_{n\to \infty}  \dfrac{1}{n} \log\Vert X_n(\omega) \cdots X_1(\omega)\Vert$$
 exists $\PP$-almost surely and
is $\PP$-almost everywhere constant;  we denote by  $\lambda_1(\mu)$ this limit, which may 
also be computed as 
$$
\lambda_1(\mu) = \lim_{n\to \infty} \dfrac{1}{n}\int_G \log\Vert g\Vert d\mu^n(g) = \inf_{n\geq 1}  \dfrac{1}{n}\int_G \log\Vert g\Vert d\mu^n(g).
$$
For more details on random matrix products,  see the survey \cite{Furman-RW}.

We give in the following proposition a lower bound for $\lambda_1(\mu)$
in terms of  a Spectral Gap Property. For this, we follow \cite[Corollary 2.2]{GuivSpec} (see also 
the proof of Corollaire~1 in \cite{GuivSpec} and  Theorem 1.19 in \cite{Furman-RW} ).

Let $\pi_V$ be the  unitary representation of $G= SL_d(\RRR)$ on
$L^2(V,m)$ given by the natural action of 
$G$ on $V$ equipped with the Lebesgue measure $m.$

\begin{proposition}
\label{Prop-RMP}
The following inequality holds:
 $$
 \lambda_1(\mu) \geq \dfrac{1}{d} \log\left( \dfrac{1}{r_\spec(\pi_V(\mu))}\right).
 $$
 \end{proposition}
 
 \proof
 For $n\geq 1,$ set $u_n = \int_G \log\Vert g\Vert d\mu^n(g).$
 So,  $\lambda_1(\mu)= \lim_{n} {u_n}/{n}.$
 
Since  $1 = \det g \leq  \Vert g\Vert^d,$
we have  $\Vert g\Vert \geq  1$  for every $g\in G.$
Fix  $\eps > 0.$ 
Let $f_\eps \in L^2(V)$ be defined by 
$f_\eps(x) = 1$ if $\Vert x\Vert \leq 1$ and $f_\eps(x) = \dfrac{1}{\Vert x\Vert^{d+\eps}}$ if $\Vert x\Vert \geq 1.$ 
For $A= \{ x\in V \ :\ 1\leq \Vert x\Vert \leq 2\}$ and $n\geq 1,$
 we have ${\mathbf 1}_A \in L^2(V)$ and
\begin{align*}
\langle  \pi_V(\mu^n) {\mathbf 1}_A, f_\eps \rangle 
&=\langle \pi_V(\check{\mu}^n)  f_\eps ,   {\mathbf 1}_A\rangle \\
&= \int_A \int_G  \dfrac{1}{\Vert g x\Vert^{d+\eps}} dm(x) d\mu^n(g)\\
& \geq  \int_A \int_G  \dfrac{1}{(\Vert g\Vert  \Vert x\Vert)^{d+\eps}} dm(x) d\mu^n(g)\\
& \geq  m(A) \dfrac{1}{2^{d+\eps}} \int_G  \dfrac{1}{\Vert g \Vert^{d+\eps}} dm(x) d\mu^n(g).
\end{align*}
Hence,  by concavity of the logarithm, we obtain
\begin{align*}
\log (\langle \pi_V(\mu^n) {\mathbf 1}_A, f_\eps \rangle) 
& \geq  \log\left(\dfrac{m(A)}{2^{d+\eps}}\right) -(d+\eps) \int_G  \log\Vert g \Vert dm(x) d\mu^n(g),\\
\end{align*}
that is,
$$(d+\eps)u_n \geq    -\log( \langle  \pi_V(\mu^n) {\mathbf 1}_A, f_\eps) +\log\left(\dfrac{m(A)}{2^{d+\eps}}\right).$$
Since
\begin{align*}
\limsup_n |\langle  \pi_V(\mu^n) {\mathbf 1}_A, f_\eps \rangle|^{1/n} 
&\leq  \limsup_n  \Vert \pi_V(\mu^n) \Vert^{1/n}(\Vert  {\mathbf 1}_A\Vert \vert f_\eps \Vert)^{1/n} =r_\spec(\pi_V(\mu)),
\end{align*}
we have therefore
$$
(d+\eps) \lim_n\dfrac{u_n}{n} \geq -\log r_\spec(\pi_V(\mu)).
$$
Letting $\eps\to 0,$ we obtain the claim.$\bsq$

The following result is an immediate consequence
of the previous proposition in combination with Theorem~\ref{Theo-Guiv}.
\begin{corollary}
\label{Cor-RMP}
Let $\mu$ be probability measure on $SL_d(\RRR)$ and
denote by $\Ga(\mu)$  the subgroup  generated by the support  of $\mu.$
Assume that the linear action of $\Ga(\mu)$ on $V=\RRR^d$ is not co-amenable.
Then $ \lambda_1(\mu) >0.$
 \end{corollary}
 So, under the assumption of the previous
 corollary, the norm $\Vert S_n(\omega)\Vert$ grows exponentially 
 almost surely, where   $S_n(\omega)=X_n (\omega)\cdots X_1(\omega)$ is the product
 of independent random unimodular  matrices  identically distributed
according to $\mu.$

Applying  Theorem~\ref{Theo-FurstenbergLemma}, we recover the following 
result of Furstenberg (Theorem 8.6 in \cite{Furst-NRP}).
 \begin{corollary}
 \label{Cor-NRP2}
 \textbf{(\cite{Furst-NRP})}
 Assume that $\Ga(\mu)$ is not  bounded and that
the action of $\Ga(\mu)$ on $V$ is totally irreducible.
Then $ \lambda_1(\mu) >0.$
\end{corollary}

\begin{example}
Let $\mu$ be the  probability measure on $G=SL_2(\RRR)$ 
given by $\mu= \dfrac{1}{4} (\delta_a + \delta_b+\delta_{a^{-1}} +\delta_{b^{-1}}  )$ for the
matrices
$$
a=\left( \begin{array}{cc}
1&2\\
0&1
\end{array}
\right) \qquad \text{and}\qquad
b=\left( \begin{array}{cc}
1&2\\
0&1
\end{array}
\right).
$$
We claim that
$$
\lambda_1(\mu) \geq  \dfrac{1}{{2}} \log\left(\sqrt{\dfrac{2}{\sqrt{3}}}\right)\approx 0.015617
$$
for the corresponding  top Lyapunov exponent $\lambda_1(\mu).$

Indeed, the subgroup $\Ga$ generated by the support of $\mu$ is 
the free group on $a$ and $b$ and is a subgroup (of index 12) of $SL_2(\ZZ).$
Moreover,  the representation $\pi_V$ of $G=SL_2(\RRR)$ on $L^2(V),$ for 
$V=\RRR^2,$
is weakly contained in the regular representation $\pi_G;$ indeed, 
$G$ acts transitively on $V\setminus\{0\}$ with stabilizers conjugated to 
$N= \left\{ \left( \begin{array}{cc}
1&*\\
0&1
\end{array}
\right) \right \}.$
It follows that $\pi_V$ is equivalent to the quasi representation $\pi_{G/N}\cong \ind_N^G(1_N)$ 
on $L^2(G/N).$ Since $N$ is amenable, $\pi_V$ is weakly contained in $\pi_G$
(see Theorem F.3.5 in \cite{BHV}).
So, $\pi_V$ is strongly $L^{p}$ for every $p>2$ and hence (see Proposition~\ref{Pro-TemperedRep}) 
$\pi_V \otimes \overline{\pi_V}$ is contained in a multiple of $\pi_G$.
Therefore,  by Proposition~\ref{Pro-Nevo},
$$\Vert \pi_V(\mu)\Vert \leq \Vert \pi_G(\mu)\Vert^{1/2}.$$
Now, since $\Ga$ is a discrete subgroup of $G,$  
the restriction of $\pi_G$ to $\Ga$ is a multiple of   the regular representation $\pi_\Ga$
 and so 
$$\Vert \pi_V(\mu)\Vert \leq \Vert \pi_\Ga(\mu) \Vert^{1/2}.$$
 Finally, as $\Ga$ is a free group on two generators, we have 
 $\Vert \pi_\Ga(\mu) \Vert= {\sqrt{3}}/{2}$ (see Section~\ref{SS:KestenFree})
and the claim follows from Proposition~\ref{Prop-RMP}.

\end{example}
\section*{Acknowledgments}
Above all, I thank  Yves Guivarc'h for countless discussions concerning this survey;
thanks are also due to Pierre de la Harpe for helpful comments and suggestions.
I would like to express my gratitude to both of them as well
as to  my other co-authors Yves de Cornulier, Jean-Romain Heu, 
Alex Lubotzky, and Alain Valette, for  joint work on which a substantial part of  this paper is based. 
I would also like  to thank  Athanase Papadopoulos  for the invitation to write this survey.
Finally, I am grateful to Shin Nayatani for his invitation to  the Rigidity School in Tokyo
in January 2013, during which part of this work was done.

\end{document}